\documentclass[12pt]{article}
\usepackage{amscd,amsfonts,amssymb,amscd,amsmath,latexsym}
\usepackage[pdftex,backref]{hyperref}
\makeatletter
\renewcommand{\subsubsection}{\@startsection
{subsection}
{1}
{0mm}
{0mm}
{0mm}
{\normalfont\normalsize\itshape}}
\makeatother
  
\usepackage[all]{xy}
\usepackage{mathrsfs}
\usepackage[final]{showkeys}
\newtheorem{theorem}{Theorem}[section] 
\newtheorem{prop}[theorem]{Proposition}
\newtheorem{lem}[theorem]{Lemma}
\newtheorem{ddd}[theorem]{Definition}
\newtheorem{kor}[theorem]{Corollary}

\newcommand{\bR}{\mathbf{R}}

\newcommand{\Map}{{\tt Map}}

\newcommand{\bp}{\mathbf{p}}

\newcommand{\K}{\mathbb{K}}

\newcommand{\tensor}{\otimes}
\newcommand{\into}{\hookrightarrow}
\newcommand{\onto}{\twoheadrightarrow}
\newcommand{\iso}{\cong}

\newcommand{\forget}[1]{}

\numberwithin{equation}{section}
\newcommand{\Z}{\mathbb{Z}}

\newcommand{\proof}{{\it Proof.$\:\:\:\:$}}

\newcommand{\R}{\mathbb{R}}

\newcommand{\bF}{{\bf F}}

\newcommand{\cH}{\mathcal{H}}

\newcommand{\cV}{\mathcal{V}}

\newcommand{\cW}{\mathcal{W}}

\newcommand{\cK}{\mathcal{K}}
\newcommand{\cA}{\mathcal{A}}

\newcommand{\Hom}{{\tt Hom}}

\newcommand{\Ext}{{\tt Ext}}

\newcommand{\im}{{\tt im}}

\newcommand{\cF}{\mathcal{F}}

\newcommand{\bc}{\mathbf{c}}

\newcommand{\id}{{\tt id}}

\newcommand{\nat}{\mathbb{N}}

\def\imath{{i}}

\def\hB{\hspace*{\fill}$\Box$ \newline\noindent}

\def\hB{\hspace*{\fill}$\Box$ \\[0cm]\noindent}

\newcommand{\Mat}{{\tt Mat}}
\newcommand{\cP}{\mathcal{P}}

\newcommand{\bH}{\mathbf{H}}

\newcommand{\pr}{{\tt pr}}

\newcommand{\bz}{\mathbf{z}}

\newcommand{\ev}{{\tt ev}}

\newcommand{\bh}{{\mathbf{h}}}

\newcommand{\bu}{\mathbf{u}}

\title{The topology of $T$-duality for $T^n$-bundles}
\author{Ulrich Bunke\thanks{Mathematisches Institut,
    Georg-August-Universit\"at G\"ottingen, bunke@uni-math.gwdg.de},
  Phillip Rumpf\thanks{Fakult\"at f\"ur Mathematik, 
Universit\"at M\"unster, Einsteinstr.~62, 48149 M{\"u}nster, Germany,
p\_rump01@math.uni-muenster.de} and Thomas Schick
\thanks{Mathematisches Institut, Georg-August-Universit{\"a}t G{\"o}ttingen,
Bunsenstr. 3, 37073 G{\"o}ttingen, Germany, 
schick@uni-math.gwdg.de}
}

\begin{document}
 
%\begin{abstract}
%\end{abstract}

\maketitle 
\tableofcontents
\parskip3ex

\section{Introduction}

\subsubsection{}
String theory is a part of mathematical quantum physics. Its ultimate goal is
the construction of quantum theories
modeling the basic structures of our universe.
More specifically, a string theory should associate a quantum field theory to
a target consisting of a manifold equipped with further geometric structures
like metrics, complex structures, vector bundles with connections, etc. A schematic
picture is $$\xymatrix{\mbox{target}\ar[rr]^{\scriptsize
    \hspace{-1.2cm}\mbox{string theory}}&&\mbox{quantum field theory}}\ .$$
The target is thought of to encode fundamental properties of the universe.
Actually there are several types of string theories, the most important ones
for the present paper are called of type $IIA$ and $IIB$ (see \cite[Ch. 10]{pol1}).

\subsubsection{}

$T$-duality is a relation between two string theories on the level of quantum field theories
to the effect that two different targets can very well lead to the same quantum field theory.
 The simplest example is the duality of bosonic string theories on the circles
 of radius $R$ and $R^{-1}$ (see \cite[Ch. 8]{pol1}). A relevant problem
 is to understand the factorization of the $T$-duality given on the level of 
 quantum theory  
 through $T$-duality on the
 level of targets. Schematically it is the problem of understanding the dottet
 arrow in 
$$\xymatrix{\mbox{target}\ar@{.>}[d]^{\scriptsize\mbox{target level $T$-duality}}\ar[rrr]^{\scriptsize\hspace{-1cm}\mbox{string theory, e.g. IIA}}&&&\mbox{quantum field theory}\ar[d]^{\scriptsize\mbox{quantum level $T$-duality}}\\\mbox{target}\ar[rrr]^{\scriptsize\hspace{-1cm}\mbox{string theory, e.g IIB}}&&&\mbox{quantum field theory}&}\ .$$
The problem starts with the question of existence, and even of the meaning
of such an arrow.
\subsubsection{}

$T$-duality on the target level is an intensively studied object in physics as well as in mathematics.
We are not qualified to review the extensive relevant literature here, but let
us mention mirror symmetry as one prominent aspect, mainly studied in
algebraic geometry (see e.g. \cite{SYZ}).

\subsubsection{}

In general, the target of a string theory is a manifold equipped with further geometric structures which in physics play the role of low-energy effective fields. The problem of topological $T$-duality 
can be understood schematically as the question of studying the dotted arrow in the following diagram.

$$\xymatrix{\mbox{target}\ar[d]^{\scriptsize\mbox{target level
      $T$-duality}}\ar[rr]^{\text{\hspace{-1.8cm}\scriptsize forget geometry}}&&\mbox{underlying topological space}\ar@{.>}[d]^{\scriptsize\mbox{topological $T$-duality}}\\
\mbox{target}\ar[rr]^{\scriptsize\mbox{\hspace{-1.8cm}forget geometry}}&&\mbox{underlying topological space}}\ .$$

\subsubsection{}

At this level one faces the following natural problems.
\begin{enumerate}
\item How can one characterize the topological $T$-dual of a topological
  space?  It is not a priori clear that this is possible at all.
\item If one understands the characterization of $T$-duals on the topological level, then one  wonders if a given space admits a $T$-dual.
\item Given a satisfactory characterization of topological $T$-duals one 
  asks for a classification of $T$-duals of a given space.
%\item Are their ways to improve $T$-duality to a functorial construction?
\end{enumerate}
As long as string theory is not part of rigorous mathematics the answer to the first question has to be found by physical reasoning and is part of the construction of mathematical models. Once an answer has been proposed the remainig two questions can be studied rigorously by methods of algebraic topology.

This is the philosophy of the present paper. For a certain class of spaces to
be explained below we propose a mathematical characterization of topological $T$-duals. On this basis we then present  a thourough and rigorous study of the existence and classification problems.

\subsubsection{}

The expression ``space'' has to be understood in a somewhat generalized sense since we consider 
targets with additional non-trivial $B$-field background. There are several possibilities to model
these backgrounds mathematically.  In the present paper we use an axiomatic
approach going under the notion of a twist, see \ref{twistsexplain}.
% Concrete realizations can be based on the notion of $U(1)$-banded gerbes \cite{math.GT/0508550}, bundles of compact operators (see \cite{mr}, \cite{mr1}) or maps to an Eilenberg-MacLane space $K(\Z,2)$, see \ref{twistsexplain}.

\subsubsection{}

Topological $T$-duality in the presence of non-trivial $B$-field backgrounds has been studied mainly
in the case of $T^n$-principal bundles (\cite{bhm}, \cite{bhm2}, \cite{bhm3}, \cite{bem1}, \cite{bem2}, \cite{mr}, \cite{mr1}, \cite{mr2}). Our proposal for the characterizations of $T$-duals in terms of $T$-duality triples is strongly based on the analysis made in these papers.

\subsubsection{}

The quantum field theory level $T$-duality predicts transformation rules for the low-energy effective fields which are objects of classical differential geometry like metrics and connections on the $T^n$-bundle, but also more exotic objects like a connective structure and a curving of the $B$-field background (these notions are explained in the framework of gerbes e.g. in \cite{hit}).
These transformation rules
are known as Buscher rules \cite{Bus}, \cite{Bus1}.

\subsubsection{}\label{xghagsx}

The Buscher rules provide local rules for the behaviour of the geometric
objects under $T$-duality on the target level. The underlying spaces of the
targets (being principal bundles on manifolds) are locally isomorphic. Therefore, topological $T$-duality is really
interesting only on the global level. 
The idea for setting up a chacterization of a topological $T$-dual comes from the desire to realize the
Buscher transformation rules globally. The analysis of this transition from geometry to topology has been started in the case of circle bundles, e.g.~\cite{bem2} and continued including the higher dimensional case with
\cite{bhm},\cite{bhm3}, \cite{bem1}, without stating a precise mathematical definition of topological $T$-duality there.

\subsubsection{}
Currently, such a precise mathematical definition of topological $T$-duality 
has to be given in an ad-hoc manner.
For $T^n$-bundles with twists we know three possibilities:
\begin{enumerate}
\item A definition in the framework of non-commutative geometry can be extracted from the works
\cite{mr}, \cite{mr1}, \cite{mr2} and will be explained
in \ref{sxhaxja}.
\item The homotopy theoretic definition used in the present paper is based on the notion of a
$T$-duality triple (see Definition \ref{dhqdjwq}).
\item Following an idea of T. Pantev, in a forthcoming paper \cite{bsst}
we propose a definition of topological  $T$-duality  for $T^n$-bundles with
twists using Pontrjagin duality for topological group stacks. 
\end{enumerate}
Surprisingly, all three definitions eventually lead to equivalent theories of
topological $T$-duality for $T^n$-bundles with twists (the equivalence of (1)
and (2) is shown in \cite{schn}, and the equivalence of (2) and (3) is shown
in \cite{bsst}). This provides strong evidence for the fact that these
definitions for topological T-duality correctly reflect the $T$-duality on the
target or even quantum theory level.

\subsubsection{}

If two spaces (with twist, i.e. $B$-field background) are in $T$-duality then
this has strong consequences on certain of their  topological invariants. For example, there are distinguished  isomorphisms (called $T$-duality isomorphisms, see Definition \ref{trd}) between their twisted cohomology groups and twisted $K$-groups.
The existence of these $T$-duality isomorphisms has already been observed in
\cite{bem1}, \cite{mr} and their follow-ups. The desire for a $T$-duality
isomorphism actually was one of our main guiding principle which led to the
introduction of the notion of a $T$-duality triple
and therefore our mathematical definition of topological $T$-duality.

\subsubsection{}

Having understood $T$-duality on the level of underlying topological spaces one can now lift back to the geometric level. We hope that the topological classification results (and  their natural generalizations to topological stacks in order to include non-free $T^n$-actions) will find applications to mirror symmetry in algebraic geometry and string theory.

\section{Topological $T$-duality via $T$-duality triples}

\subsubsection{}

In this Section we propose a mathematical set-up for topological $T$-duality of total spaces of $T^n$-bundles with twists and give detailed statements of our classification results. We will also shed some light on the relation with other pictures in the literature. 

\subsubsection{}

In the present paper we will use elements of the homotopy classification theory
of principal fibre bundles \cite[Ch 4]{hus}. Therefore, spaces in the present paper are always assumed to be Hausdorff and paracompact. 

\subsubsection{}

Let us fix a base space $B$ and $n\in \nat$. By $T^n:= \underbrace{U(1)\times\dots\times U(1)}_{\mbox{\scriptsize $n$-factors}}$ we denote the $n$-torus.
The fundamental notion of the theory is that of a pair.
\begin{ddd}\label{fhsdjfshwie}
A pair $(E,h)$ over $B$ consists of a principal $T^n$-bundle $E\to B$ and a cohomology class
$h\in H^3(E,\Z)$. An isomorphism of pairs $\phi\colon (E,h)\to (E^\prime,h^\prime)$ is an isomorphism
$$\xymatrix{E\ar[dr]\ar[rr]^\phi&&E^\prime\ar[dl]\\&B&}$$
of $T^n$-principal bundles such that $\phi^*h^\prime=h$.
We let $P(B)$ denote the set of isomorphism classes of pairs over $B$.
\end{ddd}
We can extend $P$ to a functor
$$P:\{\mbox{Spaces}\}^{op}\to \{\mbox{Sets}\}\ .$$
Let $f:B^\prime\to B$ be a continuous map and
$(E,h)\in P(B)$. Then we define 
$(E^\prime,h^\prime):=P(f)(E,H)$ as  the pull-back of $(E,h)$.  More precisely, the $T^n$-bundle $E^\prime\to B^\prime$ is defined by the pull-back diagram
$$\xymatrix{E^\prime\ar[r]^F\ar[d]&E\ar[d]\\B^\prime\ar[r]^f&B}\ ,$$
and $h^\prime:=F^*h$.

\subsubsection{}

% \subsubsection{}
% 
% Duality is a frequently occurring principle for understanding  a mathematical structure. The Fourier transformation is a classical example which leads to the Pontryagin duality of topological abelian groups. A more geometric example is the Fourier-Mukai transformation which in its simplest case relates the derived categories of an abelian variety and its dual. 
% 
% In the present paper we add a thorough investigation of a duality which
% relates  $T^n$-principal bundles equipped with a twist
% (think of a three-dimensional integral cohomology class) with their duals. 
% This duality is closely related to the examples mentioned above as we will explain later in this introduction. The role of the Fourier transformation or Fourier-Mukai transformation is played by the $T$-duality transformation which relates the twisted cohomology groups
% of the pair of dual $T^n$-bundles. 
% 
% 
% 
% 
% 
% 
% 
% 
% 
% 
% 
% 
% 
% 
% 
% 
%   
% 
% 
% 
% 
% 
% 
% 
% 
% 
% \subsubsection{}
%  
% The notion of $T$-duality has its origin in string theory.
% It asserts that a certain type of string theory on a target $F$
% partially compactified on a torus is equivalent to another type of
% string theory on a dual target $\hat F$. In the presence of an
% $H$-flux (the physicist's name for twist) the topological aspects of $T$-duality are of particular interest.
%  
% In the present paper we concentrate on the topological aspects of
% $T$-duality, and we neglect geometric questions and physical
% interpretations. 

The study of topological $T$-duality started with the case of circle bundles, i.e. $n=1$.
Guided by the experience obtained in \cite{bhm}, \cite{bem1}, \cite{bem2}, and
\cite{mr}, a
mathematical definition of topological $T$-duality for pairs in the case $n=1$ was given in  \cite{bunkeschickt}. In the latter  paper 
$T$-duality appears in two flavors. 

On the one hand, $T$-duality is a
relation (see \cite[Def. 2.9]{bunkeschickt}) which may or may not be satisfied by two pairs $(E,h)$ and $(\hat E,\hat h)$ over $B$.  The relation has a cohomological characterization. We will not recall the details of the definition here since it will be equivalent to Definition \ref{dhqdjwq} in terms of $T$-duality triples (reduced to the case $n=1$).

On the other hand we construct in \cite{bunkeschickt} a $T$-duality
transformation, a natural automorphisms of functors of order two 
\begin{equation}\label{hwdqjdw}
T\colon P\to P\ ,
\end{equation} which assigns to each pair $(E,h)$ a
specific $T$-dual $(\hat E,\hat h):=T (E,h)$.

The existence of such a transformation is a special property of the case $n=1$.
It has already been observed in
\cite{bhm}, \cite{mr}, and \cite{bunkeschickt} that such a transformation can not
exist for general higher dimensional torus bundles. 
The first reason is that for $n\ge 2$ not every pair admits a $T$-dual which
implies that $T$ in (\ref{hwdqjdw}) could at most be partially defined. An
additional obstruction (to a partially defined transformation) is the non-uniqueness of $T$-duals.

\subsubsection{}\label{cpudef}

In order to describe topological $T$-duality in the higher-dimensional ($n > 1$) case we introduce the notion of a $T$-duality triple. To this end we must categorify the third integral cohomology using the notion of twists.

% 
% 
% 
%  In the present paper we describe the, in a certain sense, optimal generalization of the results from $U(1)$ to $T^n$-bundles for higher $n$. 
% An important
% part is the precise characterization of those pairs $(F,h)$ of $T^n$-bundles
% $F$ and classes $h\in H^3(F,\Z)$ for which
% a $T$-dual pair exists, and a description of the set of such duals.
% In order to generalize the relation of $T$-duality we can not stick to cohomology.
% Note that $H^3(F,\Z)$ corresponds to the  set of  isomorphism classes of twists. In the present paper the  characterization of $T$-duality will  be  spelled out on the level of the twists themselves. It will be encoded in the notion of a $T$-duality triple.  
% 
%  

There are various models for twists, some of them are reviewed in
\ref{twistsexplain}. The reader not familiar with the concept of
twists and twisted cohomology theories is advised to consult this
appendix. The results of the present paper are independent 
of the choice of the model. Therefore, let us once and for all fix a
model for twists. 

Let us recall the essential properties of twists used in the constructions below. First of all we have a transformation
$$\xymatrix{\{\mbox{category of twists over $B$}\}/{\scriptsize \mbox{isomorphism}}\ar[rr]^{\hspace{2cm}\cong}&&H^3(B,\Z)}$$
which is natural in $B$. For a twist $\cH$ we let $[\cH]$ denote the cohomology
class corresponding to the isomorphism class of $\cH$.
Furthermore, given isomorphic twists $\cH,\cH^\prime$, the set
$\Hom_{\scriptsize\mbox{Twists}}(\cH,\cH^\prime)$ is a torsor over $H^2(B,\Z)$,
and this structure is again compatible with the functoriality in $B$. In this paper we frequently identify the based set of automorphisms
$\Hom_{\text{\scriptsize{Twists}}}(\cH,\cH)$  with
  $H^2(B,\Z)$.

For a twist $\cH$ over $B$ we will use the schematic notation
$$\xymatrix{\cH\ar@{.>}[r]&B}$$
which aquires real sense if one realizes twists as gerbes or bundles of compact operators over $B$.

\subsubsection{}\label{tr56512u}

We fix  
 an integer $n\ge 1$ 
% We start with the definition of an  $n$-dimensional $T$-duality triple.
% Let $T^n:=U(1)^n$ denote the $n$-dimensional torus. 
and a connected
base space $B$ with a base point $b\in B$. A $T^n$-principal
bundle $\pi:F\rightarrow B$ is classified by an $n$-tuple of Chern classes $c_1,\dots,c_n\in
H^2(B,\Z)$. Let $\hat \pi:\hat F\rightarrow B$ be a second
$T^n$-principal bundle with Chern classes $\hat c_1,\dots,\hat c_n\in
H^2(B,\Z)$.
 
 Let $\cH$ be a twist on $F$ such that its characteristic
class lies in the second filtration step of the Leray-Serre spectral sequence
filtration, i.e.~satisfies  $[\cH]\in \cF^2 H^3(F,\Z)$ (see \ref{lerayserreexplain} for notation). Furthermore we assume that its 
leading part fulfills
\begin{equation}\label{filt111}[\cH]^{2,1}=[\sum_{i=1}^n y_i\otimes \hat
  c_i]\in {}^\pi E_\infty^{2,1}\ ,\end{equation}
where $y_i$ are generators of the cohomology of the fibre $U(1)^n$ of $F$,
compare again Appendix \ref{lerayserreexplain}.
% See Appendix \ref{lerayserreexplain} for the notation used in this paper connected with the
% Leray-Serre spectral sequence. 
Similarly, let $\hat \cH$ be a twist
on $\hat F$ such that $[\hat \cH]\in \cF^2 H^3(\hat F,\Z)$
and (with similar notation) 
\begin{equation}\label{filt112}[\hat \cH]^{2,1}=[\sum_{i=1}^n \hat y_i\otimes c_i]\in {}^{\hat
  \pi} E_\infty^{2,1}\ .\end{equation}
  
We assume that we have an isomorphism of twists
$u:\hat p^*\hat \cH\rightarrow p^*\cH$.
as indicated in the diagram
\begin{equation}\label{rdef4}
\xymatrix{&p^*\cH\ar[dl]\ar@{.>}[dr]&&\hat p^*\hat \cH\ar[ll]^u\ar@{.>}[dl]\ar[dr]&\\\cH\ar@{.>}[dr]&&F\times_B\hat F\ar[dd]^r\ar[dl]^p\ar[dr]^{\hat p}&&\hat \cH\ar@{.>}[dl]\\&F\ar[dr]^\pi&&\hat F\ar[dl]^{\hat \pi}&\\&&B&&}\ .\end{equation}

% \begin{array}{ccccc}
% &&F\times_B\hat F&&\\
% &p\swarrow&&\searrow \hat p\\
% F&&r\downarrow &&\hat F\\
% &\pi\searrow&&\swarrow\hat \pi&\\
% &&B&&\end{array}\ .

We require that this isomorphism satisfies the condition
$\cP(u)$ which we now describe.
Let $F_b$ and $\hat F_b$ denote the fibers of $F$ and
$\hat F$ over
$b\in B$ and consider the induced diagram
$$\xymatrix{&F_b\times\hat F_b\ar[dl]^{p_b}\ar[dr]^{\hat p_b}&\\
F_b&&\hat F_b}\ .$$
The assumptions on  $\cH$ and $\hat \cH$ imply the existence of isomorphisms $v:\cH_{|F_b}\stackrel{\sim}{\rightarrow} 0$ and $\hat v:0\stackrel{\sim}{\rightarrow} \hat
\cH_{|\hat F_b}$.
We now consider the composition
\begin{equation}\label{uudef}u(b):=(0\stackrel{\hat p_b^*\hat v}{\rightarrow}\hat p_b^*
 \hat
\cH_{|\hat F_b}\stackrel{u_{|F_b\times \hat F_b}}{\rightarrow}
p_b^*\cH_{|F_b}\stackrel{p_{b}^*v}{\rightarrow }0)\in
H^2(F_b\times \hat F_b,\Z)\ .\end{equation}
The condition $\cP(u)$ requires that
\begin{equation}
[u(b)]=[\sum_{i=1}^n y_i\cup \hat y_i]\in H^2(F_b\times \hat
F_b,\Z)/(\im(p_b^*)+\im(\hat p_b^*)).\label{eq:condp}
\end{equation}

The class $[u(b)]$ in this
quotient is well-defined independent of the choice of $v$
and $\hat v$. 

\begin{ddd}\label{tr13}
An $n$-dimensional $T$-duality over $B$ triple is a triple
$$((F,\cH),(\hat F,\hat \cH),u)$$
consisting of $T^n$-bundles $\pi:F\to B$, $\hat \pi:\hat F\to B$, twists
$\xymatrix{\cH\ar@{.>}[r] &F}$, $\xymatrix{\hat \cH\ar@{.>}[r] &\hat F}$ satisfying (\ref{filt111}) and (\ref{filt112}), respectively, and an isomorphism $u:\hat p^*\hat \cH\stackrel{\sim}{\to} p^*\cH$ (for notation see (\ref{rdef4})) which satisfies condition $\cP(u)$.
\end{ddd}

\begin{ddd}
We will say that the triple $((F,\cH),(\hat F,\hat \cH),u)$
extends the pair $(F,[\cH])$ and connects the two pairs $(F,[\cH])$ and $(\hat F,[\hat \cH])$.
\end{ddd}

\subsubsection{}

We can now define our notion of topological $T$-duality based on $T$-duality triples.
\begin{ddd}\label{dhqdjwq}
Two pairs
$(F,h)$ and $(\hat F,\hat h)$ over $B$ are in $T$-duality
if here is a $T$-duality triple connecting them. 
\end{ddd}
The main results of the present paper concern the following problems:
\begin{enumerate}
\item classification of isomorphism  classes of $T$-duality triples over $B$
\item classification of $T$-duality triples which connect two given pairs
\item existence and classification of $T$-duality triples extending a given pair

\end{enumerate}

%The present paper represents an improvement over \cite{bhm} since it
%introduces  and studies $T$-duality integrally. 
%Moreover, we make precise the role of the identification of
%twists $u$.  

\newcommand{\Triple}{{\tt Triple}}
\subsubsection{}\label{tripleintro}

There is a natural notion of an isomorphism of $T$-duality triples. Its details will be spelled out in Definition \ref{isodd}. If $f:B\rightarrow B^\prime$ is a continuous map, and $x:=((F^\prime,\cH^\prime),(\hat F^\prime,\hat \cH^\prime),u^\prime)$ is a $T$-duality triple over $B^\prime$, then one
defines a $T$-duality triple  $((F,H),(\hat F,\hat H),u)=f^*x$ over $B$ in a canonical way. First of all the underlying $T^n$-bundles are given by the pull-back diagrams
$$\xymatrix{F\ar[r]^\phi\ar[d]&F^\prime\ar[d]\\ B\ar[r]^f&B^\prime} ,\quad \xymatrix{\hat F\ar[r]^{\hat \phi}\ar[d]&\hat F^\prime\ar[d]\\ B\ar[r]^f&B^\prime }\ .$$ Then we define the twists $\cH:=\phi^*\cH^\prime$ and $\hat \cH:=\hat \phi^*\hat\cH^\prime$. Finally we consider the induced map
$\psi:=(\phi,\hat \phi):F\times_B\hat F\to F^\prime\times_{B^\prime} \hat F^\prime$ and define $u$ as the composition
$$\hat p^*\hat \cH\cong \psi^*(\hat p^\prime)^* \hat \cH^\prime\stackrel{\psi^*u^\prime}{\to} \psi^*(p^\prime)^* \cH^\prime\cong p^*\cH$$
of natural isomorphisms and the pull-back of $u^\prime$ via $\psi$.

\begin{ddd}
We define the functor
$$\Triple_n:\{\mbox{spaces}\}^{op}\to \{\mbox{sets}\}$$
which associates to a
space $B$ the set of isomorphism classes $\Triple_n(B)$ of
$n$-dimensional $T$-duality triples over $B$. 
\end{ddd}

\subsubsection{}

In Lemma \ref{homo32} we will observe  that the functor $\Triple_n$  is homotopy invariant. In general, given a contravariant homotopy invariant functor from spaces to sets one asks whether it can be represented by a classifying space. If this is the case, then the  functor can be studied by applying methods of algebraic topology to its classifying space. Our study of the functor $\Triple_n$ follows this philosophy.

\subsubsection{}

In the following we describe a space $\bR_n$ which will turn out to be a classifying space of the functor $\Triple_n$ by Theorem \ref{main}.

Consider the product of two copies of the Eilenberg-MacLane space
$K(\Z^n,2)\times K(\Z^n,2)$ with canonical generators
$x_1,\dots,x_n$ and $\hat x_1,\dots, \hat x_n$ of the second
integral cohomology. We consider the class
$q:=\sum_{i=1}^n x_i\cup \hat x_i$ as a map
$q:K(\Z^n,2)\times K(\Z^n,2)\rightarrow K(\Z,4)$.
\begin{ddd}
Let $\bR_n$ be the homotopy fiber of $q$.
\end{ddd}

We consider the two components of the map $(\bc,\hat \bc):\bR_n\rightarrow 
K(\Z^n,2)\times K(\Z^n,2)$ as the classifying maps of two
$T^n$-principal bundles $\pi_n:\bF_n\rightarrow \bR_n$ and $\hat \pi_n:\hat
\bF_n\rightarrow \bR_n$. By a calculation of the cohomology of
$\bF_n$, $\hat \bF_n$ and $\bF_n\times_{\bR_n}\hat \bF_n$ we show the
following Theorem.  
\begin{theorem}[Theorem \ref{ex980}]
There exists a unique isomorphism class of $n$-dimensional $T$-duality triples
$[x_{n,univ}]=[(\bF_n,\cH_n),(\hat \bF_n,\hat \cH_n),\bu_n]\in
\Triple_n(\bR_n)$ with underlying $T^n$-bundles isomorphic to $\bF_n$ and $\hat \bF_n$. 
\end{theorem}
Let $P_n$ denote the set-valued functor classified by $\bR_n$.
%  where $[X,Y]$ denotes the homotopy
% classes of continuous maps from $X$ to $Y$. 
This functor associates to
$B$ the set $P_n(B)$ of homotopy classes $[f]$ of maps $f:B\to \bR_n$.
The universal triple $[x_{n,univ}]$ induces a natural transformation
of functors $\Psi_B:P_n\rightarrow \Triple_n(B)$ by
$$\Psi_B([f]):=\Triple_n(f)[x_{n,univ}] = f^*[x_{n,univ}]\ .$$
The following theorem characterizes $\bR_n$ as a classifying space of the functor $\Triple_n$.
\begin{theorem}[Theorem \ref{main12}]\label{main}
The natural transformation $\Psi$ is an isomorphism of functors.
\end{theorem}

\subsubsection{}\label{action4}

%\Kommentar{Although the next results are of interest in their own,
%  they are perhaps not important enough to take up so much space in
%  the introduction. I would suggest to remove all of this
%  subsubsection.}

In order to prove Theorem \ref{main} we must investigate the fine
structure of the functor $\Triple_n$. Of particular importance is the following
action of $H^3(B,\Z)$ on $\Triple_n(B)$ (see  
 \ref{ert2}). Let $x:=((F,\cH),(\hat F,\hat \cH),u)$ represent a class
$[x]\in \Triple_n(B)$, and let $\alpha\in H^3(B,\Z)$. We choose a twist
$\cV$ in the class $\alpha$ and set
$x+\cV:=((F,\cH\otimes\pi^*\cV),(\hat F,\hat \cH\otimes \hat
\pi^*\cV),u\otimes r^*\id_{\cV})$ (see (\ref{rdef4}) for  the definition of $r$).  Then we define $[x]+\alpha:=[x+\cV]$.

We now consider the set $\Triple_n^{(F,\hat F)}(B)$  of isomorphism
classes of $n$-dimensional $T$-duality triples over fixed
$T^n$-bundles $F$ and $\hat F$ (see \ref{overdef}).
The group $H^3(B,\Z)$ acts naturally
on $\Triple_n^{(F,\hat F)}(B)$ by the same construction as above.

\begin{prop}[Proposition \ref{freetrans}]\label{prop:freetrans0}
$\Triple_n^{(F,\hat F)}(B)$ is an $H^3(B,\Z)$-torsor.
\end{prop}

\subsubsection{}

In terms of the classifying spaces, fixing $F$ and $\hat F$ corresponds
to fixing classifying maps $(c,\hat c):B\rightarrow K(\Z^n,2)\times
K(\Z^n,2)$. The set 
$\Triple_n^{(F,\hat F)}(B)$ then corresponds to the set of 
homotopy classes of lifts in the diagram
$$\xymatrix{&\bR_n\ar[d]^{(\bc,\hat \bc)}\\
B\ar@{.>}[ur]^f\ar[r]_{\hspace{-2cm}(c,\hat c)}&K(\Z^{n},2)\times K(\Z^n,2)}\ .$$
Since the homotopy fiber of $(\bc,\hat \bc)$ has the homotopy type of a $K(\Z,3)$-space it is
clear by obstruction theory that $H^3(B,\Z)$ acts freely and transitively on the set of such
lifts. In combination with Proposition \ref{prop:freetrans0} this leads to the key step in the
proof that $\bR_n$ is the correct classifying space.

\subsubsection{}\label{eins7}

Let now $\psi$ and $\hat \psi$ be bundle automorphisms of $F$ and $\hat
F$. We can realize $\psi$ and $\hat \psi$ as right multiplication by
maps $\psi,\hat \psi:B\rightarrow T^n\cong K(\Z^n,1)$.
In this way the homotopy classes of $\psi$ and $\hat \psi$ can be
considered as classes $[\psi],[\hat \psi]\in H^1(B,\Z^n)$. 
Let $x:=((F,\cH),(\hat F,\hat \cH),u)$ be an $n$-dimensional
$T$-duality triple. Then we form the triple
$x^{(\psi,\hat \psi)}:=((F,\psi^*\cH),(\hat F,\hat\psi^*\hat
\cH),(\psi,\hat\psi)^*u)$. We introduce the notation
$\hat c\cup[\psi]:=\sum_{i=1}^n \hat c_i\cup [\psi]_i\in H^3(B,\Z)$,
where $\hat c_1,\dots,\hat c_n$
are the components of the Chern class of $\hat F$, and
$[\psi]_1,\dots,[\psi]_n$ are the components of $[\psi]$.
We define $c\cup[\hat \psi]$ similarly. Then we show:
\begin{prop}[Proposition \ref{autact}] In $\Triple_n^{(F,\hat F)}(B)$ we
have $$[x^{(\psi,\hat \psi)}]=[x]+\hat c\cup [\psi]+c\cup [\hat
\psi]\ .$$
\end{prop}

There is a natural forgetful map $$\Psi:\Triple_n^{(F,\hat F)}(B)\rightarrow
\Triple_n(B)\ .$$ Recall the definition (\ref{rdef4}) of the map $r$ and note that
$\im({}^r d_2^{2,1})\subseteq H^3(B,\Z)$ (see \ref{lerayserreexplain} for notation) is exactly the subgroup of elements which can be written in the form $c\cup a+\hat c\cup b$ for $a,b\in H^1(B,\Z^n)$.
Proposition \ref{autact} immediately implies:
\begin{kor}
If $\alpha\in \im({}^r d_2^{2,1})\subseteq H^3(B,\Z)$, then we have
$\Psi([x]+\alpha)=\Psi([x])$.
 \end{kor}

\subsubsection{}

Let $(e_1,\dots,e_n,\hat e_1,\dots,\hat e_n)$ be the standard basis of $\Z^{2n}$.
Let $O(n,n,\Z)\subset GL(2n,\Z)$ be the subgroup of transformations which
fix the quadratic  form $q:\Z^{2n}\rightarrow \Z$ with
$q(\sum_{i=1}^n a_i e_i+b_i\hat e_i):=\sum_{i=1}^n a_ib_i$.
 
\begin{prop}[Lemma \ref{re}]
The group $O(n,n,\Z)$ acts by homotopy equivalences on $\bR_n$.
We have an induced action of $O(n,n,\Z)$ on the functor $\Triple_n$ by
automorphisms.
\end{prop}

In the literature this group is sometimes called the \emph{$T$-duality
group}.

\subsubsection{}\label{fff564}

% A pair $(F,h)$ over $B$ consists of a $T^n$-principal bundle
% $\pi:F\rightarrow B$ and a class $h\in H^3(F,\Z)$. Moreover, if $f:B\rightarrow
% B^\prime$ is a continuous map, and $(F,h)$ is a pair over $B^\prime$,
% one defines a pair $f^*(F,h)$ in a natural way. We
% obtain a set-valued functor $B\mapsto \tilde P_{(0)}$ which associates
% to each space $B$ the set of isomorphism classes of pairs $\tilde
% P_{(0)}(B)$ over $B$. The functor $\tilde P_{(0)}$ is homotopy
% invariant.

Recall the definition \ref{fhsdjfshwie} of the functor $B\mapsto P(B)$ which associates to a space $B$ the set of isomorphism classes of $n$-dimensional pairs over $B$. We will write $\tilde P_{(0)}:=P$ since this functor appears at the lowest level
of a tower of functors $\tilde P_{(0)}\leftarrow \tilde P_{(1)}\leftarrow \dots$ (see \ref{ghasjhasd}).
In the notation for these functors we will not indicate the dimension $n$ of the torus $T^n$ explicitly.

The functor $\tilde P_{(0)}$ is homotopy invariant (the proof of \cite[Lemma
2.2]{bunkeschickt} extends from the case $n=1$ to arbitrary $n\ge 1$). 
Generalizing again the approach of \cite{bunkeschickt}
from the case $n=1$ to general $n\ge 1$ we construct a 
classifying space $\tilde \bR_n(0)$ for the functor $\tilde
P_{(0)}$ as follows.
Let $U^n\rightarrow K(\Z^n,2)$ be the universal $T^n$-bundle. Then we define
$$\tilde \bR_n(0):=U^n\times_{T^n} \Map(T^n,K(\Z,3))\ .$$
The natural map $\tilde \bR_n(0)\rightarrow K(\Z^n,2)$ classifies
a $T^n$-principal bundle $\tilde \bF_n(0)\rightarrow \tilde \bR_n(0)$ which
admits a natural map $\tilde \bF_n(0)\rightarrow K(\Z,3)$. We interpret
the homotopy class of this map as a class $\tilde \bh(0)\in H^3(\tilde
\bF_n(0),\Z)$. The isomorphism class of the universal pair $[\tilde \bF_n(0),
\tilde \bh(0)]\in \tilde P_{(0)}(\tilde \bR_n(0))$ induces a natural transformation of functors
$\tilde v_B:[B,\tilde \bR_n(0)]\rightarrow \tilde P_{(0)}(B)$ (see
Lemma \ref{cl1}) which turns out to be an isomorphism.

\subsubsection{}

In Section \ref{onecon} we introduce the one-connected cover $\tilde
\bR_n\rightarrow \tilde \bR_n(0)$. It is the universal covering of a
certain connected component of $\tilde \bR_n(0)$.
The first entry of the universal triple $((\bF_n,\cH_n),(\hat
\bF_n,\hat \cH_n),u_n)$ over $\bR_n$ gives rise to a classifying map
$$f(0):\bR_n\rightarrow \tilde \bR_n(0)\ .$$
We shall see (Lemma \ref{var}) that $f(0)$ has a factorization 
\begin{equation}\label{gshgdgasdszu}\xymatrix{&\tilde \bR_n\ar[d]\\
\bR_n\ar@{.>}[ur]^f\ar[r]^{f(0)}&\tilde \bR_n(0)} \ .\end{equation}
Note that the factorization $f$ is not unique.
\begin{theorem}[Theorem \ref{var}]\label{z56t}
The map $f:\bR_n\rightarrow \tilde \bR_n$ is a weak homotopy
equivalence.
\end{theorem}

\subsubsection{}

There are two natural transformations of functors
$$\xymatrix{&\Triple_n\ar[dr]^{\hat s}\ar[dl]^s&\\P&&P}\ ,$$
where 
$$s((F,\cH),(\hat F,\hat \cH),u):=(F,[\cH])\ ,\quad \hat s((F,\cH),(\hat F,\hat \cH),u):=(\hat F,[\hat \cH])\ .$$
The problem of the existence and the classification of $T$-duals of a pair
$(F,h)\in P(B)$ is essentially a question about the fibre $s^{-1}(F,h)\subseteq\Triple_n(B)$. The transformation $s$ is realized on the level of classifiying spaces by the map
$$\tilde \bR_n\to \tilde \bR_n(0)$$ in (\ref{gshgdgasdszu}).
This allows to translate questions about the fibres of $s$ to homotopy theory.

\subsubsection{}\label{after43}

Consider a pair $(F,h)$ over a space $B$. The representatives of elements of $s^{-1}(F,h)$ will be called extensions of $(F,h)$. 
\begin{ddd}\label{exexdef}
An extension of $(F,h)$ to an $n$-dimensional $T$-duality triple is an
$n$-dimensional $T$-duality triple $((F,\cH),(\hat F,\cH),u)$ over $B$ such
that
$[\cH]=h$.
\end{ddd}
The difference between the notions of an extension of $(F,h)$ and an element
in the fibre $s^{-1}(F,h)$ is seen on the level of the notion of an isomorphism of extensions (see Definition \ref{isoextensiondef}).
Roughly speaking, an isomorphism of extensions of $(F,h)$
is an isomorphism of triples such that the underlying bundle
isomorphism of $F$ is the identity. 
\begin{ddd}
We let $\Ext(F,h)$ denote the set
of isomorphism classes of extensions of $(F,h)$ to $n$-dimensional
$T$-duality triples. 
\end{ddd}
We have a natural surjective map
$$\Ext(F,h)\to s^{-1}(F,h)$$
which in general may not be injective.

\subsubsection{}

We then consider the following two
problems.
\begin{enumerate}
\item Under which conditions does $(F,h)$ admit an extension, i.e.~is the  set $\Ext(F,h)$ non-empty?
\item Describe the set $\Ext(F,h)$.
\end{enumerate}
Answers to these questions settle the problem of existence and classification of $T$-duals of $(F,h)$ in the following sense.
\begin{enumerate}
\item The pair $(F,h)$ admits a $T$-dual if and only if $\Ext(F,h)$ is not empty.
\item The set of $T$-duals of $(F,h)$ can be written as
$\hat s(\Ext(F,h))\subseteq P(B)$.
\end{enumerate}

\subsubsection{}

As a consequence of Theorem \ref{z56t} we derive the following answer
to the first question.
\begin{theorem}[Theorem \ref{exex}]
The pair $(F,h)$ admits an extension to a $T$-duality triple $((F,\cH),(\hat F,\hat \cH),u)$ if and
only if $h\in \cF^2H^3(F,\Z)$.  
\end{theorem}

In particular, the condition $h\in \cF^2H^3(F,\Z)$ is a necessary and
sufficient condition for the existence of a $T$-dual to $(F,h)$.
If we write out the leading part of $h$ as
$h^{2,1}=[\sum_{i=1}^n y_i\otimes \hat c_i]\in {}^\pi E^{2,1}_\infty$,
then we can read off some information about the
Chern classes $\hat c_1,\dots,\hat c_n$ of the $T$-dual bundle $\hat F$.
In fact we have ${}^\pi E^{2,1}_\infty={}^\pi E^{2,1}_2/\im({}^\pi
d^{0,2}_2)$, and ${}^\pi
d^{0,2}_2(\sum_{i<j} A_{i,j} y_i\cup y_j)=\sum_{i,j} (A_{i,j}
-A_{j,i}) y_j\otimes c_i$,
where we set $A_{i,j}:=0$ for $i\ge j$. 

It follows that the Chern classes $\hat c_i$ of the dual bundle $\hat F$ are determined by the pair $(F,h)$ up to a change
$\hat c_i\mapsto \hat c_i +\sum_{j=1}^nB_{i,j} c_j$ for some
antisymmetric matrix $B\in \Mat(n,n,\Z)$.  Of course, the classes $\hat c_i$ are completely determined  by the choice of an extension of $(F,h)$.

We fix a pair $(F,h)$. If $n\ge 2$, then even the topology of the $T$-dual bundle $\hat F$
may depend on the choice of the extension of the pair $(F,h)$ to a
$T$-duality triple. We have already  demonstrated this by an example in 
 \cite[Section 4.4]{bunkeschickt}.  

The discussion above gives a description of the topological invariants of
a\footnote{At the moment of writing this paper is was not clear that the notions of a $T$-dual
  used in the present paper coincides with that of \cite{mr1}. Meanwhile results in this direction have been obtained in \cite{schn}.}
 $T$-dual of $(F,h)$  in terms of  relations between the Chern classes $c,\hat c$ and the
 $H^3$-classes $[\cH]$ and $[\hat \cH]$. Our results improve the results of
 \cite{mr1}\footnote{ The calculation in \cite{mr1} is not correct since the
inclusion $X_j\rightarrow X$ in (4) of \cite{mr1} does not exist in
general.}.

\subsubsection{}\label{casjcka}

% In order to attack the problem of classification we must first make
% precise our notion of an isomorphism between extensions
% \ref{isoextensiondef}. Roughly speaking, an isomorphism of extensions
% of $(F,h)$
% is an isomorphism of triples such that the underlying bundle
% isomorphism of $F$ is the identity. 
% \begin{ddd}
% We let $\Ext(F,h)$ denote the set
% of isomorphism classes of extensions of $(F,h)$ to $n$-dimensional
% $T$-duality triples. 
% \end{ddd}
%The following theorem gives a complete
%description of $\Ext(F,h)$. 
Let
$x:=((F,\cH),(\hat F,\cH),u)$ represent a class $\{x\}\in \Ext(F,h)$. Then we let $c(x),\hat
c(x)$ denote the Chern classes of $F,\hat F$.
Note that the group $\ker(\pi^*)\subseteq H^3(B,\Z)$ acts on
$\Ext(F,h)$. Furthermore we have a homomorphism
$C:H^1(B,\Z^n)\rightarrow \ker(\pi^*)$ given by
$C(a):=\sum_{i=1}^n c_i\cup a_i$, where $a_i$ denotes the components
of $a$. 

The first two assertions of the following theorem have already been discussed above.
\begin{theorem}\label{class76}
\begin{enumerate}
\item The set $\Ext(F,h)$ is non-empty if and only if $h\in
  \cF^2H^3(F,\Z)$.
\item If $ \{x\}\in \Ext(F,h)$,  then for any antisymmetric matrix $B\in \Mat(n,n,\Z)$
  there exists $\{x^\prime\}\in \Ext(F,h)$ such that
$\hat c_i(x^\prime)=\hat c_i(x)+\sum_{j=1}^n B_{i,j}c_j(x)$.
Vice versa, if $\{x^\prime\}\in \Ext(F,h)$, then $\hat c(x^\prime)$ is of
this form.

\item 
If $\{x\}\in \Ext(F,h)$, then 
$\{\{x^\prime\}\in\Ext(F,h)|\hat c(x^\prime)=\hat c(x)\}$ is an orbit
under the effective action of $\ker(\pi^*)/\im(C)$.
\end{enumerate}
\end{theorem}
We prove this theorem in \ref{provecla}.
Note that in terms of the Leray-Serre spectral sequence of
$\pi:F\rightarrow B$  we can identify
$\ker(\pi^*)/\im(C)\cong \im ({}^\pi d_3^{0,2})$.

\subsubsection{}

Our work on $T$-duality was inspired by \cite{bhm}. In \cite{bhm} the authors study
smooth $T^n$-bundles equipped with real valued cohomology classes
$h_\R$. Guided by the principles which we  explained in \ref{xghagsx}, in \cite{bhm} the  $T$-dual pair is constructed geometrically using differential forms. 
Furthermore it was observed in \cite{bhm} that a condition  similar to $h_\R\in \cF^2
H^3(F,\R)$ is necessary and sufficient for the construction of a $T$-dual to work.
The problem of non-uniqueness of the construction was not  addressed in that paper.
In the present paper we provide a precise counterpart of
\cite{bhm} including the full information of integral cohomology.

%\subsubsection{}

%Note that a $T$-dual pair of $(F,h)$ is not unique. This has already
%been observed in \cite{bunkeschickt}, Sec. 4.4. One of the results of
%the present paper is the detailed
%description of set of possible choices. 

%In order to define a $T$-dual of a pair $(F,h)$ 
%we must choose an element $x\in v_{B}^{-1}([F,h])$. 
%In general the isomorphism class $v_B(T_B(x))$ of the $T$-dual pair 
%will depend non-trivially on the choice of $x\in  v_{B}^{-1}([F,h])$.

%\subsubsection{}

%In order to study the fibers of $v_B$ we construct
%a non-canonical isomorphism (see Corollary \ref{hhll})
%with another functor $B\mapsto\tilde P_n(B)$. The elements of $\tilde
%P_n(B)$ have a geometric meaning, and the fiber of the forgetful map
%$\tilde v_B:\tilde P_n(B)\rightarrow \tilde P_{(0)}(B)$ can be studied by
%obstruction theory (see \ref{clh} and \ref{clh1}). It seems to be
%difficult to formulate a general result. But let us mention the
%following opposite extremes. 
%\begin{enumerate}
%\item If $B=S^3$, then the transformation
%$\tilde v_{S^3}\colon \tilde P_n(S^3)\rightarrow P_0(S^3)$ is an inclusion of a subset.
%\item If $B=\bR_n$ is the space  introduced in
%\ref{rt3},
%then the fibers of $\tilde v_{\bR_n}:P_n(\bR_n)\rightarrow \tilde P_{(0)}(\bR_n)$ are
%torsors over the group $H^2(T^n,\Z)$.
%\end{enumerate}

\subsubsection{}

An interesting feature of topological $T$-duality is the $T$-duality isomorphism in twisted cohomology theories.
In \cite[Definition 3.1]{bunkeschickt}, compare \ref{tad3}, we introduced axioms (in particular the concept of
admissibility) for a twisted cohomology theory 
so that one can define for two pairs in $T$-duality a
$T$-duality transformation which is then proved to be an isomorphism.
Examples of twisted cohomology theories satisfying the admissibility axioms are twisted $K$-theory and twisted
real (periodic) cohomology. See \ref{twistsexplain} for a definition of twisted $K$-theory.

The existence of a $T$-duality isomorphism has been previously observed in  
\cite{bhm}, \cite{bem1}, \cite{bem2}, \cite{mr} and was our main guiding principle
for the definition of the $T$-duality relation.
 
\subsubsection{}

Let us assume that $((F,\cH),(\hat F,\hat \cH),u)$ is a $T$-duality
triple.  Let $(X,\cH)\mapsto h(X,\cH)$ be a twisted cohomology
theory, where the notation suggests that the cohomology groups depend
on two entries in a functorial way, namely the space $X$ and the twist $\cH$.  
The following definition uses the notation of (\ref{rdef4}).
\begin{ddd}\label{trd} The $T$-duality transformation 
 $T:h(F,\cH)\rightarrow h(\hat F,\hat \cH)$ is defined by
$$T:=\hat p_!\circ u^*\circ p^*\ .$$
\end{ddd}
Note that $T$ shifts degrees by $-n$. Furthermore, it is linear over $h(B)$.

\begin{theorem}[Theorem \ref{isot}]\label{rrr}
If the twisted cohomology theory $h(\dots,\dots)$ is $T$-admissible 
%\footnote{see\ref{tad3} or \cite[Definition 3.7]{bunkeschickt}},
 then the $T$-duality transformation is an isomorphism.
\end{theorem}

%% The two most important examples of
%% $T$-admissible cohomology  theories are twisted K-theory and
%% twisted cohomology with
%% coefficients in the  graded ring $\R[x,x^{-1}]$, where $x$ has degree
%% $2$. 

\subsubsection{}

A $T$-duality isomorphism for twisted $K$-theory and twisted periodic de Rham
cohomology was also obtained in \cite{bhm}. In contrast to this paper,
we take torsion in the third cohomology into account.
We refer to \cite[Section 4.3]{bunkeschickt} for an
explicit example which shows that the torsion part plays a significant
role.

\subsubsection{}\label{sxhaxja}

The approach of  \cite{mr}, \cite{mr1}, and \cite{mr2} to $T$-duality uses ideas from noncommutative topology.    The class $h\in H^3(F,\Z)$ is interpreted as the Dixmier-Douady class
of a unique isomorphism class of a stable continuous
trace algebra $\cA:=\cA(F,h)$ with spectrum $F$. Equivalently, the twist $\xymatrix{\cH\ar@{.>}[r]&F}$ is a bundle of $C^*$-algebras with fibre
the compact operators on a separable infinite-dimensional complex Hilbert space
such that $[\cH]=h$ (see \ref{twistsexplain}). Then we can write
$\cA(F,h)\cong C_0(F,\cH)$ (we assume for simplicity that $B$ (and hence $F$) is locally compact).

The authors study the question
of lifting the $T^n$-action on $F$ to an $\R^n$-action on $\cA$
such that the Mackey invariant is trivial. In this case the crossed
product $\hat \cA:=\cA\rtimes \R^n$ is again a continuous trace
algebra with a spectrum $\hat F$ which is a $T^n$-principal bundle over $B$.
Let $\hat h\in H^3(\hat F,\Z)$ denote the Dixmier-Douady class
of $\hat \cA$. From the point of view of \cite{mr}, \cite{mr1}, \cite{mr2}
the pair $(\hat F,\hat h)$ is the $T$-dual of $(F,h)$.

There is an obvious similarity of the following notions and their role in the theory of topological $T$-duality.
\begin{itemize}
\item  $\R^n$-action on $\cA(F,h)$ lifting the $T^n$-action on $F$ with trivial Mackey obstruction
\item Extension of $(F,h)$ to a $T$-duality triple.
\end{itemize}
The equivalence of the two approaches is established in \cite{schn}.

In the approach of \cite{mr}, the $T$-duality isomorphism for twisted
$K$-theory is equivalent to an isomorphism $K(\cA)\cong K(\hat \cA)$.
In fact, this isomorphism is   Connes'  Thom isomorphism for crossed
products with $\R^n$.
.  

Using the approach via noncommutative topology,
the natural two problems are to decide under which conditions
the required $\R^n$-action on $\cA$ exists, and to study the set
of choices for such an action. A satisfactory picture can be obtained
in the cases $n=1$ and $n=2$. The case $n=1$ is easy and has been
reviewed in \cite{bunkeschickt}. The main results of  \cite{mr} deal with
the case $n=2$. The necessary and sufficient condition for the
existence of the $\R^n$ action is again that $h\in \cF^2 H^3(F,\Z)$.
It is then claimed in \cite{mr}, that the action is unique. This is
not always true. In fact, 
it follows from the diagram given in \cite{mr}, Theorem 4.3.3, and the
observation that $d_2^{\prime\prime}$ (we use the notation of \cite{mr}) factors over
$p_!:H^2(F,\Z)\rightarrow H^0(B,\Z)$, that
the group
$H^0(B,\Z)/\im(p_!)$ acts freely on the set of $\R^n$-actions
with trivial Mackey invariant lifting the $T^n$-action on $F$.

%For general $n\ge 3$ the relation between the  
%the classification of the
%$\R^n$-actions on $\cA(F,h)$ with trivial Mackey invariant and the classification of the elements in
%$v_B^{-1}([F,h])$ is not yet clear to us.

\subsubsection{}

Let us mention a very interesting aspect of \cite{mr} which goes beyond the theory covered by the present paper. As explained above the necessary and sufficient condition for the existence of a $T$-dual of $(F,h)$ is $h\in \cF^2 H^3(F,\Z)$. The lift of the $T^n$ action on $F$ to an $\R^n$-action on $\cA(F,h)$ exists if $h\in \cF^1 H^3(F,\Z)$, but under this weaker condition
one might encounter a non-trivial Mackey obstruction. In this case the crossed product
$\hat \cA:=\cA\rtimes \R^n$ is not the algebra of sections of a bundle of compact operators
on a $T^n$-principal bundle over $B$. It has been observed in \cite{mr} (case $n=2$)
that one can interpret $\hat \cA$ as an algebra of sections of a bundle of noncommutative tori over $B$. In other words, the $T$-dual of $(F,h)$ can be realized as a bundle of noncommutative tori.
For further discussion of this phenomenon see also \cite{mr1} and \cite{mr2}.

\subsubsection{}
In \cite{bhm3} the point of view of \cite{mr,mr1} is generalized even
further by considering torus bundles with a completely arbitrary H-flux differential
$3$-form with integral periods\footnote{In the present paper's language  \cite{bhm2} considers pairs $(F,h)$ without any condition on $h\in H^3(F,\Z)$}. It is then argued that the resulting dual should be a bundle of
non-associative non-commutative tori.

\subsubsection{}

Assume that $F=B\times T^n$ is the trivial $T^n$-bundle, and
that we consider the trivial twist $h=0$.
Then the $T$-dual bundle is again the trivial bundle, $\hat F=B\times
T^n$, and the dual twist vanishes: $\hat h=0$.
In this situation the $T$-duality transformation
$T:K(B\times T^n)\rightarrow K(B\times T^n)$
is a $K$-theory version of the Fourier-Mukai transformation
(see e.g. \cite{muk}).

Note that the algebraic geometric analog is more precise.
In this case $F$ and $\hat F$ are bundles of dual abelian varieties.
On $F\times_B\hat F$ one has the so-called Poincar{\'e} sheaf $\cP$. Its
first Chern class $c_1$ (considered as an automorphism of the trivial twist)
satisfies the condition $\cP(c_1)$.
The Fourier-Mukai transformation is a functor $T:D^b(F)\rightarrow D^b(\hat F)$ between bounded derived
categories of coherent sheaves given
on objects by $T(X)=R\hat p_*(\cP\stackrel{L}{\otimes} p^*X)$.
Thus the $T$-duality transformation considered
in the present paper is a coarsification of the Fourier-Mukai transformation since it takes in a certain
sense only the isomorphism classes of objects into account. The tensor product with the Poincar{\'e} sheaf plays the role of an
automorphism
of the trivial twist.

A bundle of abelian varieties has a section. Therefore this case
corresponds to the case of trivial $T^n$-bundles in the present
paper. Non-trivial bundles can be interpreted as bundles of
torsors. In this case a good analog of the Poincar{\'e}
bundle such that $\cP(c_1)$ (see \ref{eq:condp}) is satisfied may not exist. In the topological situation
we must replace the Poincar{\'e} bundle by an isomorphism $u$ of
non-trivial twists in order to satisfy $\cP(u)$, and to have a $T$-duality {\em isomorphism}.
In algebraic geometry a similar observation is known (see e.g \cite{math.AG/0306213}), where twists are
represented by Azumaya algebras.

% As a specific example in the topological case, let  $F\rightarrow B$ be some possibly non-trivial
% $T^n$-bundle,  and let the  twist on $F$ be  trivial (i.e. $h=0$). Then
% the pair $(F,0)$ admits a $T$-dual pair $(\hat F,\hat h)$. 
% This is not unique in general
% (see  \cite[Section 4.4]{bunkeschickt}), but one choice is
% the trivial bundle
% $\hat F:=B\times T^n$. In order to describe
% the class $\hat h\in H^3(\hat F,\Z)$, we consider $\hat F$ as the
% fiber product of $n$-copies of a trivial $T^1$-bundle $F_1:=B\times T^1$.
% Let $\pr_i:\hat F\rightarrow F_1$ denote the corresponding
% projections. Then one possible choice of $\hat h$ is $\hat h=\sum_{i=1}^n \pr_i^* h_1$,
% where $h_i=\bc_i\times \orient_{T^1}\in H^3(B\times T^1,\Z)$,
% $\bc_i\in H^2(B,\Z)$ denotes the Chern classes of $F$, and
% $\orient_{T^1}\in H^1(T^1,\Z)$ is the orientation class.
% This example shows that it is necessary to introduce twists in order
% to generalize the Fourier-Mukai transformation to a non-trivial torus
% bundle.

\section{The space $\bR_n$}\label{fd}

\subsubsection{}\label{expla}

If $G$ is an abelian group and $k\in \nat$, then we consider the
homotopy type $K(G,k)$ of the Eilenberg MacLane space. It is characterized by $\pi_i(K(G,k))\cong 0$
for $i\not=k$, and $\pi_k(K(G,k))\cong G$. We denote a $CW$-complex
of this homotopy type by the same symbol.

The Eilenberg-MacLane space $K(G,k)$ classifies the cohomology functor
$H^k(\dots,G)$. In fact, there is a universal class $z\in
H^k(K(G,k),G)$ such that
$f\mapsto f^*(z)$ induces a natural isomorphism $[B,K(G,k)]\rightarrow
H^k(B,G)$, where $[B,K(G,k)]$ denotes homotopy classes of maps.

Occasionally, we will interpret $K(\Z,2)$ also as the
classifying space of the topological group $T^1:=U(1)$. An explicit model is $U/T^1$,
where $U$ is the unitary group of a separable infinite dimensional complex 
Hilbert space with the strong topology. The bundle $U\rightarrow U/T^1$ is the universal
$T^1$-principal bundle. Note further that
$K(\Z^n,2)\cong K(\Z,2)^n$ has the homotopy type of $BT^n$, and
that this space carries an universal $T^n$-bundle $U^n\rightarrow K(\Z,2)^n$.

\subsubsection{}\label{rndefd}

Let $x\in H^2(K(\Z,2),\Z)$ denote the canonical generator.
The product of Eilenberg-MacLane spaces
$K(\Z^n,2)\times K(\Z^n,2)$
can be written in the form
$$K(\Z^n,2)\times K(\Z^n,2)\cong K(\Z,2)^n\times K(\Z,2)^n\ .$$
We let
$$p_i,\hat p_i:K(\Z^n,2)\times K(\Z^n,2)\to K(\Z,2)\ ,\quad i=1,\dots,n$$ denote the projections onto the components of the first and second factor. Then we define the canonical generators
 $x_i,\hat x_i\in H^2(K(\Z^n,2)\times K(\Z^n,2),\Z)$, $i=1,\dots,n$,  by $x_i:=p_i^*(x)$ and $\hat x_i:=\hat p_i^*(x)$, respectively.
 
Let $q:K(\Z^n,2)\times K(\Z^n,2)\rightarrow K(\Z,4)$ be the map
classifying 
$$x_1\cup \hat x_1 + \dots +x_n\cup \hat x_{n}\in H^4(K(\Z^n,2)\times K(\Z^n,2),\Z)\ .$$
\begin{ddd} \label{rt3}
We define the homotopy type $\bR_n$ by the homotopy pull-back diagram
$$\begin{array}{ccc}
\bR_n&\stackrel{(\bc,\hat \bc)}{\rightarrow} &K(\Z^n,2)\times K(\Z^n,2)\\
\downarrow&&q\downarrow\\
*&\rightarrow&K(\Z,4)\end{array} \ .$$
\end{ddd}
In other words, $\bR_n$ is defined as the homotopy fiber of $q$.

\subsubsection{}

For later use we determine the homotopy groups of
$\bR_n$.
\begin{lem}\label{hh}
The homotopy groups of $\bR_n$ are given by
$$\begin{array}{|c|c|c|c|c|c|}\hline
i&0&1&2&3&\ge 4\\\hline
\pi_i(\bR_n)&*&0&\Z^{2n}&\Z&0\\\hline\end{array}
\ .$$\end{lem} 
\proof
The homotopy fiber of $(\bc,\hat \bc)$ is homotopy equivalent to the
one of
$*\rightarrow K(\Z,4)$, i.e to $K(\Z,3)$.
The assertion follows immediately from the long exact 
sequence of homotopy groups.
\hB

\subsubsection{}\label{tz1}

We write $\bc=(\bc_1,\dots,\bc_n)$ and $\hat \bc=(\hat \bc_1,\dots,\hat
\bc_n)$, i.e. we let $\bc_i$ and $\hat\bc_i$ denote the components of
$\bc$ or $\hat\bc$, respectively. 
\begin{lem}\label{gh}
We have
$$
\begin{array}{|c|c|c|c|c|c|}  
\hline
i&0&1&2&3&4\\\hline
H^i(\bR_n,\Z)&\Z&0&\Z^{2n}&0&\Z^{n(2n+1)-1}\\\hline \end{array}\ .
$$
Here $H^2(\bR_n,\Z)$ is freely generated by the components of $\bc$ and $\hat
\bc$, and $H^4(\bR_n,\Z)$ is generated by all possible products of the
components of $\bc$ and $\hat \bc$ subject to one relation
$$0=\bc_1\cup \hat \bc_1+\dots+\bc_n\cup \hat \bc_n\ .$$
\end{lem}
\proof
Recall from the proof of Lemma \ref{hh} that the homotopy fiber of
$\bR_n\rightarrow K(\Z^n,2)\times K(\Z^n,2)$ is a $K(\Z,3)$. The relevant part of the second page of the Leray-Serre
spectral sequence (see \ref{lerayserreexplain}) ${}^{(\bc,\hat \bc)}E_2^{p,q}\cong
H^p(K(\Z^{n},2)\times K(\Z^n,2),H^q(K(\Z,3),\Z))$ therefore
becomes
$$\begin{array}{|c|c|c|c|c|c|}
\hline3&\Z&0&*&0&*\\\hline
2&0&0&0&0&0\\\hline
1&0&0&0&0&0\\\hline
0&\Z&0&\Z^{2n}&0&\Z^{\frac{2n(2n+1)}{2}}\\ \hline
q/p&0&1&2&3&4\\\hline \end{array}\ .$$
We read off that $H^2(\bR_n,\Z)\cong \Z^{2n}$
is generated by the components of $\bc$ and $\hat \bc$.
The group ${}^{(\bc,\hat \bc)}E_2^{4,0}\cong H^4(K(\Z^n,2)\times K(\Z^n,2),\Z)$ is freely generated by all
possible products of the components of $\bc$ and $\hat \bc$.

Let $z_3\in H^3(K(\Z,3),\Z)\cong  {}^{(\bc,\hat \bc)}E_2^{0,3}$ be the canonical generator.
It also generates the group ${}^t E_2^{0,3}$ of the Leray-Serre
spectral sequence of the homotopy fibration $t:*\rightarrow K(\Z,4)$. A
part of its second page is 
$$\begin{array}{|c|c|c|c|c|c|}
\hline3&\Z&0&0&0&*\\\hline
2&0&0&0&0&0\\\hline
1&0&0&0&0&0\\\hline
0&\Z&0&0&0&\Z\\ \hline
q/p&0&1&2&3&4\\\hline \end{array}\ .$$
Since $H^3(*,\Z)\cong 0\cong H^4(*,\Z)$ we conclude that
${}^t d^{0,3}_2(z_3)=z_4\in H^4(K(\Z,4),\Z)\cong {}^t E_2^{4,0}$ is the generator.
Now by construction
$q^*z_4=x_1\cup \hat x_1+\dots+x_n\cup \hat x_n$,
and by naturality of the spectral sequences ${}^{(\bc,\hat \bc)}d_2^{0,3}(z)=q^*z_4$.
This implies the assertion about $H^4(\bR_n,\Z)$.\hB

\subsubsection{}
Recall that we consider $K(\Z^n,2)\cong BT^n$ (see  \ref{expla}).
\begin{ddd} We define  $\pi_n:\bF_n\rightarrow \bR_n$ to be the
  $T^n$-bundle which is classified by 
 $\bc:\bR_n\rightarrow K(\Z^n,2)$. 
\end{ddd}
  
Let $U\rightarrow K(\Z,2)$ be the universal $T^1$-bundle.
The $n$-fold product $U^n\rightarrow K(\Z,2)^n\cong K(\Z^n,2)$ is the
universal $T^n$-bundle. By definition we get a pull-back diagram
$$\begin{array}{ccc}
\bF_n&\rightarrow&U^n\\
\pi_n\downarrow&&\downarrow\\
\bR_n&\stackrel{\bc}{\rightarrow}&K(\Z,2)^n\end{array}\ .$$

\subsubsection{}\label{cd1}

In the following we use the spectral sequence notation introduced in \ref{lerayserreexplain} 
\begin{lem}\label{re1}
\begin{enumerate}
\item
We have 
$$
\begin{array}{|c|c|c|c|c|}  
\hline
i&0&1&2&3\\\hline
H^i(\bF_n,\Z)&\Z&0&\Z^{n}&\Z\\\hline \end{array}\ .
$$
\item
Here the group $H^2(\bF_n,\Z)$ is freely generated by the components
of $\pi_n^*\hat \bc$. 
\item In particular, restriction to the fiber of
$\bF_n\to \bR_n$ induces the zero homomorphism on
$H^2(\bF_n,\Z)$. 
\item Furthermore, $H^3(\bF_n,\Z)$ is generated by a class
$\bh_n\in \cF^2H^3(\bF_n,\Z)$ which is characterized by
$[\bh_n]^{2,1}=\sum_{i=1}^n[y_i\otimes \hat \bc_i] \in {}^{\pi_n}E_\infty ^{2,1}$. \end{enumerate}
\end{lem}
\proof
We write out the second page
${}^{\pi_n}E_2^{p,q}$. 
$$\begin{array}{|c|c|c|c|c|c|}
\hline3&\Z^{\frac{n(n-1)(n-2)}{6}}&&&&\\\hline
2&\Z^{\frac{n(n-1)}{2}}&0&\Z^{n^2(n-1)}&0&\\\hline
1&\Z^n&0&\Z^{2n^2}&0&\\\hline
0&\Z&0&\Z^{2n}&0&\Z^{{n(2n+1)}{}-1}\\ \hline
q/p&0&1&2&3&4\\\hline \end{array}\ .$$
We know that 
${}^{\pi_n}d_2^{0,1}(y_i)=\bc_i\in {}^{\pi_n}E_2^{2,0}\cong H^2(\bR_n,\Z)$.
%To this end we consider the universal
%$T^1$-bundle
%$q:U\rightarrow K(\Z,2)$.
%Let $y\in H^1(T^1,\Z)$ be the canonical generator.
%It also generates the group
%${}^q E_2^{0,1}$ of the associated Leray-Serre spectral sequence
%$$\begin{array}{|c|c|c|c|}\hline
%1&\Z&0&*\\\hline
%0&\Z&0&\Z\\ \hline
%q/p&0&1&2\\\hline \end{array}\ .$$
%Since $U$ is contractible ${}^q d_2^{0,1}(y)\in
%{}^q E_2^{2,0}\cong H^2(K(\Z,2),\Z)$ is a generator.
%The bundle
%$\pi_n:\bF_n\rightarrow \bR_n$ is obtained as a pull-back
%of an $n$-fold product of the bundle $U\rightarrow K(\Z,2)$ via the
%components
%of $\bc$.
%We conclude by the naturality of the spectral sequences that
%$d_2^{0,1}(y_i)\in\{\bc_i,-\bc_i\}$.
It follows that $ {}^{\pi_n} d_2^{0,1}$ is an isomorphism onto the subgroup of
${}^{\pi_n} E_2^{0,2}$ generated
by $\bc$ so that
${}^{\pi_n} E_3^{0,1}\cong 0$ and $ {}^{\pi_n}E_3^{2,0}$ is freely generated by the components
of $\hat \bc$. We see already that $H^1(\bF_n,\Z)\cong 0$.

The group ${}^{\pi_n} E_2^{0,2}$ is freely generated by all products
$y_i\cup y_j$, $i<j$. 
We now use the multiplicativity of the Leray-Serre spectral sequence
in order see that ${}^{\pi_n}d_2^{0,2}(y_i\cup y_j)=y_j\otimes \bc_i -y_i\otimes
\bc_j$. We conclude that ${}^{\pi_n}d_2^{0,2}$ is injective.
This implies that ${}^{\pi_n}E_3^{0,2}\cong 0$, and it follows that
$H^2(\bF_n,\Z)$ is freely generated by the components of
$\pi_n^*\hat \bc$.

We have ${}^{\pi_n}E_2^{4,0}\cong H^4(\bR_n,\Z)$, 
${}^{\pi_n}d^{1,2}_2(y_i\otimes \bc_j)=\bc_i\cup \bc_j$, and 
${}^{\pi_n}d^{1,2}_2(y_i\otimes \hat \bc_j)=\bc_i\cup \hat \bc_j$.

In order to  calculate $\ker({}^{\pi_n}d_2^{1,2})$ we recall the relation
$\bc_1\cup \hat \bc_1+\dots+\bc_n\cup
\hat \bc_n=0$. Let $$h:=y_1\otimes \hat
\bc_1+\dots+y_n\cup \hat \bc_n \ .$$ Then we have ${}^{\pi_n}d_2^{1,2}(h)=0$.
We claim that $\ker({}^{\pi_n}d_2^{1,2})\cong \Z h\oplus \im({}^{\pi_n}d_2^{0,2})$.
Let $t:=\sum_{i,j=1}^n a_{i,j} y_i\otimes \bc_j+ b_{i,j} y_i \otimes \hat
\bc_j$  for $a_{i,j},b_{i,j}\in \Z$ and assume that
${}^{\pi_n}d_2^{1,2}(t)=0$. Then
$\sum_{i,j=1}^n a_{i,j} \bc_i\cup\bc_j+ b_{i,j} \bc_i \cup \hat
\bc_j=0$.
This implies that $a_{i,j}+a_{j,i}=0$ for all pairs $i,j$, $b_{i,j}=0$ for $i\not=j$, and
that there exists $b\in \Z$ such that 
$b_{i,i}=b$ for all $i=1,\dots,n$.
But then we can write
$t=\sum_{i<j} a_{i,j} {}^{\pi_n}d_{2}^{0,2}(y_j\cup y_i) + bh$.
It follows that
${}^{\pi_n}E_3^{2,1}\cong \Z$ is generated by the class of $h$.

The group ${}^{\pi_n}E_2^{3,0}$ is freely generated by
the products $y_i\cup y_j\cup y_k$, $i<j<k$. Furthermore, 
${}^{\pi_n}d_2^{0,3}(y_i\cup y_j\cup y_k)=y_j\cup y_k\otimes \bc_i - y_i\cup
y_k\otimes \bc_j + y_i\cup y_j\otimes \bc_k$.
We thus calculate that ${}^{\pi_n}d_2^{0,3}$ is injective.
We conclude that $H^3(\bF_n,\Z)\cong \cF^2 H^3(\bF_n,\Z)$ is generated
by the class $\bh_n$ represented by $h\in {}^{\pi_n} E_2^{1,2}$.
\hB 
 
\subsubsection{}\label{f16}

In the proof of Lemma \ref{hh} we have found a homotopy Cartesian
square
$$\begin{array}{ccc}K(\Z,3)&\stackrel{i}{\rightarrow}&\bR_n\\
\downarrow&&(\bc,\hat \bc)\downarrow\\
*&\rightarrow&K(\Z^n,2)\times K(\Z^n,2)\end{array}\ .$$
We pull this square back along the map
$\psi:K(\Z^n,2)\cong U^n\times K(\Z^n,2)\rightarrow K(\Z^n,2)\times K(\Z^n,2)$
and obtain a cube of homotopy Cartesian squares
$$\begin{array}{ccccccc}
i^*\bF_n&&&\stackrel{I}{\longrightarrow}&&&\bF_n\\&\pr\searrow&&&&\pi_n\swarrow&\\
&&K(\Z,3)&\stackrel{i}{\rightarrow}&\bR_n&&\\\lambda \downarrow&&\downarrow&&(\bc,\hat
\bc)\downarrow&&\kappa\downarrow\\
&&*&\rightarrow&K(\Z^n,2)\times K(\Z^n,2)&&\\&\nearrow&&&&\psi\nwarrow&\\K(\Z,1)^n&&&\longrightarrow&&&K(\Z^n,2)
\end{array}\ .$$
\begin{lem}\label{g16}
We have $I^*\bh_n=\pm \pr^*z_3$,
where $z_3\in H^3(K(\Z,3),\Z)$ is the canonical generator.
\end{lem}
\proof
The fiber of $\kappa$ is equivalent to $K(\Z,3)$.
The second page ${}^\kappa E_2$ of the corresponding  Leray-Serre spectral
sequence 
has the form
$$\begin{array}{|c|c|c|c|c|c|}
\hline3&\Z&0&*&*&*\\\hline
2&0&0&0&0&*\\\hline
1&0&0&0&0&*\\\hline
0&\Z&0&\Z^{n}&0&\Z^{\frac{n(n+1)}{2}}\\ \hline
q/p&0&1&2&3&4\\\hline \end{array}\ .$$
We know that $H^3(\bF_n,\Z)\cong \Z$ is freely generated by $\bh_n$.
We see that ${}^\kappa d_2^{0,3}=0$ and ${}^\kappa E_2^{0,3}\cong H^3(\bR_n,\Z)$
is generated by $\bh_n$. On the other hand
the group ${}^\kappa E_{2}^{0,3}\cong H^3(K(\Z,3),\Z)$ is freely generated
by $z_3$. Therefore $\bh_n=\pm z_3$.

Note that $i^*\bF_n\cong K(\Z,3)\times K(\Z,1)^n$ so that
the Leray-Serre spectral sequence ${}^\lambda E$ of $\lambda$ degenerates.
Let $I^*:{}^\kappa E_2\rightarrow {}^\lambda E_2$ be
the induced map of the second pages  and note that
${}^\lambda E_2$ has the form
$$\begin{array}{|c|c|c|c|c|c|}
\hline3&\Z&\Z^n&*&*&*\\\hline
2&0&0&0&0&0\\\hline
1&0&0&0&0&0\\\hline
0&\Z&\Z^n&\Z^{\frac{n(n-1)}{2}}&\Z^{\frac{n(n-1)(n-2)}{6}}&\Z^{\frac{n(n-1)(n-2)(n-3)}{24}}\\ \hline
q/p&0&1&2&3&4\\\hline \end{array}\ ,$$
where ${}^\lambda E_2^{0,3}$ is freely generated by $z_3\in
H^3(K(\Z,3),\Z)$. 
The map $I$ induces an equivalence of the fibers of $\lambda$ and
$\kappa$. In particular, it induces an isomorphism 
$I^*:{}^\kappa E_2^{0,3}\rightarrow {}^\lambda E_2^{0,3}$ identifying the generators
above. This implies that $I^*\bh_n=\pm \pr^*z_3$. \hB

\section{The $T$-duality group and the universal triple}

\subsubsection{}\label{tz}

Let $(e_1,\dots,e_n,\hat e_1,\dots,\hat e_n)$ be the standard basis of $\Z^{2n}$.
Let $G_n\subset GL(2n,\Z)$ be the subgroup of transformations which
fix the form $q:\Z^{2n}\rightarrow \Z$ given by 
$q(\sum_{i=1}^n a_i e_i+b_i\hat e_i):=\sum_{i=1}^n a_ib_i$.
Usually denoted $O(n,n,\Z)$, $G_n$ will here be
called the group of $T$-duality transformations.

\subsubsection{}

Each $g\in G_n$ induces an equivalence
$g:K(\Z^{n},2)\times K(\Z^{n},2)\rightarrow K(\Z^{n},2)\times K(\Z^{n},2)$.
\begin{lem}\label{re}
There exist a unique homotopy class of lifts $\tilde g$ in the diagram
$$
\begin{CD}
  \bR_n @>{\tilde g}>> \bR_n\\
  @VV{(\bc,\hat \bc)}V @VV{(\bc,\hat \bc)}V\\
  K(\Z^{n},2)\times K(\Z^{n},2) @>{g}>> K(\Z^{n},2)\times K(\Z^{n},2)
\end{CD}
%\begin{array}{ccc}& &\bR_n\\
%&\tilde g \nearrow&\downarrow\\
%\bR_n&\stackrel{g\circ (\bc,\hat \bc)}{\rightarrow}&K(\Z,2)^{2n} \end{array}
\ .$$
\end{lem}
\proof
We apply obstruction theory to the problem of existence and
classification of lifts $\tilde g$. The situation is simple because the fibre of 
$(\bc,\hat \bc)$ is a $K(\Z,3)$.
In fact, the obstruction against the existence of a lift is the class $(\bc,\hat \bc)^*g^*
(\sum_{i=1}^n x_i\cup \hat x_i)\in
H^4(\bR_n,\Z)$
which vanishes since $g$ preserves $q$. Therefore a lift $\tilde g$ exist.
The set of homotopy classes of lifts is a torsor
over $H^3(\bR_n,\Z)$. Since $H^3(\bR_n,\Z)\cong  \{0\}$, the lift $\tilde g$ is unique.
\hB 

The correspondence $G_n\ni g\mapsto \tilde g$ induces a homotopy action of $G_n$
on $\bR_n$.

\subsubsection{}

We let $t\in G_n$ be the transformation
given by $t(e_i)=\hat e_i$ and $t(\hat e_i)=e_i$.
\begin{ddd}\label{tt1}
The universal $T$-duality is the lift
$T:=\tilde t:\bR_n\rightarrow \bR_n$ of $t$ according to Lemma
\ref{re}.
\end{ddd}
Note that $T \circ T=\id_{\bR_n}$ since $t^2=1\in G_n$.

\subsubsection{}

\begin{ddd}\label{rtz}
The universal dual $T^n$-bundle is defined by the pull-back
 $$\begin{array}{ccc}
\hat \bF_n&\stackrel{\tilde T}{\rightarrow}& \bF_n\\
\hat \pi_n\downarrow&&\pi_n\downarrow\\
\bR_n&\stackrel{T}{\rightarrow}&\bR_n \end{array}
\ .$$
Furthermore, we define $\hat \bh_n:=\tilde T^*\bh_n\in H^3(\hat \bF_n,\Z)$.
\end{ddd}

\newcommand{\br}{{\bf r}}
\subsubsection{}\label{adef}

We consider the pull-back diagram
$$\begin{array}{ccccc}
&&\bF_n\times_{\bR_n} \hat \bF_n&&\\
&\bp_n\swarrow&&\hat \bp_n\searrow&\\
\bF_n&&\br_n\downarrow&&\hat \bF_n\\
&\pi_n\searrow&&\hat \pi_n\swarrow&\\
&&\bR_n&&
\end{array}\ .$$

\begin{lem}\label{aeq}
We have 
$$\begin{array}{|c|c|c|c|c|}\hline
i&0&1&2&3\\\hline
H^i(\bF_n\times_{\bR_n} \hat \bF_{n},\Z)&\Z&0&0&\Z\\\hline\end{array}
\ .$$
Moreover, $\pi_n^* \bh_n=\hat \pi_n^* \hat \bh_n$, and this element
generates $H^3(\bF_n\times_{\bR_n} \hat \bF_n,\Z)$.
\end{lem}
\proof
We use the Leray-Serre spectral sequence ${}^{\br_n} E$.
 The relevant part of its second page has the form
$$\begin{array}{|c|c|c|c|c|c|}
\hline3&\Z^{\frac{2n(2n-1)(2n-2)}{6}}&&&&\\\hline
2&\Z^{\frac{2n(2n-1)}{2}}&0&\Z^{2n^2(2n-1)}&0&\\\hline
1&\Z^{2n}&0&\Z^{4n^2}&0&\\\hline
0&\Z&0&\Z^{2n}&0&\Z^{{n(2n+1)}-1}\\ \hline
q/p&0&1&2&3&4\\\hline \end{array}\ .$$
Now ${}^r E_2^{0,1}$ is freely generated by $y_i,\hat y_i$, and ${}^r E_2^{2,0}$
is freely generated by $\bc_i,\hat \bc_i$. We know that
${}^r d_2^{0,1}(y_i)=\bc_i$ and ${}^r d_2^{0,1}(\hat y_i)=\hat \bc_i$.
We conclude that ${}^r d_2^{0,1}$ is injective and $H^1(\bF_n\times_{\bR_n}
\hat \bF_{n},\Z)\cong 0$.

The group ${}^{\br_n} E_2^{0,2}$ is freely generated by  the products $y_i\cup
y_j$, $\hat y_i\cup \hat y_j$, $i<j$ and
$y_i\cup \hat y_j$.
We have ${}^{\br_n} d_2^{0,2}(y_i\cup y_j)=y_j\otimes \bc_i- y_i\otimes \bc_j$,
$ {}^{\br_n} d_2^{0,2}(\hat y_i\cup \hat y_j)=\hat y_j\otimes \hat\bc_i-\hat y_i\otimes
\hat \bc_j$, and
${}^{\br_n} d_2^{0,2}(y_i\cup \hat y_j)=\hat y_j\otimes  \bc_i-y_i\otimes
\hat \bc_j$. It follows that ${}^r d_2^{0,2}$ is injective.

In a similar way we see that ${}^{\br_n} d_2^{0,3}$ is injective.

We now calculate ${}^{\br_n} E_3^{2,1}$. We claim that
this group is freely generated by one class
which can be represented by
$\sum_{i=1}^n y_i\otimes \hat \bc_i\in {}^{\br_n} E^{2,1}_2 $ or alternatively by $\sum_{i=1}^n \hat y_i\otimes
\bc_i\in {}^{\br_n} E^{2,1}_2$.
In view of the construction of $\bh_n$ and $\hat \bh_n$, and by
naturality of the Leray-Serre spectral sequence this would
imply the assertion of the lemma about the third cohomology.

The group ${}^{\br_n} E_2^{4,0}$ is generated by the products  $\bc_i\cup \bc_j, \hat \bc_i\cup
\hat \bc_j$, $i\le j$,  and all products  $\bc_i\cup \hat \bc_j$ subject to one
relation
$\sum_{i=1}^n \bc_i\cup \hat \bc_i=0$.
Let $t=\sum_{i,j}a_{i,j} y_i\otimes \bc_j+\sum_{i,j}b_{i,j} y_i\otimes
\hat \bc_j+\sum_{i,j}c_{i,j} \hat y_i\otimes \bc_j+\sum_{i,j}d_{i,j}
\hat y_i\otimes \hat \bc_j$ and assume that ${}^{\br_n} d_2^{2,1}(t)=0$.
Then we have $\sum_{i,j}a_{i,j} \bc_i\otimes \bc_j+\sum_{i,j}b_{i,j} \bc_i\otimes
\hat \bc_j+\sum_{i,j}c_{i,j} \hat \bc_i\otimes \bc_j+\sum_{i,j}d_{i,j}
\hat \bc_i\otimes \hat \bc_j=0$.
This implies that $a_{i,j}+a_{j,i}=0$, $d_{i,j}+d_{j,i}=0$, for all
$i,j$, 
$b_{i,j}+ c_{i,j}=0$ for all $i\not=j$,
and that there exists a unique $e\in \Z$ such that $b_{i,i}+c_{i,i}=e$ for all
$i=1,\dots,n$.
We can now write
$t=\sum_{i<j} a_{i,j}{}^{\br_n}d_{2}^{0,2}(y_j\cup y_i)+\sum_{i<j}
d_{i,j}{}^{\br_n} d_{2}^{0,2}(\hat y_j\cup \hat y_i)
- \sum_{i\not=j} b_{i,j}  {}^{\br_n} d_2^{0,2}(y_i\cup \hat y_j)-
\sum_{i=1}^n b_{i,i}{}^{\br_n} d_{2}^{0,2} (y_i\cup \hat y_i) +
e \sum_{i=1}^n \hat y_i\otimes \bc_i$.
This already shows that
${}^{\br_n} E_{3}^{2,1}$ is freely generated by the class of
$\sum_{i=1}^n \hat y_i\otimes \bc_i$.
Finally note that
$\sum_{i=1}^n {}^{\br_n} d_2^{0,2}(y_i\cup \hat y_i)=\sum_{i=1}^n \hat y_i\otimes
\bc_i-
\sum_{i=1}^n y_i\otimes \hat \bc_i$.
This finishes the proof of the claim.
\hB 

\subsubsection{}
 
Let $\cH_n\in T(\bF_n)$  be a twist with isomorphism class
$[\cH_n]=\bh_n\in H^3(\bF_n,\Z)$. Set further $\hat \cH_n:=\tilde T^*\cH_n \in T(\hat
\bF_n)$,
where $\tilde T$ was defined in \ref{rtz}.
Then $[\hat \cH_n]=\hat \bh_n$.

Recall the definition of $\bp_n$ and $\hat \bp_n$ in \ref{adef}.
Since $\bp_n^*\bh_n=\hat  \bp_n^*\hat \bh_n$ by Lemma \ref{aeq}, we conclude that
there exists an isomorphism of twists $\bu_n:\hat \bp_n^*\hat
\cH_n\rightarrow \bp_n^* \cH_n$. Since  $H^2(\bF_n\times_{\bR_n}\hat
\bF_n,\Z)\cong 0$ by Lemma \ref{aeq} this isomorphism is unique.

\subsubsection{}

The notion of a $T$-duality triple over a space $B$ was
introduced in Definition \ref{tr13}. Here we clarify the notion of an
isomorphism between such triples  $x:=((F,\cH),(\hat F,\hat
\cH),u)$ and $x^\prime:=((F^\prime,\cH^\prime),(\hat F^\prime,\hat
\cH^\prime),u^\prime)$. 
\begin{ddd}\label{isodd}
We say that $x$ and $x^\prime$ are isomorphic if there exist
underlying 
bundle isomorphisms
$$\begin{array}{ccc}
F&\stackrel{\psi}{\rightarrow}&F^\prime\\
\downarrow&&\downarrow\\
B&=&B\end{array}\ ,\quad \begin{array}{ccc}
\hat F&\stackrel{\hat\psi}{\rightarrow}&\hat F^\prime\\
\downarrow&&\downarrow\\
B&=&B\end{array}$$
and isomorphisms of twists
$$v:\psi^*\cH^\prime\rightarrow \cH\ ,\quad \hat v:\hat
\psi^*\hat\cH^\prime\rightarrow \hat \cH$$
such that the composition
$$\hat p^*\hat \cH\stackrel{\hat p^*\hat v^{-1}}{\rightarrow} \hat p^* \hat
\psi^*\hat \cH^\prime\cong (\psi,\hat\psi)^*(\hat p^\prime)^*\hat
\cH^\prime
\stackrel{(\psi,\hat \psi)^*u^\prime}{\rightarrow}(\psi,\hat\psi)^*(p^\prime)^*
\cH^\prime\cong p^*\psi^* \cH^\prime\stackrel{p^* v}{\rightarrow}
p^*\cH$$
is equal to $u:\hat p^*\hat \cH\rightarrow p^*\cH$.
Here $(\psi,\hat \psi):F\times_B\hat F\rightarrow F^\prime\times_B\hat
F^\prime$ is the induced map (and the other notation is as in \ref{tr56512u}). 
\end{ddd}

\subsubsection{}

\begin{theorem}\label{ex980}
$((\bF_n,\cH_n),(\hat \bF_n,\hat \cH_n),\bu_n)$ represents the unique
isomorphism class of 
$T$-duality triples over $(\bF_n,\hat \bF_n)$.
\end{theorem}
\proof
We must only
verify that $((\bF_n,\cH_n),(\hat \bF_n,\hat \cH_n),\bu_n)$ is a
$T$-duality triple, because uniqueness is clear from the conditions on
$T$-duality triples of \ref{tr56512u}, from uniqueness of
$[\cH_n]$ and $[\hat \cH_n]$ by Lemma \ref{re1} (1),  and the uniqueness of
$\bu_n$ by Lemma \ref{aeq}.

The classes $[\cH]$ and $[\hat \cH]$ have the required properties by
construction. It remains to show that condition $\cP(u)$ is satisfied
(see \ref{cpudef}). We fix a base point $*\in \bR_n$. 
Note that by Lemma \ref{var} and Proposition \ref{var1} (which are
independent of the result to be proved)
we have a canonical equivalence $\bR_1\cong \tilde \bR_1$.
In \cite{bunkeschickt} we studied in detail the topology
of $\bR_1$ and the associated $T$-duality. The idea of the proof is
to reduce the present task to the case $n=1$.

We consider the $n$-fold product
$F\rightarrow B$ of the $T^1$-bundle $\bF_1\rightarrow \bR_1$
i.e. we set $F:=\bF_1^n$ and  $B:=\bR_1^n$.
Let $p_i:B\rightarrow \bR_1$, $i=1,\dots,n$, denote the projections.
Let $(z,\hat z):B\rightarrow K(\Z^n,2)\times K(\Z^n,2) $ be the map
whose  components classify $z_i:=p_i^*x$ and $\hat z_i:=p_i^*\hat
x$, and where $x,\hat x\in H^2(\bR_1,\Z)$ are the canonical generators.

We now apply obstruction theory to the lifting problem
$$\begin{array}{ccc}
&&\bR_n\\
&f\nearrow&(\bc,\hat\bc)\downarrow\\
B&\stackrel{(z,\hat z)}{\rightarrow}&K(\Z^n,2)\times K(\Z^n,2) 
\end{array}\ .$$
The relation $x\cup \hat x=0$ in $H^4(\bR_1,\Z)$ implies that$z_i\cup \hat z_i=0$, and hence that
$\sum_{i=1}^n z_i\cup \hat z_i=0$. In view of the definition \ref{rt3}
this diagram admits a lift $f$.
 Since
$H^3(B,\Z)\cong 0$ the homotopy class of lifts of the lift  $f$ is in fact uniquely determined.

We therefore have a pull-back diagram of principal $T^n$-bundles
$$\begin{array}{ccc}
F&\stackrel{\tilde f}{\rightarrow}&\bF_n\\
\pi\downarrow&&\pi_n\downarrow\\
B&\stackrel{f}{\rightarrow}&\bR_n\end{array}\ .$$
We define $h:=\tilde f^*\bh_n$.
In this way we obtain a pair
$(F,h)$ over $B$.
We further have 
natural projections $\pr_i:F\rightarrow \bF_1$, $i=1,\dots,n$.
Using the characterizations of $\bh_n$ and $\bh_1$
in Lemma \ref{re1}, the naturality of Leray-Serre spectral
sequences, and $H^3(B,\Z)\cong 0$, we see that 
\begin{equation}\label{defc}h=\sum_{i=1}^n \pr^*_i \bh_1\ .\end{equation}

Let $T_B:B\rightarrow B$ be the product of the $T$-duality
transformations $T:\bR_1\rightarrow \bR_1$ on each factor. Since $T$ is a homotopy
equivalence, there is a unique homotopy classes of lifts $\alpha$ in
$$\begin{array}{ccc}
&&\bR_n\\
&\alpha\nearrow&T\downarrow\\
B&\stackrel{f\circ T_B}{\rightarrow}&\bR_n 
\end{array}\ .$$

Since the composition $B\stackrel{\alpha}{\rightarrow}\bR_n\rightarrow
K(\Z^n,2)\times K(\Z^n,2)$ coincides with $(z,\hat z)$ we have $\alpha=f$, and the
following diagram commutes upto homotopy
$$\begin{array}{ccc}
B&\stackrel{f}{\rightarrow}&\bR_n\\
T_B\downarrow&&T\downarrow\\
B&\stackrel{f}{\rightarrow}&\bR_n\end{array}\ .$$
This shows that we have a pull-back diagram
$$\begin{array}{ccc}
\hat F&\stackrel{{\hat{f}}}{\rightarrow}&\hat \bF_n\\
\hat \pi \downarrow&&\hat \pi_n\downarrow\\
B&\stackrel{f}{\rightarrow}&\bR_n\end{array}\ ,$$
where $\hat F\rightarrow B$ is the $n$-fold product of
$\hat \bF_1\rightarrow \bR_1$.

We get the commutative diagram
$$\begin{array}{ccccc}
F_b\times \hat F_b&\rightarrow&F\times_B\hat F&\stackrel{\tilde f\times_B \hat f}{\rightarrow}&\bF_n\times_{\bR_n}\hat \bF_n\\
\downarrow&&r\downarrow&&\br_n\downarrow\\
\{b\}&\rightarrow&B&\stackrel{f}{\rightarrow}&\bR_n\end{array}\ ,$$
where we assume that $f(b)=*$.

It follows from Lemma \ref{aeq} (in the case $n=1$) and the K{\"u}nneth
formula that
$H^2(F\times_B \hat F,\Z)\cong 0$. Therefore 
$(\tilde f\times_B \hat f)^*u:(\tilde f\times_B \hat f)^*\hat
\bp_n^*\hat \cH\rightarrow (\tilde f\times_B \hat f)^*\bp_n^*\cH$ is the unique isomorphism.

Let $\cV\in T(\bF_1)$ and $\hat \cV\in T(\hat \bF_1)$ be twists in the
classes
$\bh_1$ and $\hat \bh_1$.
Then by (\ref{defc}) there exists isomorphisms of twists $\kappa:\tilde
f^*\cH\rightarrow  \sum_{i=1}^n\pr_i^*\cV$ and
$\hat \kappa:\hat{f}^*\hat \cH\rightarrow  \sum_{i=1}^n \pr_i^*\hat \cV$.
The choice of these isomorphisms is not unique. However, their
restrictions  to the fiber over $\{b\}\rightarrow B$ is unique. This follows
from the structure of $H^2(F,\Z)$ (and of $H^2(\hat F,\Z)$) implied by Lemma \ref{re1} (1)
in the case $n=1$ and the K{\"u}nneth formula.

At this moment we fix some choices of $\kappa$ and $\hat\kappa$.

Let $q_i:F\times_B\hat F\rightarrow \bF_1\times_{\bR_1} \hat \bF_1$
be the projection onto the $i$th component.
Note that $\tilde f\circ p=\bp_n\circ  (\tilde f\times_B\hat f)$,
${\hat {f}}\circ \hat p=\hat \bp_n\circ (\tilde f\times_B\hat f)$, 
$\bp_1\circ q_i=\pr_i\circ p$ and $\hat \bp_1\circ q_i=\pr_i\circ \hat p$.
We now have fixed isomorphisms
\begin{equation*}
p^* \kappa\colon (\tilde f\times_B\hat f)^*\bp_n^*\cH\cong \sum_{i=1}^n
q_i^* \bp_1^* \cV\ ,\qquad\hat p^*\hat\kappa: (\tilde
f\times_B\hat f)^*\hat
\bp_n^*\hat \cH\cong \sum_{i=1}^n q_i^* 
\hat \bp_1^* \hat \cV.
\end{equation*}
Note that there is a unique isomorphism
$\psi:\hat \bp_1^*\hat \cV\rightarrow \bp_1^*\cV$.
This induces another isomorphism
$$\Phi:(\tilde f\times_B\hat f)^*\hat \bp_n^*\hat \cH\stackrel{\hat
  p^*\hat \kappa}{\cong} \sum_{i=1}^n q_i^* \hat \bp_1^* \hat
\cV\stackrel{\sum_{i=1}^n q_i^* \psi }{\cong}
 \sum_{i=1}^n q_i^* \bp_1^*
\cV\stackrel{(p^*\kappa)^{-1}}{\cong} (\tilde f\times_B\hat f)^*\bp_n^*\cH\ .
$$
It follows that
$(\tilde f\times_B\hat f)^*u= \Phi$.  We can now restrict
$\Phi$ to the fiber $F_b\times \hat F_b$.

It was shown in \cite[3.2.4]{bunkeschickt} that the restriction of
$\psi$ to the fiber $T^1\times \hat T^1$ is
classified by a generator of $H^2(T^1\times \hat T^1,\Z)$, namely
by $y\cup \hat y$ in the canonical basis of $H^1(T^1\times \hat
T^1,\Z)$. If we restrict the whole composition defining $\Phi$ to the
fiber the we see that $[(\tilde f,\hat f)^*u(b)]=[\sum_{i=1}^n y_i\cup \hat y_i]\in H^2(F_b\times \hat F_b,\Z)/(\im(p_b^*)+\im(\hat p_b^*))$, as required.

\hB

\section{Pairs and triples}\label{onecon}

\subsubsection{}

Recall the construction of the classifying space $\tilde \bR_n(0)$ of
pairs \ref{fff564}. The goal of the present section is the
identification of $\bR_n$
with a one-connected covering of $\tilde \bR_n(0)$.

We start with the universal $T^n$-bundle $U^n\rightarrow K(\Z^n,2)$
and the $T^n$-space $\Map(T^n,K(\Z,3))$, where $T^n$ acts by
reparametrization. 
In a first step we form
the associated bundle
$$p:\tilde \bR_n(0):=U^n\times_{T^n} \Map(T^n,K(\Z,3))\rightarrow K(\Z^n,2)\ .$$
We define a $T^n$-bundle $\tilde\bF_n()$ via pull-back
$$\begin{array}{ccc}
\tilde \bF_n(0)&\rightarrow&U^n\\
\downarrow&&\downarrow\\
\tilde \bR_n(0)&\stackrel{p}{\rightarrow}&K(\Z^n,2)\end{array}
\ .$$
There is a canonical map
$$\tilde \bh_n(0)\colon \tilde \bF_n(0)\rightarrow K(\Z,3)\ .$$
It is
given by
$\tilde \bh_n(0)([v,\phi],u):=\phi(s)$, where $s\in T^n$ is the unique
element such that
$sv=u$. Here $u,v\in U^n$, $\phi\in \Map(T^n,K(\Z,3))$, 
$[v,\phi]\in \tilde \bR_n(0)$, and $([v,\phi],u)\in \tilde \bF_n(0)$.

\subsubsection{}

Recall that a pair $(F,h)$ over a space $B$ consists of a $T^n$-bundle
$F\rightarrow B$ and a class $h\in H^3(F,\Z)$.
An isomorphism between pairs $(F,h)$ and $(F^\prime,h^\prime)$
is given by a diagram
$$\begin{array}{ccc}
F&\stackrel{\Phi}{\rightarrow}&F^\prime\\ 
\downarrow&&\downarrow\\
B&=&B
\end{array}\ ,$$
where $\Phi$ is a $T^n$-bundle isomorphism
such that $\Phi^*h^\prime=h$.

Given a map $f:B^\prime\rightarrow B$ of spaces we can form the
pull-back
$$\begin{array}{ccc}
F^\prime&\stackrel{\tilde f}{\rightarrow}&F\\ 
\downarrow&&\downarrow\\
B^\prime&\stackrel{f}{\rightarrow}&B
\end{array}\ .$$
We define the pair
$f^*(F,h):=(F^\prime,\tilde f^*h)$ over $B^\prime$. Pull-back preserves
isomorphism classes of pairs.

\subsubsection{}

Let $\tilde P_{(0)}$ be the contravariant set-valued functor which associates to
each space $B$ the set $\tilde P_{(0)}(B)$ of isomorphism classes of  pairs.

\begin{lem}\label{cl1}
The space $\tilde \bR_n(0)$ is a classifying space for $\tilde P_{(0)}$. More precisely,
the pair $[\tilde \bF_n(0),\tilde \bh_n(0)]\in P_{(0)}(\tilde \bR_n(0))$ induces a natural
isomorphism $\tilde v_B:[B,\tilde \bR_n(0)]\rightarrow  \tilde
P_{(0)}(B)$ such that $\tilde v_B(f)=f^* [\tilde \bF_n(0),\tilde \bh_n(0)]$.
\end{lem}
\proof
This is completely analogous to the proof of \cite[Proposition
2.6]{bunkeschickt}. We therefore refrain from repeating the proof here.
\hB

\subsubsection{}\label{ghasjhasd}

By Lemma \ref{cl1} we have an isomorphism 
$\pi_0(\tilde \bR_n(0))\cong  \tilde P_{(0)}(*)$.
Furthermore, note that
$P_{(0)}(*)\cong H^3(T^n,\Z)$ in a canonical way.
We define
$\tilde \bR_n(1)\subset \tilde \bR_n(0)$ to be the component
which corresponds to $0\in H^3(T^n,\Z)$.
Restricting the pair $[\tilde \bF_n(0),\tilde \bH_n(0)]$
gives a pair $[\tilde \bF_n(1),\tilde \bh_n(1)]$ over $\tilde \bR_n(1)$.
We let $\tilde P_{(1)}$ be the functor classified by $\tilde \bR_n(1)$.
Observe that $\tilde P_{(1)}(B)\subset \tilde P_{(0)}(B)$ is the set of isomorphism classes
of pairs $[F,h]$ such that the restriction of $h$ to the fibers
of $F$ vanishes.

\subsubsection{}\label{cl2}

By Lemma \ref{cl1} we have $\tilde P_{(0)}(S^1)\cong H^3(S^1\times T^n)$.
By the K{\"u}nneth formula
$ H^3(S^1\times T^n,\Z)\cong H^3(T^n,\Z)\oplus H^2(T^n,\Z)$, and 
$\pi_1(\tilde \bR_n(1))\cong \tilde P_{(1)}(S^1)\cong
H^2(T^n,\Z)$ corresponds to the second summand. On can check that this
bijection is a group homomorphism.

We consider the
isomorphism
$\phi:\pi_1(\tilde \bR_n(1))\stackrel{\sim}{\rightarrow} H^2(T^n,\Z)$ as a cohomology class
$\phi\in H^1(\tilde \bR_n(1),H^2(T^n,\Z)))$, i.e. as a homotopy class of maps
$\phi:\tilde \bR_n(1)\rightarrow K(H^2(T^n,\Z),1)$.
We define $\tilde \bR_n$ as the homotopy pullback
$$ \begin{array}{ccc}
\tilde \bR_n&\rightarrow&\tilde \bR_n(1)\\
\downarrow&&\phi\downarrow\\
*&\rightarrow&K(H^2(T^n,\Z),1)\end{array}\ .$$
Furthermore, we consider the pull-back
$$\begin{array}{ccc}\tilde \bF_n&\rightarrow&\tilde \bF_n(1)\\
\downarrow&&\downarrow\\
\tilde \bR_n&\rightarrow &\tilde \bR_n(1)\end{array}\ ,
$$
and we let $\tilde \bh_n\in H^3(\tilde \bF_n,\Z)$ be the pullback of $\tilde \bh_n(1)$.
Note that by construction and naturality,
$\tilde \bh_n$ pulls back to zero on the fiber of $\tilde \bF_n\to
\tilde \bR_n$. Since $\tilde \bR_n$ is simply connected, it then even
belongs to the second step $\cF^2 H^3(\tilde
\bF_n,\Z)$ of the Leray-Serre filtration associated with
$\tilde
\bF_n\rightarrow \tilde
\bR_n$. Let $x_1,\dots,x_n\in H^2(\tilde\bR_n,\Z)$ be the Chern classes of
$\tilde \bF_n$.

\subsubsection{}

\begin{lem}\label{ght}
The homotopy groups of $\tilde \bR_n$ are given by
$$\begin{array}{|c|c|c|c|c|c|}\hline
i&0&1&2&3&\ge 4\\\hline
\pi_i(\tilde \bR_n)&*&0&\Z^{2n}&\Z&0\\\hline\end{array}
\ .$$\end{lem}
\proof
By construction, $\tilde \bR_n$ is connected
and simply connected.
The homotopy fiber
of $\tilde \bR_n\rightarrow \tilde \bR_n(1)$ is equivalent to 
the homotopy fiber of
$*\rightarrow K(H^2(T^n,\Z),1)$, i.e.~to $K(H^2(T^n,\Z),0)$.
Hence 
this map induces an isomorphism
$\pi_i(\tilde \bR_n)\cong \pi_i(\tilde \bR_n(1))$ for $i\ge 2$.

Note that $K(\Z,k)$ is an $H$-space for each $k$. Hence have an equivalence
$$\Map(T^1,K(\Z,k))\simeq K(\Z,k)\times \Omega K(\Z,k)\simeq
K(\Z,k)\times K(\Z,k-1)\ .$$ We use the exponential law to write
$\Map(T^n,K(\Z,3))$ as an iterated mapping space, and we obtain in this way an equivalence
$$\Map(T^n,K(\Z,3))\simeq K(\Z,3)\times K(\Z,2)^{n}\times
K(\Z,1)^{\frac{n(n-1)}{2}}\times K(\Z,0)^{\frac{n(n-1)(n-2)}{6}}\ .$$
The long exact sequence of homotopy groups 
for 
$$\Map(T^n,K(\Z,3))\rightarrow \tilde \bR_n(0)\rightarrow K(\Z,2)^n$$
and the fact that $\pi_3(K(\Z,2)^n)\cong 0\cong \pi_1(K(\Z,2)^n)$ and
$\pi_2(K(\Z,2)^n)\cong  \Z^n$ yields the exact sequence
$$0\rightarrow \Z^n\rightarrow \pi_2(\tilde \bR_n(1))\rightarrow
\Z^n\stackrel{\delta}{\rightarrow} \Z^{\frac{n(n-1)}{2}}\stackrel{\alpha}{\rightarrow} 
\pi_1(\tilde \bR_n(1))\rightarrow 0\ .$$

Furthermore, we observe that $\pi_3(\tilde \bR(1))\cong \Z$ and $\pi_i(\tilde \bR_n(1))=0$ for $i\ge 4$.

We have seen in \ref{cl2} that $\pi_1(\tilde\bR_n(1))\cong H^2(T^2,\Z)\cong
\Z^{\frac{n(n-1)}{2}}$. We conclude that $\alpha$ must be
surjective. Consequently it is injective and $\delta=0$.

We therefore have an exact sequence
$$0\rightarrow \Z^n\rightarrow \pi_2(\tilde \bR_n(1))\rightarrow
\Z^n\rightarrow 0\ ,$$
and this implies that
$\pi_2(\tilde \bR_n(1))\cong  \pi_2(\tilde \bR_n) \cong \Z^{2n}$.
\hB 

By duality, the classes $x_1,\dots,x_n\in H^2(\tilde\bR_n,\Z)$ introduced
above form half of a $\Z$-basis for $H^2(\tilde\bR_n,\Z)\cong \Z^{2n}$.

\subsubsection{}\label{bc1}
 
Since the fiber of $p\colon\tilde \bR_n\rightarrow \tilde \bR_n(1)$ is
equivalent to $K(H^2(T^n,\Z),0)\cong \pi_1(\tilde \bR_n(1))$
we can consider this map
as a universal covering of $\tilde \bR_n(1)$.
We now consider the problem of existence and classification
of lifts in the diagram
$$\begin{array}{ccccc}
&&\tilde \bR_n & \rightarrow &*\\
&\tilde f\nearrow&p\downarrow && \downarrow\\
B&\stackrel{f}{\rightarrow}&\tilde \bR_n(1) & \xrightarrow{\phi} & K(H^2(T^n,\Z),1)\end{array}\ .
$$ 
It follows from the construction of $\tilde \bR_n$ that a lift $\tilde
f$ exists if and only if $\phi\circ f$ is homotopic to a constant map. The lift
itself depends on the choice of an explicit homotopy. 
If a lift exists, the  set of homotopy classes of lifts is a torsor
over $H^0(B,H^2(T^n,\Z))$. 

\subsubsection{}\label{clh}

The classification of homotopy classes $\tilde f$ (considered just as
maps, not as lifts)
lifting a homotopy class $f$ is more subtle. 
In order to study this problem we assume that $B$ is path connected
and equipped with a base point $b\in B$.
Let $\tilde f_0$ be a lift of $f$ and consider
$x\in \pi_1(\tilde \bR_n(1))\cong H^0(B,H^2(T^n,\Z))$.
Then we consider the lift $\tilde f_1=x\tilde f_0$, i.e.~the
composition of $\tilde f_0$ with the deck transformation associated to
$x$.

Assume that $\tilde f_0$ and $\tilde f_1$ are homotopic.
Let $H:I\times B\rightarrow \tilde \bR_n$ be a homotopy.
Then $p\circ H:S^1\times B\rightarrow \tilde \bR_n(1)$ can be
considered as a map $h:B\rightarrow \Map(S^1,\tilde \bR_n(1))$.
We have the following diagram
\begin{equation}\label{de9}\begin{array}{ccc}
\{b\}&\stackrel{x}{\rightarrow}&\Map(S^1,\tilde \bR_n(1))\\
\downarrow&h\nearrow &\ev_1\downarrow\\
B&\stackrel{f}{\rightarrow}&\tilde \bR_n(1)\end{array}\ ,\end{equation}
where $\ev_1:\Map(S^1,\tilde \bR_n(1))\rightarrow \tilde \bR_n(1)$
is the evaluation at $1\in S^1$. 

Vice versa, if $x\in \pi_1(\tilde \bR_n(1))$ is such that the diagram
above admits a lift $h$, then $\tilde f_0$ and $x\tilde f_0$ are
homotopic.

\subsubsection{}\label{clh1}

The existence problem for a lift $h$ can be studied using obstruction
theory. The fiber of the map $\ev_1$ is $\Omega \tilde \bR_n(1))$. In the proof of Lemma \ref{ght} we have
seen that the homotopy groups of $\tilde \bR_n(1)$ are given by
$$\begin{array}{|c|c|c|c|c|c|}\hline
i&0&1&2&3&\ge 4\\\hline
\pi_i(\tilde \bR_n(1))&*&\Z^{\frac{n(n-1)}{2}}&\Z^{2n}&\Z&0\\\hline\end{array}
\ .$$\
It follows that the homotopy groups of $\Omega  \tilde \bR_n(1)$
are given by 
$$\begin{array}{|c|c|c|c|c|}\hline
i&0&1&2&\ge 3\\\hline
\pi_i(\Omega \tilde \bR_n)&Z^{\frac{n(n-1)}{2}}&\Z^{2n}&\Z&0\\\hline\end{array}
\ .$$
We therefore have obstructions in
$H^1(B,\Z^{\frac{n(n-1)}{2}})$, $H^2(B,\Z^{2n})$ and $H^3(B,\Z)$, and a
general discussion seems to be complicated.

\subsubsection{} \label{p2}

Let $\tilde P_n$ be the set-valued functor classified by $\tilde \bR_n$.
For each space $B$  we have a natural transformation
$p_B:\tilde P_n(B)\rightarrow \tilde P_{(1)}(B)$ induced by composition with the
map $p:\tilde \bR_n\rightarrow
\tilde \bR_n(1)$.

We conclude from Subsection \ref{bc1} that the fibers of $p_B$ are homogeneous spaces over 
$H^0(B,H^2(T^n,\Z))$. 

Consider a class $[F,h]\in \tilde P_{(1)}(B)$. By Subsection
\ref{bc1}, it belongs to the image of $p_B$ if and only the
restriction
$[F,h]_{|B^{(1)}}$
to a $1$-skeleton $B^{(1)} \subset B$ is trivial.
Since every principal torus bundle on a $1$-dimensional complex is trivial, this is
equivalent to the condition that the restriction of $h$ to $F_{|B^{(1)}}$ is
trivial, or equivalently
 $h\in \cF^2H^3(F,\Z)$. 
Let us fix a map $f:B\rightarrow \tilde \bR_n(1)$ representing a pair
$[F,h]$ with this property. If we choose a homotopy from $\phi\circ f$ to the constant
map, then we distinguish an element in the fiber $p_B^{-1}([F,h])$.

\subsubsection{}

In Section \ref{fd} we have introduced a pair $(\bF_n,\bh_n)$ over
$\bR_n$ such that $\bh_n\in \cF^2H^3(\bF_n,\Z)$. 
The isomorphism class of this pair gives rise to a classifying map
$f(1):\bR_n\rightarrow \tilde \bR_n(1)$, well defined upto homotopy.

\begin{lem}\label{var}
The set of homotopy classes of maps $f$ which are lifts of $f(1)$ in the diagram
$$\begin{array}{ccc}
&&\tilde \bR_n\\
&f\nearrow&p\downarrow\\
\bR_n&\stackrel{f(1)}{\rightarrow}&\tilde \bR_n(1)\end{array}
$$ 
is a torsor over $H^2(T^n,\Z)$.
\end{lem}
\proof
Since $\bR_n$ is simply connected we know that lifts exist.
Let now $x\in \pi_1(\tilde \bR_n(1))$.
We choose a base point $b\in \bR_n$ and identify
$\pi_1(\tilde \bR_n(1))\cong H^2(\bF_{n,b},\Z)$,
where $\bF_{n,b}$ denotes the fiber of $\bF_n$ over $b$.
In particular we view  $x\in H^2(\bF_{n,b},\Z)$.

We must show that the existence of a lift $h$ in
the diagram \ref{de9}
(with $B$ replaced by $\bR_n$ and $f$ replaced by $f(1)$) implies that $x=0$.
Assume that a lift $h$ exists, adjoint to a map $H\colon S^1\times \bR_n\to
\tilde \bR_n(1)$. This corresponds to a $T^n$-bundle $F\rightarrow S^1\times \bR_n$ and a class
$h\in H^3(F,\Z)$. Let $\pr:S^1\times \bR_n\rightarrow \bR_n$ be the
projection. Since $\bR_n$ is simply connected $\pr$ induces an
isomorphism in second cohomology and the bundle
$F$ is the pull-back via $\pr$ of a $T^n$-bundle from $\bR_n$.
Since $H$ restricts to $f$ on $\{1\}\times \bR_n$ and the corresponding
$T^n$-bundle is $\bF_n$, necessarily $F\cong \pr^*\bF_n=S^1\times
\bF_n$. By the K{\"u}nneth formula $H^3(F,\Z)\cong H^3(\bF_n,\Z)\oplus
H^2(\bF_n,\Z)$, with corresponding decomposition
$h=\bh_n\oplus u$. By the definition of $h$ and the calculation of
$\pi_1(\tilde\bR_n(1))$ in Subsection \ref{cl2},
the restriction of $u$ to $\bF_{n,b}$ is $x$. Since the restriction $H^2(\bF_n,\Z)\rightarrow
H^2(\bF_{n,b},\Z)$ is trivial by the description of $H^2(\bF_n,\Z)$
given in Lemma \ref{re1},
it follows that $x=0$. \hB

\subsubsection{}

We now fix one choice of $f$ in Lemma \ref{var}.

\begin{prop}\label{var1}
The map $f:\bR_n\rightarrow \tilde \bR_n$ is a weak homotopy
equivalence.
\end{prop}
\proof
By  Lemma \ref{hh} and Lemma \ref{ght}
 it suffices to show that
$f$ induces  isomorphisms on $\pi_2$ and $\pi_3$.
 
Note that $H^2(\tilde \bR_n,\Z)\cong \Hom_\Z(\pi_2(\tilde
\bR_n),\Z)\cong \Z^{2n}$.
We have a natural map
$x:\tilde \bR_n\rightarrow K(\Z^n,2)$
which classifies the $T^n$-principal bundle $\tilde
\bF_n\rightarrow \tilde \bR_n$. 
We study the corresponding first Chern classes $x_i\in H^2(\tilde
\bR_n,\Z)$, $i=1,\dots,n$. 
% Now observe that the
%map $\Z^n\rightarrow  H^2( \tilde \bR_n(2),\Z)$ is dual to the
%projection  $\pi_2(\tilde \bR_n(1))\rightarrow
%\Z^n$. Therefore we can extend the sequence 
%$x_1,\dots,x_n$  to a basis of
%$H^2(\tilde \bR_n(2),\Z)$. 

Let us consider the second page of the
Leray-Serre spectral sequence ${}^{\tilde p} E_2^{p,q}$ of the fibration
$\tilde p: \tilde \bF_n\rightarrow \tilde \bR_n$.
Since $\bh_n\in \cF^2 H^3(\tilde \bF_n,\Z)$, there are elements
$\hat x_1,\dots,\hat x_n\in H^2(\tilde \bR_n,\Z)$
% and $a_{i,j}\in \Z$
such that
$\tilde \bh_n^{2,1}$ is represented by
$\tilde h:=
%\sum_{i,j} a_{i,j} y_i\otimes x_j + 
\sum_{i} y_i\otimes \hat x_i\in
{}^{\tilde p} E_2^{2,1}$, where the $y_i$ are the canonical generators of
$H^1(T^n,\Z)$ (for the fibre). 

Under the pull-back induced by $f$ the spectral sequence ${}^{\tilde
  p} E$
is mapped to the spectral sequence ${}^{\pi_n} E$ considered  in \ref{cd1}.
In particular $[\tilde h]\in {}^{\pi_n}E_\infty^{2,1}$ is mapped to
$[\sum_{i=1}^n y_i\otimes \hat \bc_i]\in {}^{\pi_n} E_\infty^{2,1}$.
%Using that $f^*x_i=\bc_i$ 
We conclude that
$f^*\hat x_i=\hat \bc_i
%+\sum_{i,j} a_{i,j} \bc_j
 $.
It follows that $f^*:H^2(\tilde \bR_n,\Z)\rightarrow
H^2(\bR_n,\Z)$ maps
$(x_1,\dots,x_n,\hat x_1,\dots,\hat x_n)$ to the basis $(\bc_1,\dots,\bc_n,\hat
\bc_1,\dots,\hat \bc_n)$ and is therefore
surjective and injective.
This implies that $f_*:\pi_2(\bR_n)\rightarrow \pi_2(\tilde \bR_n)$
is an isomorphism.

It now suffices to show that $f:\pi_3(\bR_n)\rightarrow \pi_3(\tilde
\bR_n)$ is surjective.
A generator $g\in \pi_3(\tilde \bR_n)$ is represented by
a map $g:S^3\rightarrow \tilde \bR_n$. The corresponding pair is the
trivial torus bundle $\pr_1:S^3\times T^n\rightarrow S^3$ with the
cohomology class of the form  $h=\pr_1^*z$ for some $z\in
H^3(S^3,\Z)$ which is a generator. 

It suffices to show that the isomorphism class of pairs $[S^3\times
T^n,\pr_1^*z]$ is
the pull back of $[\bF_n,\bh_n]$ on $\bR_n$. 
We consider the composition
$g:S^3\rightarrow K(\Z,3)\stackrel{i}{\rightarrow} \bR_n$, 
where the map $i$ was defined in \ref{f16} and the first map
realizes a generator of $\pi_3(K(\Z,3))$. It then follows immediately from 
Lemma \ref{g16} that $g^*[\bF_n, \bh_n]=[S^3\times
T^n,\pm \pr_1^*z]$. Choosing the opposite generator of
$\pi_3(K(\Z,3))$, if necessary, the assertion follows. \hB 

\subsubsection{}\label{epsd}

\begin{kor}\label{hhll}
The functors $P_n$ classified by $\bR_n$ and
$\tilde P_n$ classified by $\tilde \bR_n$ are naturally isomorphic.
The group $H^2(T^n,\Z)$ acts freely on the set of
such isomorphisms, and it acts transitively if we fix the composition with
$p_*\colon \tilde P_n\to \tilde P_{(1)}$.
\end{kor}

%Let $B$ be a space.
%Note that the group $G_n$ of $T$-duality transformations acts
%by automorphisms on $P(B)$. In particular each element
%$x\in P_n(B)$ has a canonical $T$-dual $\hat x:=T(x)$ (see \ref{tt1}).

%Assume that we have fixed an isomorphism $\epsilon:P_n\cong \tilde P_n$.
%Then we can transfer the $T$-duality transformations to $\tilde P_n$.
%Assume that $(F,h)$ is a pair over $B$. 

\subsubsection{}
We consider a pair $(F,h)$ over a space $B$. Recall the notion of an
extension of a pair to a $T$-duality triple \ref{exexdef}.

%\begin{ddd}\label{drt}
%We say that $(F,\cH)$ extends
%$(\hat F,\hat h)$ if there exists a lift $x\in \tilde P_n(B)$
%of $[F,h]\in \tilde P_{(0)}(B)$ such that
%$(\hat F,\hat h)$ is the underlying pair of
%$\epsilon\circ T\circ \epsilon^{-1}(x)$.
%\end{ddd}

%Recall that $\cF^2H^3(F,\Z)\subset H^3(F,\Z)$ denotes the second step
%of the Leray-Serre filtration of the cohomology associated to the
%bundle
%$F\rightarrow B$ (see \ref{filtex}).
\begin{theorem}\label{exex}
The pair $(F,h)$ admits an extension to a $T$-duality triple
$((F,\cH),(\hat F,\hat \cH),u)$
if and only if
$h\in \cF^2 H^3(F,\Z)$.
\end{theorem}
\proof
The condition is necessary by definition of a triple. We show that it is also sufficient.
Let $g:\tilde \bR_n\rightarrow \bR_n$ be a homotopy inverse of a
choice of $f$ in Proposition \ref{var1}. Then
$g^*((\bF_n,\cH_n),(\hat \bF_n,\hat \bH_n),\bu_n)$ is a $T$-duality
triple over $\tilde \bR_n$ such that 
$g^*(\bF_n,[\cH_n])$ is isomorphic to $(\tilde \bF_n,\tilde
\bh_n)$. Therefore we can consider
 $g^*((\bF_n,\cH_n),(\hat \bF_n,\hat \cH_n),\bu_n)$ as an extension of
$(\tilde \bF_n,\tilde
\bh_n)$ to a $T$-duality triple.

Let now $(F,h)$ be a pair over $B$ with classifying map
$\phi:B\rightarrow \tilde \bR_n(0)$. Every lift $\Psi$ in the diagram
$$\begin{array}{ccc}
&&\tilde \bR_n\\
&\Psi \nearrow&\downarrow\\
B&\stackrel{\phi}{\rightarrow}&\tilde \bR_n(0)\end{array}$$
yields an extension $\Psi^*g^*((\bF_n,\cH_n),(\hat \bF_n,\hat
\bH_n),\bu_n)$ of $(F,h)$ to a $T$-duality
triple.

The lifting problem can be decomposed into two stages
$$\begin{array}{ccc}
B&\stackrel{\Psi}{\rightarrow}&\tilde \bR_n\\
\|&&p\downarrow\\
B&\stackrel{\Psi_1}{\rightarrow}&\tilde \bR_n(1)\\
\|&&\downarrow\\
B&\stackrel{\phi}{\rightarrow}&\tilde \bR_n(0)\end{array}$$
The lift $\Psi_1$ exists since $h\in \cF^1 H^3(F,\Z)$.
Since $ p$ is the universal covering map, the stronger condition 
$h\in \cF^2 H^3(F,\Z)$ implies the existence of the lift $\Psi$ in the
second stage. \hB 

\section{$T$-duality transformations in twisted cohomology}

\subsubsection{}\label{w2w2w}

In \cite{bunkeschickt}, Section 3.1,  we have introduced a set of axioms
which describe the basic properties of twisted cohomology
theories.  
Below we will only use these axioms. Let $h(\dots,\dots)$ be a twisted
cohomology theory.

\subsubsection{}
Let $((F,\cH),(\hat F,\hat \cH),u)$ be a $T$-duality triple over a
space $B$. It gives rise to the 
diagram
\begin{equation}
\begin{array}{ccccc}
&&F\times_B\hat F&&\\
&p\swarrow&&\searrow \hat p\\
F&&r\downarrow &&\hat F\\
&\pi\searrow&&\swarrow\hat \pi&\\
&&B&&\end{array}\ .\end{equation}
Recall the Definition \ref{trd} of the $T$-duality transformation
$$T:=\hat p_!\circ u^*\circ p^*:h(F,\cH)\rightarrow h(\hat F,\hat
\cH)\ .$$

\subsubsection{}
There is a unique $1$-dimensional $T$-duality triple over a point.
Recall the following definition from \cite{bunkeschickt}.
\begin{ddd}\label{tad3}
The twisted cohomology theory $h$ is called $T$-admissible, if the 
$T$-duality transformation associated with the unique one-dimensional
$T$-duality triple over a point is an isomorphism (of degree $-1$).
\end{ddd}

\subsubsection{}
\begin{theorem}\label{isot}
If $((F,\cH),(\hat F,\hat \cH),u)$ is a $T$-duality triple over a finite
dimensional
CW-complex $B$, and if
$h$ is a $T$-admissible twisted cohomology theory, then the $T$-duality
transformation \ref{trd}  is an isomorphism.
\end{theorem}
\proof
Exactly as in \cite{bunkeschickt}, Proof of Theorem 3.13,  one uses
induction on the dimension of the base space, the Mayer-Vietoris exact sequence and the $5$-lemma to reduce this to the case $B=*$.
For $B=*$ one observes that the $n$-dimensional $T$-duality
transformation is an iterated $1$-dimensional $T$-duality
transformation which is an isomorphism by the definition of
$T$-admissibility.
\hB

\section{Classification of $T$-duality triples and extensions}

%\subsubsection{}

%In this section we characterize the space $\bR_n$ as the classifying
%space of $T$-duality triples. The notion of a $T$-duality triple was
%introduced in Definition \ref{tr13}. Let $x:=((F,\cH),(\hat F,\hat
%\cH),u)$ and $x^\prime:=((F^\prime,\cH^\prime),(\hat F^\prime,\hat
%\cH^\prime),u^\prime)$ be two $T$-duality triples.
%\begin{ddd}\label{isodd}
%We say that $x$ and $x^\prime$ are isomorphic, if there exist
%bundle isomorphisms
%$$\begin{array}{ccc}
%F&\stackrel{\psi}{\rightarrow}&F^\prime\\
%\downarrow&&\downarrow\\
%B&=&B\end{array}\ ,\quad \begin{array}{ccc}
%\hat F&\stackrel{\hat\psi}{\rightarrow}&\hat F^\prime\\
%\downarrow&&\downarrow\\
%B&=&B\end{array}$$
%and isomorphisms of twists
%$$v:\psi^*\cH^\prime\rightarrow \cH\ ,\quad \hat v:\hat
%\psi^*\hat\cH^\prime\rightarrow \hat \cH$$
%such that the composition
%$$\hat p^*\hat \cH\stackrel{p^*v^{-1}}{\rightarrow} \hat p^* \hat
%\psi^*\hat \cH^\prime\cong (\psi,\hat\psi)^*(\hat p^\prime)^*\hat
%\cH^\prime
%\stackrel{(\psi,\hat \psi)^*u^\prime}{\rightarrow}(\psi,\hat\psi)^*(p^\prime)^*
%\cH^\prime\cong p^*\psi^* \cH^\prime\stackrel{p^* v}{\rightarrow}
%p^*\cH$$
%is equal to $u:\hat p^*\hat \cH\rightarrow p^*\cH$.
%Here $(\psi,\hat \psi):F\times_B\hat F\rightarrow F^\prime\times_B\hat
%F^\prime$ is the induced map, and
%$p:F\times_B\hat F\rightarrow F$, $\hat p:F\times_B\hat F\rightarrow
%\hat F$, $p^\prime:F^\prime\times_B\hat
%F^\prime\rightarrow F^\prime$, and $\hat p^\prime:F^\prime\times_B\hat
%F^\prime\rightarrow \hat F^\prime$ are the projections.
%\end{ddd}

\subsubsection{}

%For each space $B$ let $\Triple_n(B)$ denote the set of isomorphism
%classes of $n$-dimensional $T$-duality triples over $B$. Given by map
%$f:B^\prime\rightarrow B$ and a triple $x$ over $B$ there is a natural
%definition of a functorial pull-back $f^*x$. We observe that pull-back
%preserves isomorphism classes. In this way obtain a contravariant
%set-valued functor $B\mapsto \Triple_n(B)$.
Recall from \ref{tripleintro} the definition of the functor $B\rightarrow \Triple_n(B)$
which associates to a space the set of isomorphism classes (see Definition
\ref{isodd})
 of $n$-dimensional $T$-duality triples. 
%The following result
%is very expected but it took a while to write out all details of a
%formal proof.

\begin{lem}\label{homo32}
The functor $B\mapsto \Triple_n(B)$ is homotopy invariant.
\end{lem}
\proof
Let $x=((F,\cH),(\hat F,\hat \cH),u)$ be a triple over $B$.
Let $H:[0,1]\times B\rightarrow B^\prime$ be a map.
Let $f_0,f_1$ denote the evaluations of $H$ at the endpoints of the
interval. We must show that $f_0^*x\cong f_1^*x$.
%In the argument below we will use repeatedly the following two facts
%about twists.  
%\begin{enumerate}
%\item Isomorphisms classes of twists over a space $X$ are classified by
%$H^3(X,\Z)$.
%\item
%The set of isomorphisms between two twists over $X$ is naturally a torsor over
%$H^2(X,\Z)$.
%\end{enumerate}
Let $\pr:[0,1]\times B\rightarrow B$ denote the projection.
In the first step we choose  bundle isomorphisms which extend the
identity at $0$
$$\begin{array}{ccc}H^*F&\stackrel{\psi}{\leftarrow}& \pr^* f_0^*F\\
\downarrow&&\downarrow \\[0pt] 
 [0,1]\times B&=& [0,1]\times B\end{array}\ , \quad \begin{array}{ccc}H^*\hat F&\stackrel{\hat
  \psi}{\leftarrow}& \pr^* f_0^*\hat F\\
\downarrow&&\downarrow\\[0pt]
[0,1]\times B&=&[0,1]\times B\end{array}\ .$$
Evaluation at $\{1\}\times B$ gives isomorphisms
$\psi_1:f_0^*F\rightarrow f_1^*F$ and
$\hat \psi_1:f_0^*\hat F\rightarrow f_1^*\hat F$. 

Let $\Phi_0:f_0^*F\rightarrow F$, $\hat \Phi_0:f_0^*\hat F\rightarrow
\hat F$, $\Phi_1:f_1^*F\rightarrow F$, $\hat \Phi_1:f_1^*\hat F\rightarrow
\hat F$, 
$\gamma:H^*F\rightarrow F$, and $\hat \gamma:H^*\hat F\rightarrow \hat
F$ denote the induced maps. Let $J_i:f_0^*F\rightarrow \pr^*f_0^*F$ and
$\hat J_i:f_0^*\hat F\rightarrow \pr^*f_0^*\hat F$ be the canonical
embeddings over $B\rightarrow \{i\}\times B\subset [0,1]\times B$,
$i=0,1$. Since $\gamma\circ \psi\circ J_0=\Phi_0$ and
$\hat \gamma\circ \hat \psi\circ \hat J_0=\hat \Phi_0$
we have canonical isomorphisms of twists
$$J_0^*\psi^*\gamma^*\cH\cong \Phi_0^*\cH,\quad\hat J_0^*\hat \psi^*
\hat \gamma^*\hat
\cH\cong \hat \Phi_0^*\hat \cH\ .$$ 
These uniquely extend to isomorphisms of twists
$$w:\psi^*\gamma^* \cH\rightarrow \pr^* \Phi_0^*\cH, \quad\hat w:\hat \psi^*\hat \gamma^*
\hat \cH\rightarrow \pr^* \hat \Phi_0^*\hat \cH\ .$$
We now restrict to $\{1\}\times B$ and obtain isomorphisms of twist
$v:\psi_1^*\Phi_1^*\cH\cong
 J_1^*\psi^*\gamma^*\cH\stackrel{J_1^*w}{\rightarrow}J_1^*\pr^*\Phi_0^*\cH\cong \Phi_0^*\cH$ and
$\hat v:\hat \psi_1^*\hat \Phi_1^*\hat \cH\cong
 \hat J_1^*\hat \psi^*\hat \gamma^*\hat \cH\stackrel{\hat J_1^*\hat
   w}{\rightarrow}\hat J_1^*\pr^*\hat \Phi_0^*\hat \cH\cong \hat
 \Phi_0^*\hat \cH$.
We have the following diagram of isomorphisms of twists
over $I\times (f_0^*F\times_B f_0^*\hat F)$:
$$\begin{array}{ccc}
\hat p^*\pr^*\hat \Phi_0^*\hat \cH&\stackrel{\pr^* (\Phi_0,\hat
  \Phi_0)^* u}{\rightarrow}& p^*\pr^*\Phi_0^*\cH
\\
\hat p^*\hat w\uparrow&& p^* w \uparrow\\
\hat p^*\hat \psi^*\hat \gamma^* \hat \cH&\stackrel{(\gamma\circ
  \psi,\hat \gamma\circ \hat \psi)^*u}{\rightarrow}&
p^*\psi^*\gamma^*\cH\end{array} \ .$$
Here $p:I\times (f_0^*F\times_B f_0^* \hat F)\rightarrow I\times f_0^*F$ and
$\hat p:I\times (f_0^*F\times_B f_0^* \hat F)\rightarrow I\times f_0^*\hat F$
are the projections.
This diagram commutes after restriction to $\{0\}\times
B$. Hence it commutes, and so does its restriction to
$\{1\}\times B$. The latter restriction gives
$$\begin{array}{ccc}
\hat p_1^*\hat \Phi_0^*\hat \cH&\stackrel{(\Phi_0,\hat
  \Phi_0)^* u}{\rightarrow}& p_1^*\Phi_0^*\cH
\\
\hat p_1^*v\uparrow&&\hat p_1^* v \uparrow\\
\hat p_1^*\hat \psi_1^* \hat \Phi_1^* \hat \cH&\stackrel{(\psi,\hat
  \psi_1)^* (\Phi_1,\hat \Phi_1)^*u}{\rightarrow}&
p_1^*\psi_1^*\Phi_1^*\cH\end{array} \ $$
with the projections
$p_1:f_0^*F\times_B f_0^* \hat F\rightarrow f_0^*F$ and
$\hat p_1:f_0^*F\times_B f_0^* \hat F\rightarrow f_0^*\hat F$.

This shows that the bundle isomorphisms
$\psi_1,\hat \psi_1$ and the isomorphisms of twists $v,\hat v$
form an isomorphism of triples $f_0^*x\cong f_1^*x$. \hB

\subsubsection{}\label{overdef}

Let us fix the $T^n$-bundles $F$ and $\hat F$ over $B$. Then we can consider triples
of the form $x:=((F,\cH),(\hat F,\hat \cH),u)$, 
$x^\prime:=((F,\cH^\prime),(\hat F,\hat \cH^\prime),u^\prime)$.
\begin{ddd} \label{def:isoextensiondef}We say that $x$ is 
isomorphic to $x^\prime$ over $(F,\hat F)$ if there exists an
isomorphism of triples $x\cong x^\prime$ such that the underlying
bundle isomorphisms are the identity maps (see \ref{isodd}).
\end{ddd}
Let $\Triple_n^{(F,\hat F)}(B)$ denote the set of isomorphism classes
of triples over $(F,\hat F)$.

\subsubsection{}\label{ert2}

We have a natural action of $H^3(B,\Z)$ on $\Triple_n^{(F,\hat F)}(B)$.
Let $\delta\in H^3(B,\Z)$ and choose a twist $\cV$ in the
corresponding isomorphism class. If $x:=((F,\cH),(\hat F,\hat \cH),u)$
represents $[x]\in \Triple^{(F,\hat F)}_n(B)$,
then we define the triple
\begin{equation}\label{eq:def_of_r}x+\cV:=((F,\cH\otimes \pi^*\cV),(\hat F,\hat \cH\otimes \hat\pi^*\cV),u\otimes r^*
\id_{\cV})\ ,\end{equation} and $[x]+\delta:=[x+\cV]$,
 where $\pi\colon F\to B$, $\hat\pi\colon \hat F\to B$, and $r\colon
 F\times_B\hat F\rightarrow B$  are the projections.

\begin{prop}\label{freetrans}
$H^3(B,\Z)$ acts freely and transitively on the set 
$\Triple^{(F,\hat F)}_n(B)$.
\end{prop}
\proof
We prove Proposition \ref{freetrans} in several steps. 
First recall the following principle, which we will employ several times. Let a group $G$ act on a set
$S$. If $\tau$ is an equivalence relation on $S$ preserved by $G$,
then $G$ acts freely and transitively on $S$ if and only if the two
following conditions hold:
\begin{enumerate}
\item 
$G$ acts
transitively on the set $S/\tau$ of equivalence classes
\item the isotropy group of one (and hence every) equivalence class
  acts freely and transitively on this equivalence class.
\end{enumerate}

\subsubsection{}
Since $F$ and
$\hat F$ are fixed at the moment, we use the notation  $[\cH,\hat
\cH,u]$ for the class $[(F,\cH),(\hat F,\hat \cH),u]$. Let us first
introduce an equivalence relation $\tau_1$ on the set
$\Triple^{(F,\hat F)}_n(B)$ by defining that $[\cH,\hat\cH,u]$ is
  equivalent to $[\cH',\hat\cH',u']$ if and only if $\cH\iso \cH'$ and
  $\hat \cH\iso \hat \cH'$. The action of $H^3(B,\Z)$
  preserves this equivalence relation.

  \begin{lem}\label{lem:only_twist_iso_classes}
     The action of $H^3(B,\Z)$ on the set of equivalence classes
     $\Triple^{(F,\hat F)}_n(B)/\tau_1$ is transitive with stabilizer
  $\ker(\pi^*)\cap \ker(\hat\pi^*)$.
  \end{lem}
  \proof
     The assertion about the stabilizer  is
     an immediate consequence of the definitions. We now verify
     transitivity of the action on the set of equivalence classes. If
     $[\cH,\hat\cH,u]\in \Triple^{(F,\hat F)}_n(B)$, then by the
     definition of a 
$T$-duality triple, the leading part  $[\cH]^{2,1}=[\sum_{i=1}^n y_i\otimes
\hat c_i]\in {}^\pi E_\infty^{2,1}$ 
of $[\cH]$ is determined by $\hat F$, and similar  $[\hat\cH]^{2,1}$ is determined by $F$. 
Let us consider a second class $[\cH',\hat \cH',u']\in \Triple^{(F,\hat F)}_n(B)$.
It follows
from the structure of the spectral sequences that there are classes
$\delta,\hat \delta\in H^3(B,\Z)$ such that
$[\cH^\prime]=[\cH]+\pi^*\delta$ and $[\hat \cH^\prime]=[\hat
\cH]+\hat \pi^*\hat \delta$.

\begin{lem}
We have $\delta-\hat \delta\in \ker(\pi^*)+\ker(\hat\pi^*)$.
\end{lem}
\proof
We assume without loss of generality that $B$ is connected and choose
a base point $b\in B$.
Because of the presence of the isomorphisms $u$ and $u'$ of the
pullbacks of the twists to $F \times_B \hat F$
we know that $r^*(\delta-\hat \delta)=0$ (with $r$ as in  \eqref{eq:def_of_r}).
We choose twists $\cV$ and $\hat\cV$ representing $\delta$ and $\hat
\delta$, respectively. For our question, we can replace $\cH'$ by
$\cH\tensor \pi^*\cV$ and $\hat\cH'$ by $\hat\cH\tensor
\hat\pi^*\hat\cV$. Thus, with the choice of an appropriate
isomorphism of twists $w:r^*\cV\rightarrow 
r^*\hat\cV$ we can write $u^\prime=u\otimes w $.
% for some $z\in
%H^2(F\times_B \hat F,\Z)$.
Note that we have  canonical isomorphisms $\cV_{|\{b\}}\cong 0$,
$\hat\cV_{|\{b\}}\cong 0$.
%\cH^\prime_{|F_b}$ and $\hat \cH_{|\hat F_b}\cong \hat
%\cH^\prime_{|\hat F_b}$. 
Therefore we can consider 
$$(w(b):0\stackrel{\sim}{\rightarrow} r^*\cV_{|F_b\times \hat F_b}\stackrel{w_{|F_b\times \hat F_b}}{\rightarrow} r^*\hat \cV_{|F_b\times \hat F_b}\stackrel{\sim}{\rightarrow}0 )\in H^2(F_b\times
\hat F_b,\Z)\ .$$ We know by Lemma \ref{lem:i891} that
${}^rd_2^{0,2}(w(b))=0$ and
${}^rd_3^{0,2}(w(b))=\delta-\hat \delta+\im({}^rd_2^{1,1})$.
Let $W:=H^2(F_b\times
\hat F_b,\Z)/(\im(p_b^*)+ \im(\hat
p_b^*))$, and let $[\dots]$ denote for the moment classes in this
quotient. Note that by \ref{uudef} we have well-defined classes
$[u(b)],[u^\prime(b)]\in W$. Since $u$ and $u^\prime$ satisfy the
condition $\cP$ we have $[u(b)]=[u^\prime(b)]$, and it follows that
$[w(b)]=[u'(b)]-[u(b)]=0$. %, where $z(b):=z_{|F_b\times \hat F_b}$.
Accordingly, we choose a decomposition $w(b)=p_b^*(d)+\hat p_b^*(\hat d)$ for
some $d\in H^2(F_b,\Z)$ and $\hat d\in H^2(\hat F_b,\Z)$.
Then
\begin{equation*}
0={}^rd_2^{0,2}(w(b))={}^rd_2^{0,2}(p_b^*(d))+{}^rd_2^{0,2}(\hat
p_b^*(\hat d )).
\end{equation*}
Since $\im({}^r d_2^{0,2} p_b^*) \cap \im({}^r d_2^{0,2}\hat p_b^*)=0$,
${}^\pi d_2^{0,2}(d)=0$ and 
${}^{\hat \pi} d_2^{0,2}(\hat d)=0$. 
Therefore
\begin{equation*}
\delta-\hat \delta+\im({}^rd_2^{1,1})={}^rd_3^{0,2}(w(b))=
{}^\pi d_3^{0,2}(d)+{}^{\hat \pi}d_3^{0,2}(\hat d)
+\im({}^rd_2^{1,1}).
\end{equation*}
Finally we use that 
$\im({}^rd_2^{1,1})=\im({}^{\pi}d_2^{1,1})+\im({}^{\hat
  \pi}d_2^{1,1})$ since $H^1(F_b\times \hat F_b,\Z)\iso
H^1(F_b,\Z)\oplus H^1(\hat F_b,\Z)$.
Using the relation between the images of ${}^\pi
d_*^{*,*}, {}^{\hat\pi}
d_*^{*,*}$ and the kernels of $\pi^*$ and $\hat \pi^*$ we obtain the assertion.\hB 

Write now $\delta-\hat\delta=-a+b$ with $a\in\ker(\pi^*)$ and
$b\in\ker(\hat\pi^*)$. Then $e:=\delta+a=\hat\delta+b$, and
$[\cH']=[\cH]+\pi^*\delta=[\cH]+\pi^*e$,
$[\hat\cH']=[\hat\cH]+\hat\pi^*\hat\delta=[\hat\cH]+\hat\pi^*e$.
It follows that $[\cH^\prime,\hat \cH^\prime,u^\prime]\sim_{\tau_1}
[\cH,\hat \cH,u]+e$. This
proves the transitivity statement of  Lemma \ref{lem:only_twist_iso_classes}. 
\hB 

\subsubsection{}

Fix now an equivalence class for $\tau_1$. In fact, we can fix the
corresponding twists $\cH$ and
$\hat \cH$. The isomorphism classes over $(F,\hat F)$ in the given
$\tau_1$-equivalence class  can be represented in the form
$[u]:=[\cH,\hat \cH,u]$ for varying $u$. By  
$\Triple(\cH,\hat\cH)$ we denote the set of these isomorphism classes.

We introduce an equivalence relation $\tau_2$ on
$\Triple(\cH,\hat\cH)$.  Let $i\colon F_b\times\hat F_b\into F\times_B\hat F$ be the inclusion of the
  fiber over the base point $b\in B$.
We declare $[u], [u']\in
\Triple(\cH,\hat \cH)$ to be  $\tau_2$-equivalent if and only if the class $w\in H^2(F\times_B\hat
F,\Z)$ which is determined by $u+w=u^\prime$ satisfies 
\begin{equation}
i^*w \in i^*(p^*H^2(F,\Z)+\hat p^*
H^2(\hat F,\Z))\ . \label{eq:cond_for_tau2}
\end{equation}
Note that $[u_1]=[u_2]\in
\Triple(\cH,\hat \cH)$ if and only if there are $a\in H^2(F,\Z)$,
$\hat a\in H^2(\hat F,\Z)$ (representing automorphisms of $\cH$ and
$\hat\cH$, respectively) such that $u_1= p^*a+\hat p^*\hat a+u_2$. It follows that, although the class $w\in H^2(F\times_B\hat F,\Z)$ is not
determined by the classes $[u],
[u']\in\Triple(\cH,\hat \cH)$, the condition
\eqref{eq:cond_for_tau2} is well defined, and it defines an equivalence relation.

\begin{lem}\label{lem:tau_2_equiv}
The subgroup
  $\ker(\pi^*)\cap \ker(\hat\pi^*)\subset H^3(B,\Z)$ preserves
  $\tau_2$,  acts
  transitively on $\Triple(\cH,\hat \cH)/\tau_2$, and the
  isotropy subgroup of each 
  $\tau_2$-equivalence class is 
  \begin{equation*}
    \im({}^\pi d_2^{1,1})\cap \im({}^{\hat\pi}d_2^{1,1})\subset
    \ker(\pi^*)\cap \ker(\hat\pi^*)\subset H^3(B,\Z).
  \end{equation*}
 \end{lem}
\proof
Choose $\delta\in \ker(\pi^*)\cap \ker(\hat\pi^*)\subset H^3(B,\Z)$
and a twist $\cV$ representing $\delta$. In order to describe the action of
$\delta$ on $\Triple(\cH,\cH')$, we choose trivializations $w\colon
\pi^*\cV\stackrel{\sim}{\rightarrow} 0$ and $\hat w\colon \hat\pi^*\cV\stackrel{\sim}{\rightarrow}
0$. If $[u]\in \Triple(\cH,\cH')$, then $[u]+\delta$ is represented by $[u\tensor
\hat p^*\hat w\circ p^* w^{-1}]$. We introduce the cohomology class 
$v:=\hat p^*\hat w\circ p^*w^{-1}\in H^2(F\times_B\hat F,\Z)$. 

Assume that $[u]\sim_{\tau_2} [u^\prime]$, i.e.  the class $w\in H^2(F\times_B\hat
F,\Z)$ defined by
$u^\prime=u+w$ satisfies
(\ref{eq:cond_for_tau2}). Then we have
$[u]+\delta\sim_{\tau_2}[u^\prime]+\delta$, since
$[u]+\delta =[u+v]$, $[u^\prime]+\delta=[u^\prime+v]$, and
$u^\prime+v=u+v+w$. This shows that the action of $\ker(\pi^*)\cap \ker(\hat\pi^*)$ preserves $\tau_2$.

We now calculate the stabilizer of the $\tau_2$-equivalence classes. We choose $\delta\in\ker(\pi^*)\cap \ker(\hat\pi^*)$ fixing
$[u]$, and continue to use the above notation.
The restriction of $\pi^*\cV$ to the fiber $F_b$ is canonically trivial.
We can therefore consider $x:=-w_{|F_b}\in H^2(F_b,\Z)$. Similarly we have
$\hat x:=\hat w_{|\hat F_b}\in H^2(\hat F,\Z)$.
We can now write
 $i^*v=p_b^*(x)+\hat p_b^*(\hat x)$. 
By Lemma \ref{lem:i891}, ${}^{\pi} d^{0,2}_2(x)=0$, ${}^{\hat \pi}
d^{0,2}_2(\hat x)=0$, 
${}^\pi d^{0,2}_3(x)=-[\delta]\in H^3(B,\Z)/\im({}^\pi d_2^{1,1})$, and
${}^{\hat \pi}d^{0,2}_3(\hat
x)=[\delta]\in H^3(B,\Z)/\im({}^{\hat \pi} d_2^{1,1})$. 
Note that $x$ is well-defined modulo restrictions of classes in $H^2(F,\Z)$ to $F_b$. Similarly, $\hat x$ is well-defined modulo restrictions of $H^2(\hat F,\Z)$ to $\hat F_b$. Therefore the condition that
${}^\pi d^{0,2}_3(x)=0$ $and$ ${}^{\hat \pi} d^{0,2}_3(\hat x)=0$ is
independent of all choices. Furthermore, it is equivalent to
$\delta\in \im({}^\pi d^{1,1}_2)\cap \im({}^{\hat \pi} d^{1,1}_2)$.
This holds if and only if  $x$ is a restriction of a class in $H^2(F,\Z)$ to $F_b$, and  $\hat x$ is the restriction of a class in $H^2(\hat F,\Z)$ to $\hat F_b$.

We see that $\delta\in \im({}^\pi d^{1,1}_2)\cap \im({}^{\hat \pi} d^{1,1}_2)$ implies 
that $i^*v\in i^*p^*H^2(F,\Z)+i^*\hat p^*H^2(\hat
F,\Z)$, and therefore $[u]+\delta\sim_{\tau_2}[u]$.
Vice versa, if $i^*v\in i^*p^*H^2(F,\Z)+i^*\hat p^*H^2(\hat
F,\Z)$, then we write
$i^*v=i^*p^*y+i^*\hat p^*\hat y$ for $y\in H^2(F,\Z)$
and $\hat y\in H^2(\hat F,\Z)$. We then have
$p_b^*(x)+\hat p_b^*(\hat x)=p_b^*(y_{|F_b})+\hat p_b^*(\hat y_{|\hat F_b})$. This implies that $x=y_{|F_b}$ and $\hat x=\hat y_{|\hat F_b}$.
We now see that ${}^\pi d^{0,2}_3(x)=0$ $and$ ${}^{\hat \pi} d^{0,2}_3(\hat x)=0$, and therefore $\delta\in \im({}^\pi d^{1,1}_2)\cap \im({}^{\hat \pi} d^{1,1}_2)$.

It remains to prove transitivity. Choose a second class $[u']\in
\Triple(\cH,\cH')$. Then we have $u'=u+h$ with $h\in H^2(F\times_B\hat F,\Z)$
such that the $c:=i^*h$ is of the form $c=p_b^*(a)+\hat p_b(\hat a)$
for $(a,\hat a)\in  H^2(F_b,\Z)\oplus H^2(\hat F_b,\Z)$. Since
$c$ is in the image of $i^*$ we have $0={}^r d_2^{0,2}(c)={}^r
d_2^{0,2}(p_b^*(a))+ {}^{r} d^{0,2}_2(\hat p_b^*(\hat a))$.  Consequently, ${}^\pi d_2^{0,2}a=0$
and ${}^{\hat \pi}d^{0,2}_2(\hat a)=0$. Set $\delta:= {}^{\hat\pi}
d^{0,2}_3(\hat a)$. Then $\delta= -{}^{ \pi} d^{0,2}_3( a)$. It
follows that
$\delta\in \ker(\pi^*)\cap \ker(\hat \pi^*)$. By Lemma
\ref{lem:i891} and similarly to the above considerations
$[u]+\delta \sim_{\tau_2} [u']$. 
\hB

\subsubsection{}

We now fix a class $[u_0]\in \Triple(\cH,\hat \cH)$ and
let $$U:=\{[u]\in \Triple(\cH,\hat \cH)|[u]\sim_{\tau_2} [u_0]\}$$
denote the corresponding $\tau_2$-equivalence class.

\begin{lem}\label{lem:obvious_transitive_action_on_tau_2}
  The subgroup $\ker(i^*)\subset
  H^2(F\times_B\hat F,\Z)$ acts transitively on $U$, and  the stabilizer of each element
  is given by  $$p^*H^2(F,\Z)\cap \ker(i^*)+ \hat p^*H^2(\hat F,\Z)\cap
  \ker(i^*)\ $$
 Alternatively, if $$H_{\tau_2}
  :=(i^*)^{-1}i^*\left(p^*H^2(F,\Z)+\hat p^*H^2(\hat F,\Z)\right)\subset
  H^2(F\times_B\hat F,\Z)\ ,$$ then $H_{\tau_2}$ acts transitively on
  $U$  with point stabilizers
  $p^*H^2(F,\Z)+\hat p^* H^2(\hat F,\Z)$.

%Let $H_{\tau_2}\subset H^2(F\times_B\hat F,\Z)$ be the
%  inverse image under $i^*$ of $H^2(F_b)+H^2(\hat F_b)\subset
%  H^2(F_b\times \hat F_b)$. Then $H_{\tau_2}$ acts transitively on each set
%  $\Triple^{\tau_2}(\cH,\hat\cH,u)$ with point stabilizers
%  $p^*H^2(E)+\hat p^* H^2(\hat E)$.
\end{lem}
\proof 
Let $[u]\in U$ and $w\in \ker(i^*)$. Then $w$ obviously satisfies (\ref{eq:cond_for_tau2}), and therefore we have $[u+h]\sim_{\tau_2}[u]$. Hence $\ker(i^*)$ acts on $U$.

 Let now $[u],[u^\prime]\in U$. Then $u^\prime=u+w$ and  $w$ satisfies $(\ref{eq:cond_for_tau2})$. Let us write $i^*w=i^*p^*x+i^*\hat p^*\hat x$ for $x\in H^2(F,\Z)$ and $\hat x\in H^2(\hat F,\Z)$. We define $w^\prime:=w-p^*x-\hat p^*\hat x$. Then we have $[u+w]=[u+w^\prime]$ on the one hand, and $i^*w^\prime=0$ on the other.
It follows that $[u^\prime]=[u+w^\prime]$. This shows that $\ker(i^*)$ acts transitively on $U$.

Assume that $w\in \ker(i^*)$ and $[u+w]=[u]$. Then we have
$w=p^*x+\hat p^*\hat x$ for some $x\in H^2(F,\Z)$ and $\hat x\in H^2(\hat F,\Z)$. From $0=i^*w=i^*(p^*x+\hat p^*\hat x)$ we conclude that
$i^*p^*x=0$ and $i^*\hat p^*\hat x=0$. This shows the assertion about the point stabilizers.

In order to deduce the second assertion of the lemma we shall observe that the natural map
$\ker(i^*)\rightarrow H_{\tau_2}$ induces an isomorphism
\begin{equation}\label{u631}\frac{\ker(i^*)}{p^*H^2(F,\Z)\cap \ker(i^*)+ \hat p^*H^2(\hat F,\Z)\cap \ker(i^*)}
\stackrel{\sim}{\rightarrow} \frac{H_{\tau_2}}{p^*H^2(F,\Z)+\hat p^* H^2(\hat F,\Z)} .\end{equation}
Injectivity is clear since $p^*H^2(F,\Z)\cap \ker(i^*)+ \hat p^*H^2(\hat F,\Z)\cap \ker(i^*)=(p^*H^2(F,\Z)+ \hat p^*H^2(\hat F,\Z))\cap \ker(i^*)$.
In order to show surjectivity we consider $w\in H_{\tau_2}$.
We write $i^*w=i^*p^*x+i^*\hat p^*\hat x$ for some
$x\in H^2(F,\Z)$ and $\hat x\in H^2(\hat F,\Z)$. Then $w^\prime:=w-p^*x-\hat p^*x$ represents the same class as $w$ on the right hand side of (\ref{u631}), and it satisfies $i^*w^\prime=0$.\hB

\subsubsection{}\label{mudef3}

We define a map 
\begin{equation}\label{eq:constr_of_mu}
  \mu\colon \im({}^\pi d_2^{1,1})\cap \im({}^{\hat \pi}d_2^{1,1}) \to
  H_{\tau_2}/(p^* H^2(F,\Z)+\hat p^* H^2(\hat F,\Z))
\end{equation}
in the following way. Represent $a\in  \im({}^\pi d_2^{1,1})\cap
\im({}^{\hat \pi}d_2^{1,1}) \subset\ker(\pi^*)\cap
\ker(\hat\pi^*)\subset H^3(B,\Z)$ by a map $f\colon B\to K(\Z,3)$. Choose a homotopy $H\colon [0,1]\times
F\to K(\Z,3)$ between the constant map and $f\circ \pi$, and a
homotopy $\hat H\colon [1,2]\times \hat F\to K(\Z,3)$ between $f\circ
\hat \pi$ and the constant map. Note that the homotopy classes of
these homotopies can be modified by concatenation with homotopies from
the constant map to the constant map, i.e.~by elements of $H^3(\Sigma
F,\Z)\iso H^2(F,\Z)$ or of $H^2(\hat F,\Z)$, respectively.

We can now pull  back $H$ and $\hat H$ to $[0,1]\times F\times_B\hat
F$ and $[1,2]\times F\times_B\hat F$. By construction, these two pull-backs
coincide at $\{1\}\times F\times_B\hat F$, and therefore they can be
concatenated to a map $[0,2]\times F\times_B\hat F\to K(\Z,3)$. This concatenation is constant
at both ends of the interval $[0,2]$ and thus factors over the
suspension $\Sigma ((F\times_B\hat F)_+)$. Therefore it represents an element in
$H^3(\Sigma((F\times_B\hat F)_+),\Z)\iso H^2(F\times_B\hat F,\Z)$.
Its class $\mu(a)\in H^2(F\times_B\hat F,\Z)/(p^*H^2(F,\Z)+\hat
p^* H^2(\hat F,\Z))$ is independent of the choice of the homotopies $H$ and $\hat H$. We still have to check that $\mu(a)$ belongs to the subspace 
$H_{\tau_2}/(p^*H^2(F,\Z)+\hat
p^* H^2(\hat F,\Z))$. This will follow from a universal construction of
$\mu$ for such classes and is therefore postponed until we discuss
this universal construction, compare \eqref{eq:mu_with_right_trager}.

Observe that this secondary operation actually makes sense for all elements in
$\ker(\pi^*)\cap \ker(\hat\pi^*)$, provided that we allow arbitrary values
in $H^2(F\times_B\hat F,\Z)/(p^*H^2(F,\Z)+\hat
p^* H^2(\hat F,\Z))$. 
 
%The pullback of this cohomology class to the fiber over $b\in B$ is
%represented by the restriction of this map to the fiber. Observe,
%however, that the restriction of $H$ and $\hat H$ to the fiber over
%$b$ is a homotopy from the constant map to the constant map and
%therefore represents an element in $H^2(F_b)+H^2(\hat F_b)$,
%i.e.~$\mu(a)\in H_{\tau_2}$, as required.

\subsubsection{}

Note that the groups $\im({}^\pi d_2^{1,1})\cap \im({}^{\hat
    \pi}d_2^{1,1})$ and $H_{\tau_2}/(p^*H^2(F,\Z)+\hat p^* H^2(\hat F,\Z))$ act on $U$.

\begin{lem}\label{lem:calculate_final_action}
  The map $\mu\colon \im({}^\pi d_2^{1,1})\cap \im({}^{\hat
    \pi}d_2^{1,1}) \to H_{\tau_2}/(p^*H^2(F,\Z)+\hat p^* H^2(\hat F,\Z))$ is
  compatible with the action of both groups on
  $U$.  
\end{lem}
\proof 
Choose $[u]\in U$ and $a\in \im({}^\pi d_2^{1,1})\cap \im({}^{\hat
    \pi}d_2^{1,1})$. As in \ref{mudef3} we represent $a$ by a map
$f:B\rightarrow K(\Z,3)$. We choose a twist $\cK$ on $K(\Z,3)$ such that
$[\cK]\in H^3(K(\Z,3),\Z)$ is the canonical generator.
The homotopies $H$ and $\hat H$ (see  \ref{mudef3}) give rise to isomorphisms
of twists $w:0\stackrel{\sim}{\rightarrow}\pi^*f^*\cK$ and $
\hat w:\hat \pi^*f^*\cK\stackrel{\sim}{\rightarrow}0$. We consider
$\hat p^*\hat w\circ p^*w^{-1}\in H^2(F\times \hat F,\Z)$ and observe that
$[u]+a=[u+\hat p^*\hat w\circ p^*w^{-1}]$. The Lemma now follows since
by construction $\hat w\circ w^{-1}$ represents $\mu(a)$.
\hB

\subsubsection{}

\begin{lem}\label{lem:final_step_for_tau_2}
  The map $ \mu$ is an isomorphism.
\end{lem}

By a combination of 
 Lemma \ref{lem:final_step_for_tau_2} with Lemma
  \ref{lem:calculate_final_action} and Lemma
  \ref{lem:obvious_transitive_action_on_tau_2}
we obtain:
\begin{kor}
    $\im({}^\pi d_2^{1,1})\cap \im({}^{\hat
    \pi}d_2^{1,1}) $ acts freely and transitively on
  $U$. \end{kor} 
Together with Lemma
  \ref{lem:only_twist_iso_classes} and Lemma \ref{lem:tau_2_equiv},
  this finishes the proof of Proposition \ref{freetrans}.
\hB

\subsubsection{}
We now prove Lemma \ref{lem:final_step_for_tau_2}.
  We have an isomorphism \begin{equation}\label{iso7r}\ker(i^*)/r^* H^2(B,\Z)\cong  \ker({}^r d_2^{1,1})\ .\end{equation}
  Let $i_b\colon
F_b\to F$ and $\hat i_b\colon \hat F_b\to \hat F$ be the inclusions of
the fibers over the basepoint. Using isomorphisms similar to (\ref{iso7r}) for 
$\pi$ and $\hat \pi$ we see that (\ref{iso7r})
induces an 
isomorphism
\begin{equation*}
  \frac{\ker(i^*)}{p^*(\ker(i_b^*))+\hat p^*(\ker(\hat i_b)^*)} \xrightarrow{\iso}
  \frac{\ker({}^r d_2^{1,1})}{\ker({}^\pi d^{1,1}_2)\oplus
    \ker({}^{\hat\pi}d^{1,1}_2)}\ .
\end{equation*}

To prove Lemma \ref{lem:final_step_for_tau_2}, we now consider the diagram
\begin{equation}\label{eq:big_diagram}
  \begin{array}{ccc}
     \frac{\ker({}^r d_2^{1,1})}{\ker({}^\pi d^{1,1}_2)\oplus
       \ker({}^{\hat\pi}d^{1,1}_2)} &\stackrel{{}^\pi d_2^{1,1}\circ \pr_1}{\rightarrow}&\im({}^\pi d_2^{1,1})\cap \im({}^{\hat
       \pi}d_2^{1,1}) \\
     \cong \uparrow&&\mu\downarrow\\
     \frac{\ker(i^*)}{p^*(\ker(i_b^*))+\hat p^*(\ker(\hat i_b^*))}
    &\xrightarrow[\cong]{\alpha }&\frac{H_{\tau_2}}{p^*H^2(F,\Z)+\hat p^* H^2(\hat F,\Z)}
  \end{array}\ .
\end{equation}
Here, the upper horizontal map is induced by the projection
$\pr_1:{}^rE_2^{1,1}\cong {}^\pi E_2^{1,1}\oplus {}^{\hat \pi} E_2^{1,1}\rightarrow {}^\pi E_2^{1,1}$ composed
with the differential. Since ${}^r d_2^{1,1}={}^\pi d_2^{1,1}+
{}^{\hat\pi} d_2^{1,1}$, on $\ker({}^r d_2^{1,1})$ this map indeed
maps to the intersection of the two images and is an
isomorphism. The map $\alpha$ is the isomorphism (\ref{u631}), where
we note that $p^*\ker(i_b^*)=p^*H^2(F,\Z)\cap \ker(i^*)$ and
$\hat p^*\ker(\hat i_b^*)=\hat p^*H^2(\hat F,\Z)\cap \ker(i^*)$.
Consequently, all maps in \eqref{eq:big_diagram} apart
from $\mu$ are isomorphisms, and Lemma \ref{lem:final_step_for_tau_2}
follows immediately from the following lemma.

\begin{lem}\label{lem:commuting_diagram}
  The diagram \eqref{eq:big_diagram} commutes (up to sign).
\end{lem}

Note that this reduces the proof of Proposition \ref{freetrans} to a
question about the cohomology of fiber bundles, more precisely about a
precise description of the secondary operation $\mu$.
% Although this
%might be well known, since we don't know of a reference we provide a
%proof of this statement.

\subsubsection{}
To study commutativity of the diagram, we can deal with one element in
the lower left corner
at a time. Since all the constructions
 are natural, it therefore suffices to study a universal situation and
prove the assertion there. To do this, we first have to identify such
a universal situation.

Given is a base space $B$ together with two $T^n$-bundles $F,\hat F$,
classified by their first Chern classes $c_1,\dots,c_n,\hat
c_1,\dots,\hat c_n\in H^2(B,\Z)$. Moreover, we have classes
$a_1,\dots,a_n,\hat a_1,\dots,\hat a_n\in H^1(B,\Z)$ such that
\begin{equation*}
  {}^r d_2^{1,1}(\sum_i (a_i\tensor y_i-\hat a_i\tensor \hat y_i) )= \sum_i(a_i\cup
  c_i-\hat a_i\cup \hat c_i)=0 \in H^3(B,\Z)\ .
\end{equation*}
 
\subsubsection{}

We now describe the universal case for this kind of data.
Let $t\colon K\to K(\Z,2)^{2n}\times K(\Z,1)^{2n}$ be the homotopy fiber of
the map $K(\Z,2)^{2n}\times K(\Z,1)^{2n}\to K(\Z,3)$ classified by
$\sum_i x_i\cup a_i-\hat x_i\cup \hat a_i$, where $x_i,\hat
x_i$, $a_i$, $\hat a_i$ are the canonical generators of the cohomology of the
factors of the source space. The universal $T^n$-bundles  classified by
$\bc_i$ and $\hat \bc_i$ pull back to bundles $\pi^u\colon F^u\to K$
and $\hat \pi^u\colon \hat F^u\to K$, and  $r^u\colon F^u\times_K \hat F^u\to K$ is the  fiber product.

A straightforward calculation with the Leray-Serre spectral sequence
of $t$ shows that we have isomorphisms
\begin{eqnarray}
   t^*:H^1(K(\Z,1)^{2n},\Z)&\stackrel{\cong}{\rightarrow}&  H^1(K,\Z) \nonumber\\ t^*:H^2(K(\Z,2)^{2n},\Z)\oplus
   H^2(K(\Z,1)^{2n},\Z)&\stackrel{\cong}{\rightarrow}&
    H^2(K,\Z) \label{eq:cohomol_of_K}\\
    \bar t^*:H^3(K(\Z,2)^{2n}\times K(\Z,1)^{2n},\Z)/(\sum_i(a_i\cup
    x_i-\hat a_i\cup\hat
    x_i)\Z)&\stackrel{\cong}{\rightarrow}&H^3(K,\Z) \nonumber
  \end{eqnarray}
%Note that these three groups are all free abelian.

Continuing with the Leray-Serre spectral sequences of $r^u$, $\pi^u$
and $\hat\pi^u$, we see that  for $k=2,3$
\begin{equation}
  \label{eq:cohom_of_univ_bunldes}
   H^k(F^u\times_K\hat F^u,\Z)=\ker(i^*)\ ,\quad
   H^k(F^u,\Z)=\ker(i_b^*)\ ,\quad H^k(\hat F^u,\Z)=\ker(\hat{i_b}^*),
 \end{equation}
 and that the projections 
 \begin{equation}\label{eq:inf_cyclic_grp}
\ker(i^*)/r^*H^2(K,\Z)\xrightarrow{\iso}
\ker({}^{r^u} d_2^{1,1}) \xrightarrow{\iso}
\ker({}^{r^u}d_2^{1,1})/(\ker({}^{\pi^u}d_2^{1,1})\oplus
\ker({}^{\hat\pi^u}d_2^{1,1}))\iso\Z
\end{equation}
all are isomorphisms.

\subsubsection{}

In the universal situation we have
$\ker(i^*)=H_{\tau_2}=H^2(F^u\times_K\hat F^u,\Z)$, and consequently
\begin{equation}
\im(\mu)\subset H_{\tau_2}/(p^*H^2(F^u,\Z)+\hat p^* H^2(\hat
F^u,\Z)). \label{eq:mu_with_right_trager}
\end{equation}
By naturality this implies that (\ref{eq:mu_with_right_trager}) holds in general.
This adds the missing detail in the construction of $\mu$
 of \eqref{eq:constr_of_mu}.

\subsubsection{}

We choose the  representative
$$g:=\sum_{i}( t^*a_i\tensor y_i-t^*\hat a_i\tensor \hat y_i)\in 
{}^{r^u}E_2^{1,1}$$ of the generator (see (\ref{eq:inf_cyclic_grp})) of $\ker({}^{r^u}d_2^{1,1})/(\ker({}^{\pi^u}d_2^{1,1})\oplus
\ker({}^{\hat\pi^u}d_2^{1,1}))$,
and let $\tilde g\in \ker(i^*)$ be a lift.
Let $[\tilde g]$ be the class represented by $\tilde g$
in the lower left corner of (\ref{eq:big_diagram}).
  Note that 
\begin{equation}
{}^{\pi^u}d_2^{1,1}\circ
\pr_1(g) =
\sum_i t^*a_i\cup t^*x_i = \sum_i
t^*\hat a_i\cup t^*\hat x_i.\label{eq:univ_d_pr}
\end{equation}
We must show that $$\mu(\sum_i t^*a_i\cup t^*c_i)=\alpha([\tilde g])\ .$$

%Consider now the map $v_0\colon B\to K(\Z,2)^{2n}\times K(\Z,1)^{2n}$ such that
%$v^*\bc_i=c_i$, $v^*\hat\bc_i=\hat c_i$, $v^*\bya_i=a_i$ and
%$v^*\hat\bya_i=\hat a_i$. It lifts to a map $v\colon B\to K$.  The element $\sum_i( t^*\bya_i\tensor y_i-t^*\hat\bya_i\tensor
%\hat y_i)\in \ker({}^{r^u}d^{1,1}_2)\subset {}^{r^u}E_2^{1,1}$ which pulls backto
%$\sum_i(a_i\tensor y_i-\hat u_i\tensor \hat y_i)$.

Since the relevant group is cyclic, it
suffices to show the equality of the leading parts in ${}^{r_u}
E_\infty^{1,1}=\ker({}^{r_u}d_2^{1,1})$, i.e. that
\begin{equation}
\mu(\sum_i t^*a_i\cup t^*x_i)^{1,1}=\alpha([\tilde g])^{1,1}=
\sum_i(t^*a_i\tensor y_i-t^*\hat a_i\tensor \hat
y_i)\label{eq:goal}\ .
\end{equation}

%modulo
%$(r^u)^* H^2(K,\Z)$ in $H^1(K)\tensor(H^1(F^u_b)\oplus H^1(\hat
%F^u_b))$, the ${}^{r^u}E^{1,1}_2$-term of the Leray-Serre spectral
%sequence for $H^2(F^u\times_K\hat F^u,\Z)$. 

\subsubsection{}

We have $(\pi^u)^*t^*x_i=0$ and $(\hat\pi^u)^*t^*\hat x_i=0$. Thus we can choose a homotopy $[0,1]\times F^u\to K(\Z,2)$
between the constant map and $(\pi^u)^*t^*x_i$, and a homotopy $[1,2]\times
\hat F^u \to K(\Z,2)$ between $(\hat\pi^u)^*t^*\hat x_i$ and the constant
map. Since the transgression of the Chern class of an $U(1)$-bundle
 is represented by a generator of the first
cohomology of the fiber 
(modulo the image of the first
cohomology of the total
space)  we can choose these homotopies in such a way that their
restrictions to $[0,1]\times F_b^u$ and $[1,2]\times\hat F_b^u$ (which
are necessarily constant at both ends of the intervals) are the suspensions of the
generators $y_i$ or $\hat y_i$, respectively.

We now take the product of the above homotopies with the corresponding
maps $(\pi^u)^*t^*a_i$ and $(\hat \pi^u)^*t^*\hat a_i$, respectively, and then the
product over the  index $i$ in order to
get 
homotopies
\begin{equation*}
  h_1\colon [0,1]\times F^u \to \prod_{i=1}^n (K(\Z,2)\times K(\Z,1))\
  ,\quad \hat
  h_1\colon [1,2]\times\hat F^u\to \prod_{i=1}^n (K(\Z,2)\times K(\Z,1)).
\end{equation*}
We finally compose with the map  $(K(\Z,2)\times K(\Z,1))^n\to
K(\Z,3)$ representing $\sum_i x_i\cup a_i$ to get the required
homotopies $H_1\colon [0,1]\times F^u\to K(\Z,3)$ and $\hat H_1\colon
[1,2]\times \hat F^u\to K(\Z,3)$ in the construction of $\mu(\sum t^* a_i\cup t^* x_i)$. 

\subsubsection{}
We consider the lifting problem
$$\begin{array}{ccc}
&&K\\
&\alpha_1\nearrow&t\downarrow\\
K(\Z,1)&\stackrel{w}{\rightarrow}&K(\Z,2)^{2n}\times
K(\Z,1)^{2n}\end{array}\ ,$$
where $w$ is the inclusion of the first $K(\Z,1)$-factor ---defined using the
base points of the remaining factors. Since $w^*(\sum_{i} a_i\cup
x_i-\hat a_i\cup \hat x_i)=0$, a lift $\alpha_1$ exists.

Since $\alpha_1^*t^*x_i$ and $\alpha_1^*t^*\hat x_i$ are constant,
it follows that the bundles $\alpha_1^*F^u$ and $\alpha_1^*\hat F^u$  are
 trivialized.  Since we
can choose the map $\sum x_i\cup a_i \colon (K(\Z,2)\times
K(\Z,1))^n \to K(\Z,3)$ to factorize over the  product of smash products
$\prod(K(\Z,2)\wedge
K(\Z,1))^n$ we see that  the pull-back of $H_1$ along $[0,1]\times \alpha_1^* F^u\rightarrow [0,1]\times F^u$ has a
factorization over the suspension $\Sigma (\alpha_1^*F^u_+)\rightarrow K(\Z,3)$.
It represents the suspension of the class $a\cup y_1\in
H^2(\alpha_1^* F^u,\Z)\cong H^2(K(\Z,1)\times T^n,\Z)$  since 
$t^*a_i$ pulls back to $\delta_{i1}a_1$, and the null homotopy for
$(\pi^u)^*t^*x_1$ gives the suspension of $y_1$ in the fiber by our
choices above. In the
same way the corresponding pull-back of $\hat H_1$ has a
factorization $\Sigma (\alpha_1^*\hat F^u_+)\rightarrow K(\Z,3)$ which
in this case is actually constant and thus represents the suspension of $0\in
H^2(\alpha_1^*\hat F^u,\Z)$.

\subsubsection{}

Let $A_1:\alpha_1^* (F^u\times_K\hat F^u)\rightarrow F^u\times_K\hat
F^u$ be the induced map over $\alpha_1$. We let $A_1^*$ denote the
induced map on the Leray-Serre spectral sequences.
Then we have by the above discussion that 
 $A_1^* \mu(\sum t^*a_i\cup x_i)^{1,1}=a\otimes y_1\in
 {}^{\alpha_1^*r_u}E_\infty^{1,1}$.
On the other hand $A_1^*(\sum_i(t^*a_i\tensor y_i-t^*\hat a_i\tensor \hat
y_i))=a\otimes y_1$.
This shows  (\ref{eq:goal}), since it is   an equality in a cyclic group .
\hB

 This finishes the proof of Lemma
\ref{lem:commuting_diagram} and therefore of Proposition \ref{freetrans}.\hB

\subsubsection{}

We consider the isomorphism class $[x_{n,univ}]\in \Triple(\bR_n)$ 
of the universal $T$-duality triple $((\bF_n,\cH_n),(\hat \bF_n,\hat
\cH_n),\bu_n)$ (see Theorem \ref{ex980}). 
Furthermore, we consider the functor
$B\mapsto P_n(B):=[B,\bR_n]$ classified by $\bR_n$.
The element $[x_{n,univ}]\in \Triple_n(\bR_n)$ induces a natural transformation of
functors
$\Psi_B:P_n(B)\rightarrow \Triple_n(B)$  which maps  $[f]\in
P_n(B)$ to $\Psi_B(f):=f^*[x_{n,univ}]$.

\begin{theorem}\label{main12}
The natural  transformation $\Psi:P_n\rightarrow \Triple_n$ is an isomorphism.
\end{theorem}

In other words, $\bR_n$ is the classifying space for $T$-duality
triples. The following is then a consequence of \ref{tz}.
\begin{kor}
For each space $B$ there is a natural action of the $T$-duality group
$G_n$ on the set of triples $\Triple_n(B)$.
\end{kor}

\subsubsection{}

\proof [of Theorem \ref{main12}]
Given a space $B$ we must show that $\Psi_B$ is a bijection of sets.
Therefore we show:
\begin{lem}\label{sur2}
$\Psi_B$ is surjective.
\end{lem}

\begin{lem}\label{inj2}
$\Psi_B$ is injective.
\end{lem}

\subsubsection{}\label{deded1}

The homotopy fiber of $(\bc,\hat \bc):\bR_n\rightarrow K(\Z^n,2)\times
K(\Z^n,2)$ is $K(\Z,3)$ (see \ref{rndefd}).
By obstruction theory the set of homotopy classes of lifts in the
diagram
$$\begin{array}{ccc}
&&\bR_n\\&f\nearrow&(\bc,\hat \bc)\downarrow\\
B&\stackrel{(c,\hat c)}{\rightarrow}&K(\Z^{n},2)\times K(\Z^n,2)\end{array}$$
is a torsor over $H^3(B,\Z)$.
So given two such lifts $f_0,f_1$ we have a well-defined difference 
element $\delta(f_1,f_0)\in H^3(B,\Z)$. 

Let $F:=c^*U^n$ and $\hat F:=\hat c^*U^n$ be the pull-backs of the
universal $T^n$-bundle over $K(\Z^n,2)$. Note further that
by construction $\bF_n=\bc^*U^n$ and $\hat \bF_n=\hat \bc^*U^n$.
Since $f$ lifts $(c,\hat c)$ we have natural identifications
$f^*\bF_n\cong F$ and $f^*\hat \bF_n\cong \hat F$.
We can thus consider $f^*[x_{n,univ}]\in \Triple_n^{(F,\hat F)}(B)$ in
a natural way. Recall the action of
$H^3(B,\Z)$ on $\Triple_n^{(F,\hat F)}(B)$ introduced in \ref{ert2}.
\begin{lem}\label{rret}
We have
$f_1^*[x_{n,univ}]=f_0^*[x_{n,univ}]+\delta(f_1,f_0)$.
\end{lem}
\proof
We can choose a model of $K(\Z,3)$ which is an abelian group.
Furthermore, we choose a model of the map
$(\bc,\hat \bc):\bR_n\rightarrow K(\Z^n,2)\times
K(\Z^n,2)$ which is a $K(\Z,3)$-principal bundle.
Let $a:\bR_n\times K(\Z,3)\rightarrow \bR_n$ be the right-action.
Let $\bz\in H^3(K(\Z,3),\Z)$ be a generator, and let $\cW$ be a twist in the
corresponding isomorphism class. We claim
that 
\begin{equation}\label{yxc1}a^* x_{n,univ}\cong \pr^*_1 x_{n,univ}+
  \pr_2^* \cW. \end{equation}
We have two maps
$\pr,m:\bF_n\times K(\Z,3)\rightarrow \bF_n$, namely the projection
and the map induced by the right-action
$a:\bR_n\times K(\Z,3)\rightarrow \bR_n$.
By the K{\"u}nneth formula and Lemma \ref{re1} we have
$H^3(\bF_n\times K(\Z,3),\Z)\cong \Z \pr^*\bh_n \oplus \Z \pr_2^*
\bz$.
Then $m^*\bh_n=d\pr^*\bh_n+b \pr_2^*\bz$ with  $d,b\in \Z$ to be determined.
By restriction to
$\bF_n\times \{1\}\subset \bF_n\times K(\Z,3)$ we see that $d=1$.
In order to calculate $b$ we restrict the torus bundle to a fiber of $(\bc,\hat \bc)$.
We identify this restriction with
$T^n\times K(\Z,3)$.
By Lemma \ref{g16} we know that the restriction of $\bh_n$ is
$1\times \bz$ (after choosing the appropriate sign of the  generator
$\bz$). It follows that the restriction of $m^*\bh_n$ is equal to
$1\times \bz\times 1+ 1\times 1\times b \bz\in H^3(T^n\times
K(\Z,3)\times K(\Z,3),\Z)$. 
Let $\mu :K(\Z,3)\times K(\Z,3) \rightarrow K(\Z,3)$
be the multiplication. Then we have $\mu^*\bz=\bz\times 1 + 1\times
\bz$. By a comparison of these two formulas  we conclude that $b=1$.

We now know that we have an isomorphism
$m^*\cH_n\iso\pr^*\cH_n\otimes \pr_2^*\cW$.
In a similar manner we have a map $\hat m:\hat \bF_n\times
K(\Z,3)\rightarrow \hat\bF_n$ and an isomorphism
$\hat m^*\hat \cH_n\cong \hat \pr^*\hat \cH_n\otimes \hat \pr_2^*\cW$.
Note that
$H^2(\bF_n\times_{\bR_n}\bF_n\times K(\Z,3)\times K(\Z,3),\Z)\cong \{0\}$.
Therefore the identification of the isomorphism classes  of the twists already implies (\ref{yxc1}).

We now consider maps $f_0,f_1:B\rightarrow \bR_n$.
Then we can write
$f_1=a\circ (f_0, g)$, where $g:B\rightarrow K(\Z,3)$
represents $\delta(f_1,f_0)\in H^3(B,\Z)$. It follows that
$f_1^*x_{n,univ}\cong (f_0,g)^*a^*x_{n,univ}\cong 
(f_0,g)^*(\pr^* x_{n,univ} + \pr_2^* \cW)\cong f_0^*x_{n,univ}+g^*\cW$.
Therefore
$f_1^*[x_{x,univ}]=f_0^*[x_{n,univ}]+\delta(f_1,f_0)$. \hB

\subsubsection{}
\proof [of Lemma \ref{sur2}].
It suffices to show the Lemma under the assumption that $B$ is connected.
Let $((F,\cH),(\hat F,\hat \cH),u)$ represent an isomorphism class
$x\in \Triple_n(B)$.
The bundles $F$ and $\hat F$ give rise to classifying maps
$c:B\rightarrow K(\Z^n,2)$ and $\hat c:B\rightarrow K(\Z^n,2)$. 
We thus have a map $(c,\hat c):B\rightarrow  K(\Z^n,2)\times
K(\Z^n,2)$.

Let $c_i,\hat
c_i$ be the components of $c,\hat c$.
It follows from the structure of $\cH$ that $$0={}^\pi d_2^{2,1}(\sum_{i=1}^n y_i\otimes \hat c_i)=\sum_{i=1}^n
c_i\cup \hat c_i\ .$$
This implies by the definition \ref{rt3} of $\bR_n$  the existence of
a lift $f$ in the diagram
$$\begin{array}{ccc}
&&\bR_n\\&f\nearrow&(\bc,\hat \bc)\downarrow\\
B&\stackrel{(c,\hat c)}{\rightarrow}&K(\Z^{n},2)\times K(\Z^n,2)\end{array}\ .$$
As already noticed in \ref{deded1},
the set of homotopy
classes of such lifts is a torsor over $H^3(B,\Z)$.
Now $x_f^\prime:= f^*x_{n,univ}=:((F^\prime,\cH^\prime),(\hat F^\prime,\hat
\cH^\prime),u^\prime)$. 
By construction we can assume that $F^\prime=F$ and $\hat
F^\prime=\hat F$ for all such lifts.
We therefore consider $[x],[x^\prime]\in \Triple_n^{(F,\hat F)}(B)$.
By Proposition \ref{freetrans} there exists a class $\alpha_f\in
H^3(B,\Z)$ such that $[x]=[x_f^\prime]+\alpha_f$. By Lemma \ref{rret}, if we 
use $f':=f+\alpha_f$, then $[x]=[x^\prime_{f'}]$.
\hB

\subsubsection{}

\proof[of Lemma \ref{inj2}]
Let $f_0,f_1:B\rightarrow \bR_n$ be given such that
$f_0^*x_{n,univ}$ and $f_1^* x_{n,univ}$ are isomorphic. We must show that $f_0$ and $f_1$
are homotopic. First we choose an isomorphism $f_0^*x_{n,univ}\cong f_1^* x_{n,univ}$. Let 
$\Psi:f^*_0\bF_n\rightarrow  f^*_1
\bF_n$ and
$\hat \Psi:f_0^*\hat \bF_n\rightarrow f_1^*\hat \bF_n$ denote the underlying isomorphisms of $T^n$-bundles. Their existence implies that the
compositions
$(\bc,\hat \bc)\circ f_0$ and $(\bc,\hat \bc)\circ f_1$ are
homotopic. We choose such a homotopy $H$
and consider the diagram
\begin{equation}\label{di610}\begin{array}{ccccc}
[0,1]\times f_0^*\bF_n\times_Bf_0^*\hat \bF_n&\stackrel{\Phi}{\rightarrow}&H^*(U^n\times U^n)&\rightarrow&U^n\times U^n\\
\downarrow&&\downarrow&&\downarrow\\
{}[0,1]\times B&\stackrel{\id}{\rightarrow} &[0,1]\times B&\stackrel{H}{\rightarrow}&K(\Z^n,2)\times K(\Z^n,2)\end{array}\ .\end{equation}
The bundle isomorphism $\Phi$ is uniquely determined up to homotopy by
the property that its restriction to $\{0\}\times B$ is the identity.
Let $\Phi_1:f_0^*\bF_n\times_Bf_0^*\hat \bF_n\rightarrow
f_1^*\bF_n\times_Bf_1^*\hat \bF_n $ be the restriction of $\Phi$ to
$\{1\}\times B$. Its homotopy class depends on the choice of $H$ in
the following way.

\subsubsection{}\label{conca32}

The homotopy $H$ can be concatenated with a map
$g:[0,1]\times B\rightarrow \Sigma (B_+)\rightarrow K(\Z^n,2)\times K(\Z^n,2)$. Homotopy classes of such maps are classified by 
$H^1(B,\Z^n)\times H^1(B,\Z^n)$. Since $T^n$ has the homotopy type of
$K(\Z^n,1)$, elements in this cohomology group can be represented by maps
$B\rightarrow T^n\times T^n$, and these maps act by automorphisms on
$f_0^*\bF_n\times_B f_0^*\hat \bF_n$. Let $H^\prime$ be obtained from $H$ by
concatenation with $g$, and let $\tilde g:B\rightarrow T^n\times T^n$
represent the corresponding homotopy class.
As a consequence of the discussion in \ref{con53}
 the homotopy class of the evaluation  $\Phi^\prime_1$ ($\Phi^\prime$ is defined by (\ref{di610}) for $H^\prime$ in the place of $H$) is obtained from $\Phi_1$ by composition with $g$. We see that we can arrange $H$ and the choice of $\Phi$ such that
$\Phi_1=(\Psi,\hat \Psi)$. For the remainder of the proof we adopt this choice.

\subsubsection{}

We choose a lift $\tilde H$ in the diagram
$$\begin{array}{ccc}
B& \stackrel{f_0}{\rightarrow}&\bR_n\\
b\mapsto (0,b)\downarrow&\tilde H\nearrow&(\bc,\hat \bc)\downarrow\\
I\times B&\stackrel{H}{\rightarrow}&K(\Z^{n},2)\times K(\Z^n,2)\end{array}
\ ,$$
and we let $\tilde f$ be the evaluation of $\tilde H$ at $\{1\}\times B$.
 Then $\tilde f$ and $f_1$ are two maps lifting
$(\bc,\hat \bc)\circ f_1$. Furthermore, $f_0$ is homotopic to $\tilde
f$. By Lemma \ref{homo32} the pull-backs $f_0^*x_{n,univ}$ and $\tilde f^*x_{univ}$ are isomorphic.
Note that the underlying $T^n$-bundles of $\tilde f^*x_{n,univ}$ are canonically isomorphic to $f_1^*\bF_n$ and $f_1^*\hat \bF_n$.
Actually, by the choice of $H$ above, the proof of \ref{homo32} gives us an isomorphism 
$f_0^*x_{n,univ}\cong \tilde f^* x_{n,univ}$ such that the underlying  bundle isomorphisms are $\Psi$ and $\hat \Psi$. The composition
$$f_1^*x_{n,univ}\cong f_0^*x_{n,univ}\cong \tilde f^* x_{n,univ}$$ is therefore an isomorphism over $(f_1^*\bF_n,f_1^*\hat \bF_n)$.
We conclude that
 $[f_1^*x_{n,univ}]=[\tilde f^* x_{n,univ}]$ holds in $\Triple_n^{(f_1^*\bF_n,f_1^*\hat \bF_n)}(B)$.

The maps $f_1$ and $\tilde f$ are lifts in the diagram
$$\begin{array}{ccc}
&&\bR_n\\&\nearrow&(\bc,\hat \bc)\downarrow\\
B&\stackrel{(\bc,\hat \bc)\circ f_1}{\rightarrow}&K(\Z^{n},2)\times
K(\Z^n,2)\end{array}\ .$$
Their difference as homotopy classes of lifts is measured by $\delta(f_1,\tilde f)\in H^3(B,\Z)$. But by Proposition \ref{freetrans} and the equality $[f_1^*x_{n,univ}]=[\tilde f^* x_{n,univ}]$ we have
 $\delta(f_1,\tilde f)=0$. 
Now $f_0$ is homotopic to $\tilde f$, and $\tilde f$ is homotopic to $f_1$ (even as a lift), so that $f_0$ is homotopic to $f_1$.
 \hB

\subsubsection{}

We fix classifying maps $(c,\hat c):B\rightarrow K(\Z^n,2)\times
K(\Z^n,2)$. We consider a lift $f$ in the diagram
$$\begin{array}{ccc}
&&\bR_n\\&f\nearrow&(\bc,\hat \bc)\downarrow\\
B&\stackrel{(c,\hat c)}{\rightarrow}&K(\Z^{n},2)\times
K(\Z^n,2)\end{array}\ .$$
Let $x:=f^*x_{n,univ}$ and write $x=((F,\cH),(\hat F,\hat \cH),u)$.
Let furthermore $\psi:F\rightarrow F$ and $\hat \psi:\hat F\rightarrow
\hat F$ be automorphisms of $T^n$-bundles. We will use the notation
introduced in \ref{eins7}
\begin{prop}\label{autact}
In $\Triple_n^{(F,\hat F)}(B)$ we have the following identity:
$$[x^{(\psi,\hat \psi)}]=[x]+\hat c\cup[\psi]+c\cup[\hat \psi]\ .$$
\end{prop}
\proof
We have $\Omega K(\Z^n,2)\cong K(\Z^n,1)$. As in \ref{conca32} we can concatenate the constant homotopy from $(c,\hat c)$ to $(c,\hat c)$ by a map corresponding to the class $([\psi],[\hat \psi])\in H^1(B,\Z^n)\times H^1(B,\Z^n)$. In this way we obtain a new homotopy $H:[0,1]\times B\rightarrow K(\Z^n,2)\times K(\Z^n,2)$  from $(c,\hat c)$ to $(c,\hat c)$.
We again consider the diagram
$$\begin{array}{ccccc}
[0,1]\times F\times_B \hat F&\stackrel{(\phi,\hat \phi)}{\rightarrow}&H^*(U^n\times U^n)&\rightarrow&U^n\times U^n\\
\downarrow&&\downarrow&&\downarrow\\
{}[0,1]\times B&\stackrel{\id}{\rightarrow} &[0,1]\times B&\stackrel{H}{\rightarrow}&K(\Z^n,2)\times K(\Z^n,2)\end{array}\ .$$

As explained in \ref{conca32} we can choose $\phi$ and $\hat \phi$ such that their restrictions $\phi_1$ and $\hat \phi_1$ to $\{1\}\times B$   coincide with $\psi$ and $\hat \psi$. We now choose a lift
$\tilde H$ in 
$$\begin{array}{ccc}
B&\stackrel{f}{\rightarrow}&\bR_n\\i_0 \downarrow&\tilde H\nearrow&(\bc,\hat \bc)\downarrow\\
{}[0,1]\times B&\stackrel{H}{\rightarrow}&K(\Z^{n},2)\times
K(\Z^n,2)\end{array}\ ,$$
where $i_0:B\cong \{0\}\times B\rightarrow [0,1]\times B$ is the
inclusion. Let $f_1:=\tilde H_{|\{1\}\times B}:B\rightarrow \bR_n$. Then we have 
by construction $f_1^*x_{n,univ}\cong x^{(\psi,\hat \psi)}$.
By Lemma \ref{rret} we have
$f_1^*[x_{n,univ}]=f^*[x_{n,univ}]+\delta(f_1,f_0)$ in $\Triple_n^{(F,\hat F)}(B)$.
Therefore it remains to show that
$\delta(f_1,f_0)=\hat c\cup[\psi]+c\cup[\hat \psi]$.
We want to apply the result of \ref{con53}. The homotopy $H$ gives rise to a map
$H^\prime:B\rightarrow L(K(\Z^n,2)\times K(\Z^n,2))$ which restricts to
$(c,\hat c)$ at $0\in S^1$. We must consider the composition
$$B\stackrel{H^\prime}{\rightarrow} L(K(\Z^n,2)\times K(\Z^n,2))\stackrel{Lq}{\rightarrow} LK(\Z,4)\stackrel{\sim }{\rightarrow}K(\Z,4)\times K(\Z,3)\stackrel{\pr_2}{\rightarrow} K(\Z^n,3)\  .$$
By (\ref{con54}) we know that
$\delta(f_0,f_1)$ is given by the cohomology class which classifies this map.
In order to calculate this class we consider the following general result.
 
\subsubsection{}

Fix integers $k,m\ge 1$.
A homotopy class $f\in [B,LK(\Z,k)]$ is characterized by
a pair of cohomology classes $(f_0,f_1)\in H^k(B,\Z)\times H^{k-1}(B,\Z)$.
Let $q:K(\Z,k)\times K(\Z,m)\rightarrow K(\Z,k)\wedge K(\Z,m) \rightarrow K(\Z,k+m)$ be the composition of the canonical projection and the cup product.
It induces a map
$Lq:LK(\Z,k)\times LK(\Z,m)\rightarrow LK(\Z,k+m)$.
We let $q^\prime:LK(\Z,k)\times LK(\Z,m)\rightarrow K(\Z,k+m-1)$
be the composition of $Lq$ with the projection 
$LK(\Z,k+m)\cong K(\Z,k+m)\times \Omega K(\Z,k+m)\rightarrow  \Omega
K(\Z,k+m)\iso K(\Z,k+m-1)$.
\begin{lem}\label{f84t}
We have $q_*^\prime(f,g)=f_0\cup g_1+f_1\cup g_0$.
\end{lem}
\proof
We consider (adjoint) representatives $f_1:\Sigma B_+\rightarrow K(\Z,k)$ and
$g_1:\Sigma B_+\rightarrow K(\Z,m)$ such that $f_1$ is constant on
$[0,1/2]\times B$ and $g_1$ is constant on $[1/2,1]\times B$.
Then we see that $q_*^\prime(f,g)=q_*(f_0,g_1)+q_*(f_1,g_0)$, where
$+$ indicates concatenation of loops. This follows  because of commutativity in the
diagram
\begin{equation*}
  \begin{CD}
    \Omega K(\Z,k-1)\wedge K(\Z,l) @>>> \Omega(K(\Z,k-1)\wedge
    K(\Z,l)) @>{\Omega \cup}>> \Omega K(\Z,k+l-1)\\
    @VV{\sim}V && @VV{\sim}V\\
    K(\Z,k)\wedge K(\Z,l) & @>{\cup}>> & K(\Z,k+l)
  \end{CD}
\end{equation*}
where the first arrow in the upper row comes from the inclusion
$K(\Z,l)\to\Omega K(\Z,l)$ as constant loops. The Lemma now
follows by taking the adjoint.
\hB

\subsubsection{}

Let us finish the proof of Proposition \ref{autact}.
Recall that $\bR_n$ is the homotopy fiber of the map
$q:K(\Z^n,2)\times K(\Z^n,2)\rightarrow K(\Z,4)$ classified by
$\sum_{i=1}^n x_i\cup\hat x_i$.
We consider the composition
$$L(K(\Z^n,2)\times K(\Z^n,2))\stackrel{Lq}{\rightarrow} LK(\Z,4)\cong
K(\Z,4)\times K(\Z,3)\stackrel{\pr_2}{\rightarrow} K(\Z,3)\ .$$
Furthermore note that
$$L(K(\Z^n,2)\times K(\Z^n,2))\cong K(\Z^n,2)\times K(\Z^n,1)\times K(\Z^n,2)\times K(\Z^n,1) .$$
We let $x_i,a_i,\hat x_i,\hat a_i$ denote the corresponding canonical
generators of the cohomology of the factors. Furthermore, we let
$\bz\in H^3(K(\Z,3),\Z)$ be the generator. We apply Lemma \ref{f84t}
to the summands of $q$ and obtain
$$Lq^*\circ \pr_2^*(\bz)=\sum_{i=1}^n (x_i\cup \hat a_i + \hat x_i\cup a_i)\ .$$
We now see that
$(H^\prime)^*x_i=c_i$, $(H^\prime)^*\hat x_i=\hat c_i$,
and  
using the discussion in \ref{con53}, $(H^\prime)^*a_i=[\psi]_i$, 
 $(H^\prime)^*\hat a_i=[\hat \psi]_i$. We get
$$(H^\prime)^*\sum_{i=1}^n (c_i\cup \hat a_i + \hat c_i\cup a_i)=c\cup [\hat
\psi]+\hat c\cup [\psi]\ ,$$
using the notation of \ref{eins7}.
\hB 

\subsubsection{}

We fix a pair $(F,h)$ over a space $B$. Then we consider 
$n$-dimensional $T$-duality triples $x=((F,\cH),(\hat F,\hat \cH),u)$
which extend $(F,h)$ (see Definition \ref{exexdef}).

\begin{ddd}\label{isoextensiondef}
An isomorphism of extensions $$x:=((F,\cH),(\hat F,\hat
\cH),u)\ , \quad x^\prime:=((F,\cH^\prime),(\hat F^\prime,\hat
\cH^\prime),u^\prime)$$ is an isomorphism of $T$-duality triples
$x\cong x^\prime$ (see \ref{isodd}) which induces the identity on
$F$. The set of such isomorphism classes is denoted $\Ext(F,h)$.
\end{ddd}

\subsubsection{}\label{provecla}
\proof[of Theorem \ref{class76}]
The first assertion of Theorem \ref{class76} is exactly Theorem \ref{exex}.
We continue with the existence statement of (2).
Let $x:=((F,\cH),(\hat \cF,\hat \cH),u)$ represent a class in
$\Ext(F,h)$. Let $B\in \Mat(n,n,\Z)$ be antisymmetric
and set $\hat c^\prime_i:=\hat c_i(x)+\sum_{j=1}^n B_{i,j} c_j(x)$.
We compute
$\sum_{i=1}^n c_i(x)\cup \hat c_i^\prime=\sum_{i=1}^n c_i(x)\cup \hat
c_i(x)+\sum_{j,i=1}^n B_{i,j}c_i(x)\cup  c_j(x)=0$.
This implies the existence of a lift $f^\prime$ in
$$\begin{array}{ccc}
&&\bR_n\\&f^\prime\nearrow&(\bc,\hat \bc)\downarrow\\
B&\stackrel{(c,\hat c^\prime)}{\rightarrow}&K(\Z^{n},2)\times
K(\Z^n,2)\end{array}\ .$$
Then $\{x^\prime\}:=\{(f^\prime)^*x_{n,univ}\}\in \Ext(F,h)$ has the
required properties. The remaining part of (2) was shown in Subsection
\ref{after43}.

We now show (3). We first consider the set of classes of triples
$[x^\prime]\in \Triple_n^{(F,\hat F)}(B)$ which extend $(F,h)$. It follows
from Proposition \ref{freetrans} that $\ker(\pi^*)$ acts freely and
transitively on
this set. It follows that $\ker(\pi^*)$ acts transitively on any subset
of $\Ext(F,\cH)$ with fixed $c(\hat F)$. But note that the equivalence relation in
$\Triple_n^{(F,\hat F)}(B)$ is stronger than  in $\Ext(F,h)$ since
isomorphisms must induce the identity on $\hat F$. In $\Ext(F,h)$  we admit
non-trivial bundle automorphism of $\hat F$. In view of Proposition
\ref{autact} we see that the quotient $\ker(\pi^*)/\im(C)$ acts
freely on the set of isomorphisms classes of extensions of $(F,h)$
with prescribed isomorphism class of $\hat F$.
\hB

\begin{appendix}
  
\section{Twists, spectral sequences and other conventions}

\subsection{Twists}
\label{sec:twists}

{}\label{twistsexplain}

We start with a description of twists. The literature contains various models
for twists. Therefore we first describe the
common core of these models and in particular the properties which we use in
the present paper. Then we exhibit two of these models explicitly. 

In each case one considers a pre-sheaf of monoidal groupoids $B\mapsto
T(B)$ on the category of spaces. The objects of $T(B)$ are called
twists.  The unit of the monoidal structure is called the trivial
twist and denoted by $0$. 
If $f:A\rightarrow B$ is a map of spaces, then there is a monoidal functor
$f^*:T(B) \rightarrow T(A)$. Functoriality is implemented by natural
transformations. We refer to \cite[Section 3.1]{bunkeschickt} for more details.
The following three requirements provide the coupling to topology.
\begin{enumerate}
\item
We require that there is a natural monoidal transformation
$T(B)\rightarrow H^3(B,\Z)$, $\cH\mapsto [\cH]$, (the group
$H^3(B,\Z)$ is considered as a monoidal  category which has only identity morphisms) which classifies the isomorphism
classes of $T(B)$ for each $B$.
\item
If $\cH,\cH^\prime\in T(B)$ are equivalent objects, then we require that
$\Hom_{T(B)}(\cH,\cH^\prime)$ is a
$H^2(B,\Z)$-torsor such that  the composition with a
fixed morphism  form $\cH'$ to $\cH''$ or from
$\cH''$ to $\cH$, respectively, gives isomorphisms of torsors.
Furthermore, we require that the torsor structure is
compatible with the pull-back. Note that we have natural bijections
$\Hom(\cH,\cH)\stackrel{\sim}{\rightarrow} H^2(B,\Z)$ which map compositions
to sums. In order to simplify the notation we will (and have in the text)
frequently identify automorphisms of twists and cohomology classes.
The map
$\Hom_{T(B)}(\cH_0,\cH^\prime_0)\times
\Hom_{T(B)}(\cH_1,\cH^\prime_1)
\rightarrow \Hom_{T(B)}(\cH_0\otimes \cH_1,\cH^\prime_0\otimes
\cH^\prime_1)$,
$(v,w)\mapsto v\otimes w$ is bilinear.
\end{enumerate}
We now discuss two explicite realizations.
\begin{enumerate}
\item
Let $K$ be the algebra of compact operators on a separable complex
Hilbert space. The group of automorphisms of $K$ is the projective
unitary group $PU$. As a topological group we have $PU=U/U(1)$,
where the topology of the unitary group of the Hilbert space is the
strong topology. We define $T(B)$ to be the category of locally trivial
bundles with fiber $K$ such that the transition functions are continuous
functions with values in $PU$. We further define $\Hom_{T(B)}(\cH,\cH^\prime)$ as the
set of homotopy classes of algebra bundle isomorphisms. The monoidal
structure is given by the fiberwise (completed) tensor product of bundles. In
order to see that this preserves the category we use the fact that
$K\otimes K\cong K$. Isomorphism classes of objects in $T(B)$  are classified by
homotopy classes $[B,BPU]$. Since the classifying space $BPU$  has the homotopy type
$K(\Z,3)$ we see that isomorphism classes in $T(B)$ are in one-to one
correspondence with $H^3(B,\Z)$.
The group of automorphisms of a bundle $\cH\in T(B)$ can be identified
with the group of 
homotopy classes of maps $[B,PU]$ and is therefore in one-to-one
correspondence to $H^2(B,\Z)$ (this identification is not quite automatic and
uses the fact that $PU$ is homotopy commutative, to pass from the sections of
the inner automorphism bundle of the non-trivial bundle to the ones of the
trivial bundle).

This model of twists is very suitable for a quick definition of  twisted $K$-theory. Assume
that $B$ is locally compact. 
Given a twist $\cH\in T(B)$ we consider the $C^*$-algebra $C_0(B,\cH)$
of continuous sections of $\cH$ vanishing at infinity. Then we can define
$K(B,\cH):=K(C_0(B,\cH))$, where the right-hand side is the $K$-theory of
the $C^*$-algebra $C_0(B,\cH)$. 
\item
In our second model we fix an $h$-space model of $K(\Z,3)$. Note that
there exists in fact models which are topological abelian groups.
A twist $\cH\in T(B)$ is a map $\cH:B\rightarrow K(\Z,3)$. The  monoidal
structure is implemented by the $h$-space structure of $K(\Z,3)$.
The set $\Hom_{T(B)}(\cH,\cH^\prime)$ in this model is the set of
homotopy classes of homotopies from $\cH$ to $\cH^\prime$.
It is obvious from the definition that isomorphism classes in $T(B)$
are classified by $H^3(B,\Z)$. Furthermore, since $\Omega K(\Z,3)\cong
K(\Z,2)$ we see that  $\Hom_{T(B)}(\cH,\cH^\prime)$ is empty or a
torsor over $H^2(B,\Z)$.
\end{enumerate}

\subsection{Spectral sequences}\label{lerayserreexplain}
 
 The cohomology of the total space $F$ of a fiber bundle $\pi:F\rightarrow B$
has a decreasing filtration $$0\subset \cF^nH^n(F,\Z)\subset
\cF^{n-1}H^n(F,\Z)\subset \dots \subset \cF^1H^n(F,\Z)\subset \cF^0
H^n(F,\Z)=H^n(F,\Z)$$ such that
$\cF^nH^n(F,\Z)=\pi^*(H^n(B,\Z))$.  By definition, $x\in \cF^kH^n(F,\Z)$ if for any $(k-1)$-dimensional
$CW$-complex $X$ and map $\phi:X\rightarrow B$ the condition
$\Phi^*(x)=0$ is satisfied. Here $\Phi:\phi^*F\rightarrow F$ is the
induced map. The associated graded group is
calculated by the Leray-Serre spectral sequence $({}^\pi E_r^{s,t},
{}^\pi d_r^{s,t})$. If $x\in \cF^k H^n(F,Z)$, then we let
$x^{k,n-k}\in {}^\pi E^{k,n-k}_\infty$ denote  the leading part.

Let us assume that $B$ is connected and choose a base point $b\in
B$. Let $F_b:=\pi^{-1}(b)$ be the fiber over $b$. We further assume that $\pi_1(B,b)$ acts trivially on
$H^*(F_b,\Z)$.
In the present paper this assumption is always satisfied. 
Then we have isomorphisms 
${}^{\pi} E_{2}^{p,q}\cong H^p(B,H^q(F_b,\Z))$.
For various calculations we will employ naturality and
multiplicativity of the spectral sequence.

Let us now fix some notation in the case that $\pi:F\rightarrow B$ is
a $T^n$-principal bundle. We fix a generator of $
H^1(U(1),\Z)$. Since $T^n:=U(1)^n$ this
provides a natural set of generators of $H^1(T^n,\Z)$.
If we fix a base point in $F_b$, then using the right action we obtain a homeomorphism
$F_b\cong T^n$. We let $y_1,\dots,y_n\in H^1(F_b,\Z)$ denote the
natural generators obtained in this way. They do not depend on the
choice of the base point in $F_b$. 
If we consider a bundle $\hat \pi:\hat F\rightarrow B$, then we will
write
$\hat y_1,\dots , \hat y_n$ for the corresponding set of generators.
Finally, the notation $y_1,\dots,y_n,\hat y_1,\dots,y_n$ for the
generators
of $H^1(F_b\times \hat F_b,\Z)$ is self-explaining.

Note that $H^*(F_b,\Z)$ is a free $\Z$-module. Therefore we have an
isomorphism
$${}^{\pi} E_{2}^{p,q}\cong H^p(B,H^q(F_b,\Z))\cong H^q(F_b,\Z)\otimes
H^p(B,\Z)\ .$$ Let $c_1,\dots,c_n\in H^2(B,\Z)$ be the Chern classes of
the $T^n$-bundle $\pi:F\rightarrow B$. Then we know that
${}^\pi d_2^{0,1}(y_i)=c_i$. By multiplicativity this leads to a
complete calculation of ${}^\pi d_2$.

\subsubsection{}

%\Kommentar{This is an alternative to a subsection by Uli ---now
%  commented out---, stressing that
%  we only use the axioms of twists.}

Let $T\xrightarrow{i} F\xrightarrow{\pi} B$ be a fibration (of pointed
spaces) with base $B$ and fiber $T$. Choose
$\delta\in H^3(B,\Z)$ such that $\pi^*\delta=0$.
We choose a twist $\cH$ over $F$ such that $[\cH]=\delta$.
Furthermore, we choose a trivialization  $w\colon \pi^*\cH\stackrel{\sim}{\rightarrow} 0$.
Since $\pi\circ i$ is the constant map to the base point we have a natural isomorphism $0\stackrel{can}{\rightarrow} \pi^*\cH$. We obtain the automorphism
$w_b:0\stackrel{can}{\rightarrow} \pi^*\cH \stackrel{i^*w}{\rightarrow}0$ of the trivial twist. We will consider $w_b\in H^2(T,\Z)$.

Note that the set of  trivializations $w$ is a torsor over $H^2(F,\Z)$.
It follows that the class $[w_b]\in H^2(T,\Z)/i^*H^2(F,\Z)$ is independent of the choice of $w$.
If we replace $\cH$ by an isomorphic twist $\cH'$ with isomorphism
$u\colon \cH\to \cH'$, and we let $w'$ be obtained
from $w$ by conjugation with $\pi^*u$, then we have $w_b=w_b^\prime$.
Therefore the class $[w_b]\in H^2(T,\Z)/i^*H^2(F,\Z)$ is well-defined independent of the choice of $\cH$ in the class given by $\delta$. 
Note that we can consider ${}^\pi E_3^{0,2}/\ker({}^{\pi}d_3^{0,2})\subset H^2(T,\Z)/i^*H^2(F,\Z)$.

\begin{lem}\label{lem:i891}
We have $[w_b]\in {}^\pi E_3^{0,2}/\ker({}^{\pi}d_3^{0,2})$
and ${}^{\pi}d_3^{0,2}([w_b]))=[\delta]\in H^3(B,\Z)/\im({}^\pi d_2^{1,1})$.
\end{lem}
\proof
We choose a model for
$K(\Z,3)$ and a (pointed) classifying map
$f\colon B\to K(\Z,3)$ such that $\delta = f^*\bz$, where $\bz\in
H^3(K(\Z,3),\Z)$ is the canonical generator. We further choose a twist $\cK$ on $K(\Z,3)$ in the class $\bz$. Then we can assume that $\cH=f^*\cK$.

We
choose a (pointed) homotopy $H\colon [0,1]\times F \to K(\Z,3)$ between $f\circ
p$ and the constant map which exists because $\pi^*\delta=0$. The
homotopy class of such homotopies is well
defined upto the action of $H^2(F,\Z)\cong H^3(\Sigma F,\Z)$, where $\Sigma
F$ is the reduced suspension of $F$. Elements of $H^3(\Sigma  F,\Z)$
are represented by maps $\Sigma F\to K(\Z,3)$, and the action
mentioned above is defined by concatenation with the homotopy given by
the composition
$[0,1]\times F\onto \Sigma F\to K(\Z,3)$.

Note that by the axioms of twists each such map $H$ induces a well-defined isomorphism of twists 
\begin{equation}
\pi^*\cH\stackrel{\sim}{\rightarrow} H_0^*\cK\stackrel{\sim}{\rightarrow} H_1^*\cK\stackrel{\sim}{\rightarrow} 0,
\label{eq:iso_1}
\end{equation}
where
we use canonical identifications, and in particular that $H_1$ is the
constant map. The axioms also imply that the action of $H^2(F,\Z)$ on
such isomorphisms corresponds exactly to the action of $H^2(F,\Z)$ on
the homotopies $H$ as constructed above. Therefore, we can choose $H$
such that the isomorphism \eqref{eq:iso_1} is exactly the isomorphism
$w$ chosen above.

If we restrict $H$ to the fiber $T$, then we obtain a map $h\colon [0,1]\times T \to K(\Z,3)$. Since
$H_0=f\circ p$ and $p\circ i$ is the constant map to the base point, it follows that
$h_0$ is the constant map. Consequently, we have a factorization of
$h$ as $[0,1]\times T\onto \Sigma T\xrightarrow{h'} K(\Z,3)$. The
axioms of twists give us an isomorphism $h_0^*\cK\stackrel{\sim}{\rightarrow} h_1^*\cK$ which
is by naturality the pullback $i^*w$. We now  see that
$w_b\in H^2(T,\Z)$ corresponds to $(h')^*\bz \in H^3(\Sigma
T,\Z)$ under the
suspension isomorphism.

A different interpretation uses adjunction to translate $h'$ to a map
$g\colon T\to \Omega K(\Z,3)=K(\Z,2)$. Then $w_b=g^*\bu$, where $\bu$ is
the canonical generator of $H^2(K(\Z,2),\Z)$. We now consider the
Leray-Serre spectral sequence of the fibration $$\Omega K(\Z,3)\to
PK(\Z,3)\xrightarrow{p} K(Z,3)\ ,$$ where $PK(\Z,3)$ is the (contractible) space of
pointed paths in $K(\Z,3)$. Observe that ${}^p d_2^{0,2}(\bu)=0$ and ${}^p d_3^{0,2}(\bu)=\bz$. 
By naturality, $[w_b]\in {}^\pi E_3^{0,2}/\ker({}^{\pi}d_3^{0,2})$ and 
${}^{\pi} d_3^{0,2}(w_b) = [\delta]\in H^3(B,\Z)/\im({}^\pi d_2^{1,1})$.  \hB

\subsection{Classification of lifts}\label{con53}

We consider an integer $k\ge 1$ and an abelian group $G$. 
We form the loop space  $LK(G,k):=\Map(S^1,K(G,k))$. Let $h:K(G,k)\times K(G,k)\rightarrow K(G,k)$ be the $h$-space structure. Then the map
$\phi:K(G,k)\times \Omega K(G,k)\rightarrow LK(G,k)$ given by $\phi(x,l)(t)=h(x,l(t))$, $t\in S^1$, is a homotopy equivalence. We choose a homotopy inverse $\psi$ of $\phi$. Note that $\Omega K(G,k)\cong K(G,k-1)$.

We consider a space $B$ and a map
$c:B\rightarrow K(G,k)$. Furthermore  we  let $p:X\rightarrow B$ be the homotopy fiber of $c$. Then the homotopy fiber of $p$ is a $K(G,k-1)$. By obstruction theory
the set of lifts $\tilde f$ in the diagram 
$$\begin{array}{ccc}
&&X\\&\tilde f\nearrow&p\downarrow\\
Y&\stackrel{f}{\rightarrow}&B\end{array} $$
is empty or a torsor over $H^{k-1}(Y,G)$.
Given two lifts $\tilde f_0,\tilde f_1$ we have a difference element $\delta(\tilde f_0,\tilde f_1)\in H^{k-1}(Y,G)$.

In fact, we can choose a model for $K(G,k-1)$ which is an abelian group. Furthermore, we can choose a model for $p:X\rightarrow B$ which is a $K(G,k-1)$-principal bundle. Then there exists a unique map
$g:Y\rightarrow K(G,k-1)$ such that
$\tilde f_1=\tilde f_0 g$ using the right-action of $K(G,k-1)$ on $X$.
In this case the homotopy class of $g$ is classified by
$\delta(\tilde f_0,\tilde f_1)$.

We consider now a map
$H:S^1\times Y\rightarrow B$ such that $H_{|\{0\}\times Y}=f$
(we parameterize $S^1$ by $[0,1]$ with endpoints identified).
Furthermore we consider a lift $\tilde H$ in the diagram
$$\begin{array}{ccc}
Y&\stackrel{\tilde f}{\rightarrow}&X\\
i_0\downarrow&\tilde H\nearrow&p\downarrow\\
{}[0,1]\times Y&\stackrel{H\circ q}{\rightarrow}&B\end{array} ,$$
where $q:[0,1]\times Y\rightarrow S^1\times Y$ is the quotient map, and $i_0(y):=(0,y)$.
We define  $\tilde f^\prime:=\tilde H_{|\{1\}\times Y}:Y\rightarrow X$.

Let $H^\prime:Y\rightarrow LB$ be the adjoint of $H$.
By composition with $Lc:LB\rightarrow LK(G,k)$, the homotopy inverse $\psi$, and the projection to the second component we obtain a map
$$\pr_2\circ \psi\circ \circ Lc\circ H^\prime:Y\rightarrow K(G,k-1)\ .$$
The homotopy class of this map is classified by a cohomology class
$d(H)\in H^{k-1}(Y,G)$.
Several times we need the following formula:
\begin{equation}\label{con54}
\delta(\tilde f^\prime,\tilde f)=d(H)
\end{equation}

\end{appendix}

\end{document}